\documentclass[oneside, reqno] {amsart}

\usepackage{xypic}
\input xy
\xyoption{all}
\usepackage{epsfig}
\usepackage{amsthm}
\usepackage{amssymb}
\usepackage{amsmath}
\usepackage{amscd}
\usepackage{color}
\usepackage[T1]{fontenc}
\usepackage{upgreek}
\usepackage[font=scriptsize]{caption}
\usepackage{wrapfig}

\usepackage{graphicx}
\usepackage{subfigure}

\newcommand{\bg}{\begin{equation}}
\newcommand{\ed}{\end{equation}}
\newcommand{\bga}{\begin{eqnarray}}
\newcommand{\eda}{\end{eqnarray}}
\newcommand{\pf}{\textbf{Proof:\ }}

\def\cbdu{\par{\raggedleft$\Box$\par}}

\newtheorem {Theorem}  {Theorem}

\numberwithin{Theorem}{section}

\newtheorem {Lemma}[Theorem]  {Lemma}

\theoremstyle{definition}
\newtheorem{Definition}[Theorem]{Definition}
\theoremstyle{remark}
\newtheorem{Remark}[Theorem]{\bf Remark}

\expandafter\chardef\csname pre amssym.def
at\endcsname=\the\catcode`\@ \catcode`\@=11
\def\undefine#1{\let#1\undefined}
\def\newsymbol#1#2#3#4#5{\let\next@\relax
 \ifnum#2=\@ne\let\next@\msafam@\else
 \ifnum#2=\tw@\let\next@\msbfam@\fi\fi
 \mathchardef#1="#3\next@#4#5}
\def\mathhexbox@#1#2#3{\relax
 \ifmmode\mathpalette{}{\m@th\mathchar"#1#2#3}%
 \else\leavevmode\hbox{$\m@th\mathchar"#1#2#3$}\fi}
\def\hexnumber@#1{\ifcase#1 0\or 1\or 2\or 3\or 4\or 5\or 6\or 7\or 8\or
 9\or A\or B\or C\or D\or E\or F\fi}

\font\teneufm=eufm10 \font\seveneufm=eufm7 \font\fiveeufm=eufm5
\newfam\eufmfam
\textfont\eufmfam=\teneufm \scriptfont\eufmfam=\seveneufm
\scriptscriptfont\eufmfam=\fiveeufm

\catcode`\@=\csname pre amssym.def at\endcsname

\newcounter{remark}
\def  \12  {{\frac{1}{2}}}

\def\build#1_#2^#3{\mathrel{\mathop{\kern 0pt#1}\limits_{#2}^{#3}}}

\numberwithin{equation}{section}

\begin{document}

\title[2D EMHD near steady state]{Global existence of 2D electron MHD near a steady state}


\author [Mimi Dai]{Mimi Dai}

\address{Department of Mathematics, Statistics and Computer Science, University of Illinois at Chicago, Chicago, IL 60607, USA}
\email{mdai@uic.edu}



\thanks{The author is partially supported by the NSF grant DMS--2009422 and the AMS Centennial Fellowship.}

\begin{abstract}
We study the electron magnetohydrodynamics (MHD) in two dimensional geometry, which has a rich family of steady states. In an anisotropic resistivity context, we show global in time existence of small smooth solution near a shear type steady state. Convergence rate of the solution to the steady state is also obtained.

\bigskip


\hspace{0.02cm}CLASSIFICATION CODE: 35Q35, 76B03, 76D09, 76E25, 76W05.
\end{abstract}

\maketitle

\section{Introduction}

\medskip

\subsection{Overview}

The electron magnetohydrodynamics (MHD) model
\begin{equation}\label{emhd}
B_t+ \nabla\times ((\nabla\times B)\times B)= \mu\Delta B
\end{equation}
describes the motion of magnetized electrons within a background passive and neutralizing ion flow, see \cite{Bis1}. The vector $B$ represents the magnetic field and $J=\nabla\times B$ the current density, and the parameter $\mu$ measures the strength of the resistivity. The nonlinear term in (\ref{emhd}) captures the Hall effect which is believed to be responsible for the rapid magnetic reconnection phenomena in plasmas \cite{BDS}. We note that $\nabla\cdot B(t)=0$ remains for all the time if $\nabla\cdot B(0)=0$ initially.

The electron MHD has been studied extensively by physicists, mostly through numerical simulations, in order to understand the geometry configurations in the reconnection region of the plasma. For instance, we refer the reader to \cite{RDD, SDS} and references therein. In particular, the steady states of (\ref{emhd}) play a vital role in the investigations \cite{DKM, GH, WHL}. Indeed, equation (\ref{emhd}) has a rich class of equilibria, in both of the 3D and 2D settings. Due to the fact that the 2D geometry is more manageable in numerical study, many contributions \cite{CSZ, KC, WH} in the physics community have been emerged for 2D electron MHD. 
From another point of view, we observe (\ref{emhd}) is quasi-linear and supercritical (see \cite{Dai-W}) which presents obstacles for mathematical analysis of the model. Therefore we focus on the electron MHD in 2D geometry in this article.

We consider the magnetic field in the form
\begin{equation}\label{two-half}
B(x,y,t)=\nabla\times (a\vec e_z)+b\vec e_z \ \ \mbox{with} \ \ \vec e_z=(0, 0,1), \ \ (x,y)\in\mathbb R^2,
\end{equation}
where 
\[a=a(x,y,t), \ \ b=b(x,y,t)\]
are scalar-valued functions. 
It is easy to see $\nabla\cdot B=0$ and $B=(a_y, -a_x, b)$. Recasting (\ref{emhd}) into the equations of $a$ and $b$ (with more general resistivity)
\begin{equation}\label{emhd-ani}
\begin{split}
a_t+ (a_yb_x-a_xb_y)=&\ \mu_1\Delta a, \\
b_t-(a_y\Delta a_x-a_x\Delta a_y)=&\ \mu_2\Delta b
\end{split}
\end{equation}
where $\mu_1, \mu_2\geq 0$ denote the horizontal and vertical resistivity respectively. When $\mu_1\neq \mu_2$, the system is said to be anisotropic. Fluid equations with anisotropic dissipation have been studied in many works, including the ones \cite{CDGG, Pai} for the Navier-Stokes equation, and \cite{ CRW, LXZ, LZ, RWXZ} for MHD systems. 
The basic energy law of (\ref{emhd-ani}) 
\begin{equation}\notag
\frac12\frac{d}{dt}\int \left(a_x^2+a_y^2+b^2\right)\, dxdy+ \int \mu_1\left(a_{xx}^2+2a_{xy}^2+a_{yy}^2\right)+\mu_2\left(b_x^2+b_y^2\right)\, dxdy=0
\end{equation}
indicates the a priori energy estimates 
\begin{equation}\notag
\begin{split}
& a\in L^\infty(0,T; H^1)\cap L^2(0,T; H^2), \\ 
& b\in L^\infty(0,T; L^2)\cap L^2(0,T; H^1).
\end{split}
\end{equation}
According to the natural scaling of (\ref{emhd-ani}), the critical Sobolev space for $(a,b)$ is $\dot H^2\times \dot H^1$.
Thus system (\ref{emhd-ani}) is energy supercritical even it is in 2D. The question of global existence of regular solution for (\ref{emhd-ani}) in general setting remains open. In \cite{Dai-W} the authors established a regularity criterion which only relies on the low modes of the solution. 

We observe that (\ref{emhd-ani}) has many steady states $(a^*, b^*)$, including two interesting categories: shear profiles and radial profiles. A steady state of (\ref{emhd-ani}) with $\mu_1=\mu_2=0$ is referred as Hall equilibrium. A few examples of Hall equilibria are 
\begin{itemize}
\item [(i)] shear type: $a^*=a^*(x)$ and $b^*=b^*(x)$;
\item [(ii)] shear type: $a^*=a^*(y)$ and $b^*=b^*(y)$;
\item [(iii)] radial type: $a^*=a^*(r)$ and $b^*=b^*(r)$ with $r=(x^2+y^2)^{\frac12}$.
\end{itemize}
A special case of (ii) with $a^*=g(y)$ and $b^*=0$ for a class of functions $g$, and a special case of (iii) with $b^*=0$ were studied in \cite{JO}, where the authors showed ill-posedness of the system around such steady states. The authors of \cite{GH} investigated a special case of (i) numerically and discussed essential conditions for instability. 
On the other hand, numerical study in \cite{WH} suggests that stationary structures actually tend to develop in time. The aim of this paper is to verify the numerical discovery rigorously in the context of a current sheet steady state: $(a^*, b^*)=(y, 0)$, which corresponds to $B^*=(1,0,0)$. In particular, $a^*=y$ does not satisfy the conditions for the function $g$ studied in \cite{JO}.

Denote the perturbation $\psi=a-y$. The perturbation pair $(\psi, b)$ satisfies the system
\begin{equation}\label{emhd-pert}
\begin{split}
\psi_t+b_x-b_y\psi_x+b_x\psi_y=&\ \mu_1\Delta \psi,\\
b_t-\Delta \psi_x-\psi_y\Delta \psi_x+\psi_x\Delta \psi_y=&\ \mu_2\Delta b.
\end{split}
\end{equation}
The linearized system is
\begin{equation}\label{emhd-linear}
\begin{split}
\psi_t+b_x=&\ \mu_1\Delta \psi,\\
b_t-\Delta \psi_x=&\ \mu_2\Delta b.
\end{split}
\end{equation}
With $\mu_1=0$ and $\mu_2=1$, it follows from the two equations of (\ref{emhd-linear}) that
\begin{equation}\label{wave}
b_{tt}-\Delta b_t+\Delta b_{xx}=0.
\end{equation}
The Fourier transform of (\ref{wave}) is 
\begin{equation}\label{wave-f}
\widehat b_{tt}+|\xi|^2\widehat b_t+\xi_1^2|\xi|^2\widehat b=0
\end{equation}
which is a second order ODE in time variable. The characteristic equation of (\ref{wave-f}) has solutions 
\begin{equation}\notag
\lambda_{\pm}=-\frac12 |\xi|^2\left(1\pm \sqrt{1-\frac{4\xi_1^2}{|\xi|^2} }\right).
\end{equation}
Thus the solutions of (\ref{wave-f}) are given by
\begin{equation}\notag
\widehat b(\xi,t)=c_1(\xi)e^{t\lambda_{+}}+c_2(\xi)e^{t\lambda_{-}}
\end{equation}
with
\begin{equation}\notag
c_1(\xi)= \frac{\lambda_{-}(\xi)\widehat b_0(\xi)-\widehat b_1}{|\xi|^2 \sqrt{1-\frac{4\xi_1^2}{|\xi|^2} }}, \ \ \
c_2(\xi)=- \frac{\lambda_{+}(\xi)\widehat b_0(\xi)-\widehat b_1}{|\xi|^2 \sqrt{1-\frac{4\xi_1^2}{|\xi|^2} }},
\end{equation}
and $b_0=b(x,0)$, $b_1=\partial_t b(x,0)$. Note that $\lambda_{+}(\xi)\to -|\xi|^2$ as $|\xi|\to\infty$, and 
\begin{equation}\notag
\lambda_{-}(\xi)=-\frac{2\xi_1^2}{1+\sqrt{1-4\xi_1^2/|\xi|^2}} \sim -\xi_1^2.
\end{equation}
The heuristic argument for the linearized system indicates that solutions near the steady state may converge to the steady state, although the dissipation is weak in the region where $|\xi_1|$ is small.

\medskip

\subsection{Main results}
We focus on the anisotropic system (\ref{emhd-ani}) with $\mu_1=0$ and $\mu_2=1$ and its perturbed form (\ref{emhd-pert}) around the steady state $(a^*,b^*)=(y,0)$. Note that in this case, the first equation of (\ref{emhd-ani}) is a transport equation. We will show that, with only vertical dissipation, global in time small solution exists for (\ref{emhd-pert}) if the initial perturbation is small. Moreover, the solution converges to the steady state with decay rates consistent with that of the linearized system solution. 

Denote $H^{s_1,s_2}$ by the homogeneous anisotropic Sobolev space with the norm
\[\|f\|_{H^{s_1,s_2}}=\||{D_x}|^{s_1}|D|^{s_2}f\|_{L^2}.\]
\begin{Theorem}\label{thm-main}
Let $\psi=a-y$ and $s=\frac12-\epsilon$ for an arbitrarily small constant $\epsilon>0$. Let $(\psi_0,b_0)\in H^7(\mathbb R^2)\times H^6(\mathbb R^2)$ be the initial data for the perturbed system (\ref{emhd-pert}) with $\mu_1=0$ and $\mu_2=1$. There exists a small constant $\varepsilon>0$ such that if 
\[(\nabla \psi_0,b_0)\in (H^{-s,-s}(\mathbb R^2))^2\cap (H^{-s,6}(\mathbb R^2))^2\]
and 
\begin{equation}\notag
\|(\nabla \psi_0,b_0)\|_{(H^{-s,-s})^2}+\|(\nabla \psi_0,b_0)\|_{(H^{-s,6})^2}+\|(\nabla \psi_0,b_0)\|_{(H^6)^2}\leq \varepsilon,
\end{equation}
then system (\ref{emhd-pert}) has a unique global solution $(\psi, b)$ satisfying\\
(i) \[(\nabla\psi, b)\in C([0,\infty); H^6(\mathbb R^2)\times H^6(\mathbb R^2)),\]\\
(ii) \[\|\partial_x^k\nabla\psi(t)\|_{L^2}+\|\partial_x^kb(t)\|_{L^2}\leq C\varepsilon(1+t)^{-\frac{s+k}{2}}, \ \ t\in[0,\infty), \ \ k=0,1,2.\]
\end{Theorem}

\medskip

A few remarks about Theorem \ref{thm-main} are followed.
\begin{Remark}
\item [(1)] The results of Theorem \ref{thm-main} hold for system (\ref{emhd-ani}) with $\mu_1>0$ and $\mu_2>0$ as well. With full dissipation in both horizontal and vertical directions, a trivial modification of the proof of Theorem \ref{thm-main} is sufficient to show the global existence and convergence to steady state.
\item [(2)] The global existence result seems the first one in the literature for the 2D electron MHD which is supercritical, quasi-linear and non-resistive in one direction. 
\item [(3)] We conjecture similar results would be true for the anisotropic system (\ref{emhd-ani}) with $\mu_1>0$ and $\mu_2=0$. This case will be addressed in future work.
\end{Remark}

We conclude this section with the organization of the paper. In Section \ref{sec-pre} we set up notations and introduce the Littlewood-Paley decomposition theory briefly. Section \ref{sec-proof} is devoted to a proof of Theorem \ref{thm-main}, which consists of several energy estimates and decay estimates. Section \ref{sec-aux} provides details to establish a key lemma in Section \ref{sec-proof}, where we need to fully explore cancellations in the energy estimates. Section \ref{sec-app} has the collection of lengthy estimates for another important lemma from Section \ref{sec-proof}.

\bigskip

\section{Preliminaries}
\label{sec-pre}

\subsection{Notations}
We use $C$ to denote a general constant which is not necessary to be tracked; it may be different from line to line.  The symbol $\lesssim$ denotes $\leq$ up to a multiplication of constant when the constant is not of interest.

We denote the norm of the anisotropic Lebesgue space $L_x^pL_y^q(\mathbb R^2)$ by
\[\|f\|_{L_x^pL_y^q}=\left(\int_{\mathbb R} \left(\int_{\mathbb R}|f(x,y)|^q\,dy\right)^{\frac{p}{q}}\,dx\right)^{\frac1p}.\]

\subsection{Littlewood-Paley decomposition}
\label{sec:LPD}
The estimates in the later sections will be obtained mainly through frequency localization techniques. For this purpose, we recall the Littlewood-Paley decomposition theory briefly. 

As customary, $\mathcal F$ and $\mathcal F^{-1}$ denote respectively the Fourier transform and inverse Fourier transform. Let $\lambda_q=2^q$ for integers $q\geq -1$. Define the nonnegative radial function $\chi\in C_0^\infty(\mathbb R^n)$ 
\begin{equation}\notag
\chi(\xi)=
\begin{cases}
1, \ \ \mbox { for } |\xi|\leq\frac{3}{4}\\
0, \ \ \mbox { for } |\xi|\geq 1.
\end{cases}
\end{equation}
Take
\bg\notag
\varphi(\xi)=\chi(\frac{\xi}{2})-\chi(\xi),
\ed
\begin{equation}\notag
\varphi_q(\xi)=
\begin{cases}
\varphi(\lambda_q^{-1}\xi)  \ \ \ \mbox { for } q\geq 0,\\
\chi(\xi) \ \ \ \mbox { for } q=-1.
\end{cases}
\end{equation}
Denote \[\kappa=\mathcal F^{-1}\varphi, \qquad \widetilde \kappa=\mathcal F^{-1}\chi.\]
The Littlewood-Paley projection of a tempered distribution vector field $u$ is defined as
\begin{equation}\notag
\begin{split}
u_q&:=\Delta_qu=\mathcal F^{-1}(\varphi(\lambda_q^{-1}\xi)\mathcal Fu)=\lambda_q^n\int \kappa(\lambda_qy)u(x-y)dy,  \qquad \mbox { for }  q\geq 0,\\
u_{-1}&=\mathcal F^{-1}(\chi(\xi)\mathcal Fu)=\int \widetilde \kappa(y)u(x-y)dy.
\end{split}
\end{equation}
The Littlewood-Paley theory says 
\bg\notag
u=\sum_{q=-1}^\infty u_q
\ed
holds in the distributional sense. 
To simplify the notation, we denote
\bg\notag
u_{\leq q}=S_{q+1}u=\sum_{q'=-1}^qu_{q'}, \ \ u_{(q, p]}=\sum_{q'=q+1}^p u_{q'}, \ \ \widetilde u_q=\widetilde{\Delta}_q u=\sum_{|q-q'|\leq 1}u_{q'}.
\ed
We also define the Littlewood-Paley projection in the horizontal direction
\begin{equation}\notag
u_q^h:=\Delta_q^hu=\mathcal F^{-1}(\varphi(\lambda_q^{-1}\xi_1)\mathcal Fu), \ \ \ 
 u_{-1}=\mathcal F^{-1}(\chi(\xi_1)\mathcal Fu)
\end{equation}
and denote 
\bg\notag
u_{\leq q}^h=S_{q+1}^hu=\sum_{q'=-1}^qu_{q'}^h, \ \ u_{(q, p]}^h=\sum_{q'=q+1}^p u_{q'}^h, \ \ \widetilde u_q^h=\widetilde{\Delta}_q^h u=\sum_{|q-q'|\leq 1}u_{q'}^h.
\ed
In the later sections, we will use the following Bony's paraproducts often
\begin{equation}\notag
\begin{split}
\Delta_k(fg)=& \sum_{|k-k'|\leq 2}\Delta_k\left(S_{k'-1} f\Delta_{k'}g \right)+\sum_{|k-k'|\leq 2}\Delta_k\left( \Delta_{k'}fS_{k'-1}g \right)\\
&+\sum_{k'\geq k-2}\Delta_k(\widetilde f_{k'} \Delta_{k'}g),
\end{split}
\end{equation}
\begin{equation}\notag
\begin{split}
\Delta_k^h(fg)=& \sum_{|k-k'|\leq 2}\Delta_k^h\left(S_{k'-1}^h f\Delta_{k'}^hg \right)+\sum_{|k-k'|\leq 2}\Delta_k^h\left( \Delta_{k'}^hfS_{k'-1}^hg \right)\\
&+\sum_{k'\geq k-2}\Delta_k^h(\widetilde f_{k'}^h \Delta_{k'}^hg),
\end{split}
\end{equation}
\begin{equation}\notag
\begin{split}
\Delta_j\Delta_k^h(fg)=& \sum_{|j-j'|\leq2, |k-k'|\leq 2}\Delta_j\Delta_k^h\left(S_{j'-1}S_{k'-1}^h f\Delta_{j'}\Delta_{k'}^hg \right)\\
&+ \sum_{|j-j'|\leq2, |k-k'|\leq 2}\Delta_j\Delta_k^h\left(\Delta_{j'}S_{k'-1}^h fS_{j'-1}\Delta_{k'}^hg \right)\\
&+\sum_{j'\geq j-2, |k-k'|\leq 2}\Delta_j\Delta_k^h\left(\widetilde{\Delta}_{j'}S_{k'-1}^h f\Delta_{j'}\Delta_{k'}^hg \right)\\
&+\sum_{|j-j'|\leq2, |k-k'|\leq 2}\Delta_j\Delta_k^h\left(S_{j'-1}\Delta_{k'}^hf\Delta_{j'}S_{k'-1}^hg \right)\\
&+\sum_{|j-j'|\leq2, |k-k'|\leq 2}\Delta_j\Delta_k^h\left(\Delta_{j'}\Delta_{k'}^hfS_{j'-1}S_{k'-1}^hg \right)\\
&+\sum_{j'\geq j-2, |k-k'|\leq 2}\Delta_j\Delta_k^h\left(\widetilde{\Delta}_{j'}\Delta_{k'}^hf\Delta_{j'}S_{k'-1}^hg \right)\\
&+\sum_{|j-j'|\leq2, k'\geq k-2}\Delta_j\Delta_k^h(S_{j'-1}\widetilde f_{k'}^h \Delta_{j'}\Delta_{k'}^hg)\\
&+\sum_{|j-j'|\leq2, k'\geq k-2}\Delta_j\Delta_k^h(\Delta_{j'}\widetilde f_{k'}^h S_{j'-1}\Delta_{k'}^hg)\\
&+\sum_{j'\geq j-2, k'\geq k-2}\Delta_j\Delta_k^h(\widetilde{\Delta}_{j'}\widetilde f_{k'}^h \Delta_{j'}\Delta_{k'}^hg)\\
=&: \mathcal B_{j,k}^1(f,g)+...+\mathcal B_{j,k}^9(f,g).
\end{split}
\end{equation}

\begin{Definition}
A tempered distribution $u$ belongs to the anisotropic Besov space $ B^{s_1, s_2}(\mathbb R^2)$ if and only if
$$
\|u\|_{ B^{s_1,s_2}}=\sum_{q=-1}^\infty\lambda_j^{s_2}\lambda_k^{s_1}\|\Delta_j\Delta_k^hu\|_{L^2}<\infty.
$$
\end{Definition}
We also note that, 
\[
  \|u\|_{H^{s_1,s_2}} \sim \left(\sum_{q=-1}^\infty \lambda_j^{2s_2}\lambda_k^{2s_1}\|\Delta_j\Delta_k^hu\|_{L^2}^2\right)^{1/2}.
\]

\begin{Lemma}\label{le:bern} \cite{CDGG} (Bernstein's inequality.)
Let $n$ be the space dimension and $r\geq s\geq 1$. Then for all tempered distributions $u$, 
\bg\notag
\|\Delta_ju\|_{L^r}\lesssim \lambda_j^{n(\frac{1}{s}-\frac{1}{r})}\|\Delta_j u\|_{L^s}.
\ed
In particular, we have for $p,q\geq 2$
\[\|\Delta_j\Delta_k^hu\|_{L_x^pL_y^q}\lesssim \lambda_j^{(\frac{1}{2}-\frac{1}{q})}\lambda_k^{(\frac{1}{2}-\frac{1}{p})}\|\Delta_j\Delta_k^hu\|_{L^2}.\]
\end{Lemma}

We denote the commutator
\[[\Delta_j\Delta_k^h, u\nabla ]v=\Delta_q\Delta_k^h(u\nabla v)-u \Delta_q\Delta_k^h  \nabla v.\]

\begin{Lemma}\label{le-comm1}
For any $1<r_1<\infty$ and $2\leq r_2\leq \infty$ with $\frac1{r_1}=\frac1{r_2}+\frac1{r_3}$, we have
\[ \|[\Delta_j\Delta_k^h, u\nabla ]v\|_{L^{r_1}} \lesssim \|\nabla u\|_{L^{r_2}}\| v\|_{L^{r_3}}. \]
\end{Lemma}
A proof of the commutator estimate can be found in \cite{Dai-hmhd-reg}.

\bigskip

\section{Proof of Theorem \ref{thm-main}}
\label{sec-proof}

We prove Theorem \ref{thm-main} in this section following a sequence of energy estimates and decay estimates.


\subsection{Energy estimate and decay estimate for the linearized system}

Define 
\begin{equation}\notag
\begin{split}
E(t)=&\ \|b(t)\|_{L^2}^2+ \|\nabla\psi(t)\|_{L^2}^2+\|\nabla b(t)\|_{L^2}^2+ \|\nabla^2\psi(t)\|_{L^2}^2+2\varepsilon_1\left<b,\psi_x\right>(t)\\
D(t)=&\ \|\nabla b(t)\|_{L^2}^2+\|\nabla^2 b(t)\|_{L^2}^2+\varepsilon_1 \|\nabla \psi_x(t)\|_{L^2}^2-\varepsilon_1 \|b_x(t)\|_{L^2}^2-\varepsilon_1\left<\Delta b,\psi_x\right>(t)
\end{split}
\end{equation}
with a small constant $\varepsilon_1>0$ such that for some $c>0$
\begin{equation}\notag
\begin{split}
E(t)\geq&\ c\left(\|b(t)\|_{L^2}^2+ \|\nabla\psi(t)\|_{L^2}^2+\|\nabla b(t)\|_{L^2}^2+ \|\nabla^2\psi(t)\|_{L^2}^2\right),\\
D(t)\geq &\ c\left(\|\nabla b(t)\|_{L^2}^2+\|\nabla^2 b(t)\|_{L^2}^2+ \|\nabla \psi_x(t)\|_{L^2}^2\right).
\end{split}
\end{equation}

\begin{Lemma}\label{le-decay}
Let $(\psi, b)$ be the solution of the linearized system (\ref{emhd-linear}) associated with the initial data $(\nabla\psi_0,b_0)\in H^1(\mathbb R^2)\times H^1(\mathbb R^2)$. Assume  $|D_x|^{-s}\nabla \psi_0\in H^{1+s}(\mathbb R^2)$ and $|D_x|^{-s}b_0\in H^{1+s}(\mathbb R^2)$ for $s>0$. We have
\begin{equation}\notag
\|b(t)\|_{L^2}+\|\nabla\psi(t)\|_{L^2}\leq C(1+t)^{-\frac{s}{2}}
\end{equation}
for a constant $C>0$.
\end{Lemma}
\pf
In view of the linearized system (\ref{emhd-linear}) we have the energy laws through standard argument
\begin{equation}\notag
\begin{split}
&\frac12\frac{d}{dt}\left(\|b(t)\|_{L^2}^2+\|\nabla\psi(t)\|_{L^2}^2 \right)\\
=&\int_{\mathbb R^2}bb_t\, dxdy+\int_{\mathbb R^2}\nabla\psi \cdot\nabla\psi_t\, dxdy\\
=& \int_{\mathbb R^2}b\Delta\psi_x\, dxdy +\int_{\mathbb R^2}b\Delta b\, dxdy-\int_{\mathbb R^2}\nabla\psi \cdot\nabla b_x\, dxdy\\
=& \int_{\mathbb R^2}b\Delta\psi_x\, dxdy -\int_{\mathbb R^2}|\nabla b|^2\, dxdy-\int_{\mathbb R^2}b\Delta\psi_x \, dxdy\\
=& -\int_{\mathbb R^2}|\nabla b|^2\, dxdy,
\end{split}
\end{equation}
\begin{equation}\notag
\begin{split}
&\frac12\frac{d}{dt}\left(\|\nabla b(t)\|_{L^2}^2+\|\nabla^2\psi(t)\|_{L^2}^2 \right)\\
=&\int_{\mathbb R^2}\nabla b\cdot\nabla b_t\, dxdy+\int_{\mathbb R^2}\nabla^2\psi \nabla^2\psi_t\, dxdy\\
=& \int_{\mathbb R^2}\nabla b\cdot \Delta\nabla \psi_x\, dxdy +\int_{\mathbb R^2}\nabla b\cdot\Delta \nabla b\, dxdy-\int_{\mathbb R^2}\Delta\psi \Delta b_x\, dxdy\\
=& \int_{\mathbb R^2}\Delta\psi \Delta b_x\, dxdy -\int_{\mathbb R^2}|\nabla^2 b|^2\, dxdy-\int_{\mathbb R^2}\Delta\psi \Delta b_x\, dxdy\\
=& -\int_{\mathbb R^2}|\nabla^2 b|^2\, dxdy
\end{split}
\end{equation}
and
\begin{equation}\notag
\begin{split}
\frac{d}{dt}\varepsilon_1\left<b,\psi_x\right>=&\ \varepsilon_1\left<b_t,\psi_x\right>+\varepsilon_1\left<b_t,\psi_{xt}\right>\\
=&\ \varepsilon_1\int_{\mathbb R^2}\Delta\psi_x\psi_x \, dxdy+ \varepsilon_1\int_{\mathbb R^2}\Delta b\psi_x \, dxdy- \varepsilon_1\int_{\mathbb R^2}bb_{xx} \, dxdy\\
=&-\varepsilon_1\int_{\mathbb R^2}|\nabla \psi_x|^2\, dxdy+ \varepsilon_1\int_{\mathbb R^2}\Delta b\psi_x \, dxdy+\varepsilon_1\int_{\mathbb R^2}|b_x|^2\, dxdy.
\end{split}
\end{equation}
Combining the last three equations yields
\begin{equation}\label{est-ED}
\frac12\frac{d}{dt}E(t)+D(t)\leq 0.
\end{equation}

Denote 
\begin{equation}\notag
\begin{split}
E_s(t)=&\ \||\partial_x|^{-s}b\|_{L^2}^2+\||\partial_x|^{-s}\nabla \psi\|_{L^2}^2+\||D|^{1+s}|\partial_x|^{-s}b\|_{L^2}^2\\
&+\||D|^{1+s}|\partial_x|^{-s}\nabla\psi\|_{L^2}^2.
\end{split}
\end{equation}
Similar analysis as above shows that 
\begin{equation}\notag
\frac{d}{ds} E_s(t)=-\int_{\mathbb R^2}||\partial_x|^{-s} \nabla b|^2\, dxdy-\int_{\mathbb R^2}||D|^{1+s}|\partial_x|^{-s} \nabla b|^2\, dxdy\leq 0.
\end{equation}
By interpolation, we have
\begin{equation}\notag
\begin{split}
\|\nabla\psi\|_{L^2}^2\leq&\ \||\partial_x|^{-s}\nabla \psi\|_{L^2}^{\frac{2}{s+1}}\|\partial_x\nabla \psi\|_{L^2}^{\frac{2s}{s+1}}\\
\leq&\ E_s^{\frac{1}{1+s}}(t)D^{\frac{s}{1+s}}(t)\leq E_s^{\frac{1}{1+s}}(0)D^{\frac{s}{1+s}}(t), \\
\|\nabla^2\psi\|_{L^2}^2\leq&\ \||D|^{s+1}|\partial_x|^{-s}\nabla \psi\|_{L^2}^{\frac{2}{s+1}}\|\partial_x\nabla \psi\|_{L^2}^{\frac{2s}{s+1}}\\
\leq&\ E_s^{\frac{1}{1+s}}(t)D^{\frac{s}{1+s}}(t)\leq E_s^{\frac{1}{1+s}}(0)D^{\frac{s}{1+s}}(t)
\end{split}
\end{equation}
which implies 
\begin{equation}\notag
E(t)\leq c E_s^{\frac{1}{1+s}}(0)D^{\frac{s}{1+s}}(t)
\end{equation}
for some constant $c>0$. Hence 
\begin{equation}\notag
D(t)\geq c^{-1} E_s^{-\frac{1}{s}}(0) E^{1+\frac{1}{s}}(t).
\end{equation}
It then follows from (\ref{est-ED}) that
\begin{equation}\notag
E(t)\leq E(0)(1+c t)^{-s}.
\end{equation}

\cbdu

\medskip

\subsection{Energy estimates for the nonlinear system}

\begin{Lemma}\label{le-energy1}
Assume that the solution $(\psi, b)$ of (\ref{emhd-ani}) satisfies 
\begin{equation}\notag
\|b(t)\|_{H^1}+\|\nabla\psi(t)\|_{H^3}\leq c_0, \ \ \ t\in[0,T)
\end{equation}
for a sufficiently small constant $c_0>0$. There exists a constant $c>0$ such that 
\begin{equation}\notag
\frac{d}{dt}E(t)+cD(t)\leq 0, \ \ \ t\in[0,T).
\end{equation}
\end{Lemma}
\pf
Combining the estimates for the linearized system in the proof of Lemma \ref{le-decay}, we have
\begin{equation}\notag
\begin{split}
&\frac12\frac{d}{dt}E(t)+D(t)\\
=&\int_{\mathbb R^2} b\psi_y\Delta\psi_x\, dxdy-\int_{\mathbb R^2} b\psi_x\Delta\psi_y\, dxdy+\int_{\mathbb R^2} \nabla\psi\nabla(b_y\psi_x)\, dxdy\\
&-\int_{\mathbb R^2} \nabla\psi\nabla(b_x\psi_y)\, dxdy+\int_{\mathbb R^2} \nabla b\nabla(\psi_y\Delta\psi_x)\, dxdy-\int_{\mathbb R^2} \nabla b\nabla(\psi_x\Delta\psi_y)\, dxdy\\
&+\int_{\mathbb R^2} \nabla^2\psi\nabla^2(b_y\psi_x)\, dxdy-\int_{\mathbb R^2} \nabla^2\psi\nabla^2(b_x\psi_y)\, dxdy+\varepsilon_1 \int_{\mathbb R^2}\psi_y\Delta\psi_x\psi_x\, dxdy\\
&-\varepsilon_1 \int_{\mathbb R^2}\psi_x\Delta\psi_y\psi_x\, dxdy+\varepsilon_1 \int_{\mathbb R^2}b\partial_x(b_y\psi_x)\, dxdy-\varepsilon_1 \int_{\mathbb R^2}b\partial_x(b_x\psi_y)\, dxdy\\
=:&\ I_1+I_2+...+I_{12}.
\end{split}
\end{equation}
Applying integration by parts gives
\begin{equation}\notag
\begin{split}
I_1+I_2=&-\int_{\mathbb R^2} b_x\psi_y\Delta\psi\, dxdy-\int_{\mathbb R^2} b\psi_{xy}\Delta\psi\, dxdy\\
&+\int_{\mathbb R^2} b_y\psi_x\Delta\psi\, dxdy+\int_{\mathbb R^2} b\psi_{xy}\Delta\psi\, dxdy\\
=&-\int_{\mathbb R^2} b_x\psi_y\Delta\psi\, dxdy+\int_{\mathbb R^2} b_y\psi_x\Delta\psi\, dxdy\\
=&\int_{\mathbb R^2} \nabla(b_x\psi_y)\nabla\psi\, dxdy-\int_{\mathbb R^2} \nabla(b_y\psi_x)\nabla\psi\, dxdy
\end{split}
\end{equation}
and hence 
\[I_1+I_2+I_3+I_4=0.\]
Again applying integration by parts to $I_5+I_6$ and exploring cancellations yields
\begin{equation}\notag
\begin{split}
I_5+I_6=&-\int_{\mathbb R^2} \nabla b_x \nabla(\psi_y\Delta\psi)\, dxdy-\int_{\mathbb R^2} \nabla b\nabla(\psi_{xy}\Delta\psi)\, dxdy\\
&+\int_{\mathbb R^2} \nabla b_y\nabla(\psi_x\Delta\psi)\, dxdy+\int_{\mathbb R^2} \nabla b\nabla(\psi_{xy}\Delta\psi)\, dxdy\\
=&-\int_{\mathbb R^2} \nabla b_x \nabla(\psi_y\Delta\psi)\, dxdy+\int_{\mathbb R^2} \nabla b_y\nabla(\psi_x\Delta\psi)\, dxdy\\
=&\int_{\mathbb R^2} \nabla^2 b_x \nabla(\psi_y\nabla\psi)\, dxdy+\int_{\mathbb R^2} \nabla b_x \nabla(\nabla\psi_y\nabla\psi)\, dxdy\\
&-\int_{\mathbb R^2} \nabla^2 b_y\nabla(\psi_x\nabla\psi)\, dxdy-\int_{\mathbb R^2} \nabla b_y\nabla(\nabla\psi_x\nabla\psi)\, dxdy\\
=&\int_{\mathbb R^2} \nabla^2 b_x \nabla(\psi_y\nabla\psi)\, dxdy-\int_{\mathbb R^2} \nabla^2 b_y\nabla(\psi_x\nabla\psi)\, dxdy
\end{split}
\end{equation}
since we observe the cancellations
\begin{equation}\notag
\begin{split}
&\int_{\mathbb R^2} \nabla b_x \nabla(\nabla\psi_y\nabla\psi)\, dxdy-\int_{\mathbb R^2} \nabla b_y\nabla(\nabla\psi_x\nabla\psi)\, dxdy\\
=&-\frac12\int_{\mathbb R^2} \nabla b_{xy} \nabla(\nabla\psi\nabla\psi)\, dxdy+\frac12\int_{\mathbb R^2} \nabla b_{yx} \nabla(\nabla\psi\nabla\psi)\, dxdy\\
=&\ 0.
\end{split}
\end{equation}
We continue to rearrange the terms in $I_5+I_6$
\begin{equation}\label{basic-energy-56}
\begin{split}
I_5+I_6
=&\int_{\mathbb R^2} \nabla^2 b_x \nabla\psi_y\nabla\psi\, dxdy+\int_{\mathbb R^2} \nabla^2 b_x \psi_y\nabla^2\psi\, dxdy\\
&-\int_{\mathbb R^2} \nabla^2 b_y\nabla\psi_x\nabla\psi\, dxdy-\int_{\mathbb R^2} \nabla^2 b_y\psi_x\nabla^2\psi\, dxdy\\
=&\int_{\mathbb R^2} \nabla^2 b_x \psi_y\nabla^2\psi\, dxdy-\int_{\mathbb R^2} \nabla^2 b_y\psi_x\nabla^2\psi\, dxdy
\end{split}
\end{equation}
due to the cancellations
\begin{equation}\notag
\begin{split}
&\int_{\mathbb R^2} \nabla^2 b_x \nabla\psi_y\nabla\psi\, dxdy-\int_{\mathbb R^2} \nabla^2 b_y\nabla\psi_x\nabla\psi\, dxdy\\
=&-\frac12\int_{\mathbb R^2} \nabla^2 b_{xy} \nabla\psi\nabla\psi\, dxdy+\frac12\int_{\mathbb R^2} \nabla^2 b_{yx} \nabla\psi\nabla\psi\, dxdy\\
=&\ 0.
\end{split}
\end{equation}
On the other hand, we note
\begin{equation}\label{basic-energy-78}
\begin{split}
I_7+I_8
=&\int_{\mathbb R^2} \nabla^2\psi\nabla^2b_y\psi_x\,dxdy+\int_{\mathbb R^2} \nabla^2\psi b_y\nabla^2\psi_x\,dxdy\\
&+2\int_{\mathbb R^2} \nabla^2\psi\nabla b_y\nabla \psi_x\,dxdy-\int_{\mathbb R^2} \nabla^2\psi\nabla^2b_x\psi_y\,dxdy\\
&-\int_{\mathbb R^2} \nabla^2\psi b_x\nabla^2\psi_y\,dxdy-2\int_{\mathbb R^2} \nabla^2\psi\nabla b_x\nabla \psi_y\,dxdy.
\end{split}
\end{equation}
Combining (\ref{basic-energy-56}) and (\ref{basic-energy-78}) gives
\begin{equation}\notag
\begin{split}
I_5+I_6+I_7+I_8
=&\int_{\mathbb R^2} \nabla^2\psi b_y\nabla^2\psi_x\,dxdy+2\int_{\mathbb R^2} \nabla^2\psi\nabla b_y\nabla \psi_x\,dxdy\\
&-\int_{\mathbb R^2} \nabla^2\psi b_x\nabla^2\psi_y\,dxdy-2\int_{\mathbb R^2} \nabla^2\psi\nabla b_x\nabla \psi_y\,dxdy\\
=&\ 2\int_{\mathbb R^2} \nabla^2\psi\nabla b_y\nabla \psi_x\,dxdy-2\int_{\mathbb R^2} \nabla^2\psi\nabla b_x\nabla \psi_y\,dxdy,
\end{split}
\end{equation}
thanks to the cancellation
\begin{equation}\notag
\begin{split}
&\int_{\mathbb R^2} \nabla^2\psi b_y\nabla^2\psi_x\,dxdy-\int_{\mathbb R^2} \nabla^2\psi b_x\nabla^2\psi_y\,dxdy\\
=&-\frac12 \int_{\mathbb R^2} \nabla^2\psi b_{xy}\nabla^2\psi\,dxdy+\frac12 \int_{\mathbb R^2} \nabla^2\psi b_{yx}\nabla^2\psi\,dxdy\\
=&\ 0.
\end{split}
\end{equation}
We further apply integration by parts to $I_5+I_6+I_7+I_8$,
\begin{equation}\notag
\begin{split}
I_5+I_6+I_7+I_8
=&\ 2\int_{\mathbb R^2} \nabla^2\psi\nabla b_y\nabla \psi_x\,dxdy-2\int_{\mathbb R^2} \nabla^2\psi\nabla b_x\nabla \psi_y\,dxdy\\
=&-2\int_{\mathbb R^2} \nabla^2\psi\nabla b\nabla \psi_{xy}\,dxdy-2\int_{\mathbb R^2} \nabla^2\psi_y\nabla b\nabla \psi_x\,dxdy\\
&+2\int_{\mathbb R^2} \nabla^2\psi\nabla b\nabla \psi_{yx}\,dxdy+2\int_{\mathbb R^2} \nabla^2\psi_x\nabla b\nabla \psi_y\,dxdy\\
=&-2\int_{\mathbb R^2} \nabla^2\psi_y\nabla b\nabla \psi_x\,dxdy+2\int_{\mathbb R^2} \nabla^2\psi_x\nabla b\nabla \psi_y\,dxdy\\
=&-4\int_{\mathbb R^2} \nabla^2\psi_y\nabla b\nabla \psi_x\,dxdy-2\int_{\mathbb R^2} \nabla\psi_x\nabla^2 b\nabla \psi_y\,dxdy.
\end{split}
\end{equation}
Then we infer
\begin{equation}\notag
\begin{split}
\left|\int_{\mathbb R^2} \nabla^2\psi_y\nabla b\nabla \psi_x\,dxdy \right|
\leq&\ \| \nabla^2\psi_y\|_{L^\infty}\|\nabla b\|_{L^2}\|\nabla \psi_x\|_{L^2}\\
\leq &\ C \|\nabla\psi\|_{H^3}\left(\|\nabla b\|_{L^2}^2+\|\nabla \psi_x\|_{L^2}^2\right)\\
\leq&\ C\|\nabla\psi\|_{H^3}D(t),
\end{split}
\end{equation}
and 
\begin{equation}\notag
\begin{split}
\left|\int_{\mathbb R^2} \nabla\psi_x\nabla^2 b\nabla \psi_y\,dxdy \right|
\leq&\ \| \nabla\psi_y\|_{L^\infty}\|\nabla^2 b\|_{L^2}\|\nabla \psi_x\|_{L^2}\\
\leq &\ C \|\nabla\psi\|_{H^2}\left(\|\nabla^2 b\|_{L^2}^2+\|\nabla \psi_x\|_{L^2}^2\right)\\
\leq&\ C\|\nabla\psi\|_{H^2}D(t).
\end{split}
\end{equation}

We first apply integration by parts to $I_9+I_{10}$
\begin{equation}\notag
\begin{split}
I_9+I_{10}=&\ \varepsilon_1\int_{\mathbb R^2}\psi_y\Delta\psi_x\psi_x \, dxdy-\varepsilon_1\int_{\mathbb R^2}\psi_x\Delta\psi_y\psi_x \, dxdy\\
=&-\varepsilon_1\int_{\mathbb R^2}\nabla\psi_y\nabla\psi_x\psi_x \, dxdy-\varepsilon_1\int_{\mathbb R^2}\psi_y\nabla\psi_x\nabla\psi_x \, dxdy\\
&+\varepsilon_1\int_{\mathbb R^2}\nabla\psi_x\nabla\psi_y\psi_x \, dxdy+\varepsilon_1\int_{\mathbb R^2}\psi_x\nabla\psi_y\nabla\psi_x \, dxdy\\
=&-\varepsilon_1\int_{\mathbb R^2}\psi_y\nabla\psi_x\nabla\psi_x \, dxdy+\varepsilon_1\int_{\mathbb R^2}\nabla\psi_x\nabla\psi_y\psi_x \, dxdy\\
\end{split}
\end{equation}
and then
\begin{equation}\notag
\begin{split}
\left|\int_{\mathbb R^2}\psi_y\nabla\psi_x\nabla\psi_x \, dxdy \right|
\leq&\ \|\psi_y\|_{L^\infty}\|\nabla \psi_x\|_{L^2}^2\\
\leq&\ C\|\nabla\psi\|_{H^1}D(t),
\end{split}
\end{equation}
\begin{equation}\notag
\begin{split}
\left|\int_{\mathbb R^2}\nabla\psi_x\nabla\psi_y\psi_x \, dxdy \right|
\leq&\ \|\psi_x\|_{L^\infty}\|\nabla \psi_y\|_{L^2}\|\nabla \psi_x\|_{L^2}\\
\leq&\ C\|\nabla \psi_y\|_{L^2}\|\nabla \psi_x\|_{L^2}^2\\
\leq&\ C\|\nabla\psi\|_{H^1}D(t).
\end{split}
\end{equation}
Similarly we have
\begin{equation}\notag
\begin{split}
I_{11}+I_{12}=&\ \varepsilon_1\int_{\mathbb R^2}b\partial_x(b_y\psi_x) \, dxdy-\varepsilon_1\int_{\mathbb R^2}b\partial_x(b_x\psi_y) \, dxdy\\
=&- \varepsilon_1\int_{\mathbb R^2}b_xb_y\psi_x\, dxdy+\varepsilon_1\int_{\mathbb R^2}b_xb_x\psi_y \, dxdy
\end{split}
\end{equation}
and hence
\begin{equation}\notag
\begin{split}
\left|\int_{\mathbb R^2}b_xb_y\psi_x\, dxdy \right|
\leq&\ \|\psi_x\|_{L^\infty}\|b_x\|_{L^2}\|b_y\|_{L^2}\\
\leq&\ C\|\nabla\psi\|_{H^1}D(t),
\end{split}
\end{equation}
\begin{equation}\notag
\begin{split}
\left|\int_{\mathbb R^2}b_xb_x\psi_y\, dxdy \right|
\leq&\ \|\psi_y\|_{L^\infty}\|b_x\|_{L^2}^2\\
\leq&\ C\|\nabla\psi\|_{H^1}D(t).
\end{split}
\end{equation}
Choosing $c_0>0$ small enough concludes the proof of the lemma.
\cbdu

\medskip

For $l=1,2$, define 
\begin{equation}\notag
\begin{split}
E_l(t)=&\sum_{j,k}\lambda_k^{2l}\left( \|\Delta_j\Delta_k^hb(t)\|_{L^2}^2+ \|\Delta_j\Delta_k^h\nabla\psi(t)\|_{L^2}^2+\|\Delta_j\Delta_k^h\nabla b(t)\|_{L^2}^2\right. \\
& \left. + \|\Delta_j\Delta_k^h\nabla^2\psi(t)\|_{L^2}^2+2\varepsilon_1\left<\Delta_j\Delta_k^hb, \Delta_j\Delta_k^h\psi_x\right>(t)\right)\\
D_l(t)=&\sum_{j,k}\lambda_k^{2l}\left( \|\Delta_j\Delta_k^h\nabla b(t)\|_{L^2}^2+\|\Delta_j\Delta_k^h\nabla^2 b(t)\|_{L^2}^2+\varepsilon_1 \|\Delta_j\Delta_k^h\nabla \psi_x(t)\|_{L^2}^2\right.\\
&\left.-\varepsilon_1 \|\Delta_j\Delta_k^hb_x(t)\|_{L^2}^2-\varepsilon_1\left<\Delta_j\Delta_k^h\Delta b,\Delta_j\Delta_k^h\psi_x\right>(t) \right)
\end{split}
\end{equation}
where $\varepsilon_1>0$ is small enough such that
\begin{equation}\label{ED}
\begin{split}
E_l(t)\geq&\ c\left(\|\partial_x^l b\|_{L^2}^2+\|\partial_x^l \nabla b\|_{L^2}^2+\|\partial_x^l \nabla\psi\|_{L^2}^2+\|\partial_x^l \nabla^2\psi\|_{L^2}^2 \right),\\
D_l(t)\geq&\ c\left(\|\partial_x^l \nabla b\|_{L^2}^2+\|\partial_x^l \nabla^2b\|_{L^2}^2+\|\partial_x^{l+1} \nabla\psi\|_{L^2}^2 \right)
\end{split}
\end{equation}
for some constant $c>0$.

\begin{Lemma}\label{le-energy2}
Let $\epsilon>0$ be arbitrarily small and $s=\frac12-\epsilon$. Assume that the solution $(\psi, b)$ of (\ref{emhd-pert}) satisfies 
\begin{equation}\label{small}
\sup_{t\in[0,T]}\left(\|(b,\nabla\psi)\|_{H^{-s,-s}}+\|(b,\nabla\psi)\|_{H^{-s,6}}+\|(b,\nabla\psi)\|_{H^6}\right)\leq c_0
\end{equation}
for a sufficiently small constant $c_0>0$. There exists a constant $c>0$ such that 
\begin{equation}\notag
\frac{d}{dt}E_l(t)+cD_l(t)\leq 0, \ \ \ t\in[0,T], \ \ \ l=1,2.
\end{equation}
\end{Lemma}
\pf
Acting the localization operators $\Delta_j$ and $\Delta_k^h$ on (\ref{emhd-pert}) and performing energy estimate, we obtain
\begin{equation}\label{basic-energy}
\begin{split}
&\frac12\frac{d}{dt}\left(\|\Delta_j\Delta_k^hb\|_{L^2}^2+\|\Delta_j\Delta_k^h\nabla\psi\|_{L^2}^2+\|\Delta_j\Delta_k^h\nabla b\|_{L^2}^2+\|\Delta_j\Delta_k^h\nabla^2\psi\|_{L^2}^2\right.\\
&\left.+2\varepsilon_1\left<\Delta_j\Delta_k^h b, \Delta_j\Delta_k^h\psi_x\right>  \right)+\left(\|\Delta_j\Delta_k^h\nabla b\|_{L^2}^2+\|\Delta_j\Delta_k^h\nabla^2 b\|_{L^2}^2\right.\\
&\left. +\varepsilon_1 \|\Delta_j\Delta_k^h\nabla \psi_x\|_{L^2}^2 -\varepsilon_1 \|\Delta_j\Delta_k^h b_x\|_{L^2}^2-\varepsilon_1\left<\Delta_j\Delta_k^h \Delta b, \Delta_j\Delta_k^h\psi_x\right> \right)\\
=&\int_{\mathbb R^2}\Delta_j\Delta_k^hb \Delta_j\Delta_k^h(\psi_y\Delta \psi_x) \, dxdy-\int_{\mathbb R^2}\Delta_j\Delta_k^hb \Delta_j\Delta_k^h(\psi_x\Delta \psi_y) \, dxdy\\
&+\int_{\mathbb R^2}\Delta_j\Delta_k^h\nabla\psi \Delta_j\Delta_k^h\nabla(b_y \psi_x) \, dxdy
-\int_{\mathbb R^2}\Delta_j\Delta_k^h\nabla\psi \Delta_j\Delta_k^h\nabla(b_x \psi_y) \, dxdy\\
&+\int_{\mathbb R^2}\Delta_j\Delta_k^h\nabla b \Delta_j\Delta_k^h\nabla(\psi_y\Delta \psi_x) \, dxdy-\int_{\mathbb R^2}\Delta_j\Delta_k^h\nabla b \Delta_j\Delta_k^h\nabla(\psi_x\Delta \psi_y) \, dxdy\\
&+\int_{\mathbb R^2}\Delta_j\Delta_k^h\nabla^2\psi \Delta_j\Delta_k^h\nabla^2(b_y \psi_x) \, dxdy
-\int_{\mathbb R^2}\Delta_j\Delta_k^h\nabla^2\psi \Delta_j\Delta_k^h\nabla^2(b_x \psi_y) \, dxdy\\
&+\varepsilon_1 \int_{\mathbb R^2}\Delta_j\Delta_k^h\psi_x \Delta_j\Delta_k^h(\psi_y \Delta\psi_x) \, dxdy
-\varepsilon_1 \int_{\mathbb R^2}\Delta_j\Delta_k^h\psi_x \Delta_j\Delta_k^h(\psi_x \Delta\psi_y) \, dxdy\\
&+\varepsilon_1 \int_{\mathbb R^2}\Delta_j\Delta_k^hb \Delta_j\Delta_k^h\partial_x(b_y\psi_x) \, dxdy
-\varepsilon_1 \int_{\mathbb R^2}\Delta_j\Delta_k^hb \Delta_j\Delta_k^h\partial_x(b_x \psi_y) \, dxdy\\
=:&\ \tilde I_1+\tilde I_2+...+\tilde I_{12}.
\end{split}
\end{equation}
It follows from Lemmas in the Appendix that 
\begin{equation}\notag
\frac{d}{dt}E_l(t)+D_l(t)\leq Cc_0D_l(t), \ \ \ t\in[0,T), \ \ \ l=1,2
\end{equation}
for some constant $C>0$. Taking $c_0$ small enough yields the conclusion of the lemma.

\cbdu

\medskip

\subsection{Higher order energy estimates for the nonlinear system}

We define the higher order energy quantity
\begin{equation}\notag
E_{s,s_1}(t)=\|b(t)\|_{H^{-s,s_1}}^2+\|\nabla\psi(t)\|_{H^{-s,s_1}}^2+\|b(t)\|_{H^{-s,s_1+1}}^2+\|\nabla\psi(t)\|_{H^{-s,s_1+1}}^2.
\end{equation}
We also denote 
\begin{equation}\notag
\begin{split}
g_1(t)=&\ \|\nabla b(t)\|_{B^{\frac12,\frac32}}+\|\nabla\psi(t)\|_{B^{\frac32,\frac12}}+\|\nabla\psi(t)\|_{B^{\frac12,\frac32}}\\
g_2(t)=&\ \|b(t)\|_{B^{\frac32,\frac12}}\\
g_3(t)=&\ \|b(t)\|_{B^{\frac12,\frac32}}+\|b(t)\|_{B^{\frac32,\frac12}}+\|\psi(t)\|_{B^{\frac32,\frac12}}\\
&+\|\nabla\psi(t)\|_{B^{\frac32,\frac12}}+\|\psi(t)\|_{B^{\frac12,\frac32}}^2+\|\nabla\psi(t)\|_{B^{\frac12,\frac32}}^2.
\end{split}
\end{equation}

\begin{Lemma}\label{le-energy3}
Let $s\in[0,\frac12)$ and $s_1>-\frac12$. Assume that $(\psi, b)$ is a solution of (\ref{emhd-pert}) on $[0,T)$. There exist constants $c_1>0$ and $C>0$ such that if $g_1(t)\leq c_1$ for $t\in[0,T)$, we have
\begin{equation}\notag
E_{s,s_1}(t)\leq CE_{s,s_1}(0)+C\int_0^t \left(g_2(\tau)+g_3^2(\tau)\right) E_{s,s_1}(\tau)\, d\tau, \ \ \ t\in[0,T).
\end{equation}
\end{Lemma}
\pf
It follows from the basic energy inequality (\ref{basic-energy}) that
\begin{equation}\label{energy-high-1}
\begin{split}
&E_{s,s_1}(t)+c\left(\|\nabla b\|_{L^2_tH_x^{-s,s_1}}^2+\|\nabla b\|_{L^2_tH_x^{-s,s_1+1}}^2+\|\nabla \partial_x\psi\|_{L^2_tH_x^{-s,s_1}}^2\right)\\
\leq&\ E_{s,s_1}(0)+\sum_{j,k}\lambda_j^{2s_1}\lambda_k^{-2s}\int_0^t\left(\tilde I_1(\tau)+...+  \tilde I_{12}(\tau)\right) \, d\tau.
\end{split}
\end{equation}
On the other hand, combining the auxiliary Lemma \ref{le-high-aux12}- Lemma \ref{le-high-aux11} we obtain
\begin{equation}\notag
\begin{split}
&\sum_{j,k}\lambda_j^{2s_1}\lambda_k^{-2s}\int_0^t\left(\tilde I_1(\tau)+...+  \tilde I_{12}(\tau)\right) \, d\tau\\
\leq &\left(C\sup_{\tau\in[0,t)}g_1(\tau)+\frac12c \right)\left(\|\nabla b\|_{L^2_tH_x^{-s,s_1}}^2+\|\nabla b\|_{L^2_tH_x^{-s,s_1+1}}^2+\|\nabla \partial_x\psi\|_{L^2_tH_x^{-s,s_1}}^2\right)\\
&+C\int_0^t\left(g_2(\tau)+g_3^2(\tau)\right)E_{s,s_1}(\tau)\, d\tau.
\end{split}
\end{equation}
For small enough constant $c_1$, the conclusion of the lemma follows from the last two inequalities.

\cbdu

\begin{Lemma}\label{le-high-decay}
Let $(b,\psi)$ be a solution of (\ref{emhd-pert}) satisfying the smallness condition (\ref{small}) and 
\begin{equation}\label{assu-E}
\sup_{t\in[0,T), l=0,1,2}\left(E_{s,0}(t)+E_{s,l+s}(t)\right) \leq C\varepsilon^2
\end{equation}
for a constant $C>0$. We have
\begin{equation}\notag
\begin{split}
\|\partial_x^lb\|_{L^2}+\|\partial_x^l\nabla b\|_{L^2}+\|\partial_x^l\nabla\psi\|_{L^2}+\|\partial_x^l\nabla^2\psi\|_{L^2}&\leq C\varepsilon (1+t)^{-\frac{s+l}{2}},\\
\int_0^t\left(\|\partial_x^l\nabla b(\tau)\|_{L^2}^2+\|\partial_x^l\nabla^2 b(\tau)\|_{L^2}^2 +\|\partial_x^{l}\nabla \psi_x(\tau)\|_{L^2}^2\right)\, d\tau &\leq C\varepsilon^2
\end{split}
\end{equation}
for $t\in[0,T)$ and $l=0,1,2$ and a constant $C>0$.
\end{Lemma}
\pf
It is shown in Lemma \ref{le-energy1} and Lemma \ref{le-energy2} that
\begin{equation}\label{energy-ED}
\frac{d}{dt}E_l(t)+cD_l(t)\leq 0, \ \ \ t\in[0,T), \ \ \ l=0, 1,2.
\end{equation}
It then follows from (\ref{ED}) 
\begin{equation}\notag
E_l(t)+c\int_0^t\left(\|\partial_x^l\nabla b(\tau)\|_{L^2}^2+\|\partial_x^l\nabla^2 b(\tau)\|_{L^2}^2 +\|\partial_x^{l}\nabla \psi_x(\tau)\|_{L^2}^2\right)\, d\tau \leq E_l(0).
\end{equation}
Since the assumption (\ref{assu-E}) implies $E_l(0)\leq C\varepsilon^2$, we deduce 
\begin{equation}\notag
\int_0^t\left(\|\partial_x^l\nabla b(\tau)\|_{L^2}^2+\|\partial_x^l\nabla^2 b(\tau)\|_{L^2}^2 +\|\partial_x^{l}\nabla \psi_x(\tau)\|_{L^2}^2\right)\, d\tau \leq C\varepsilon^2.
\end{equation}

Applying interpolation yields
\begin{equation}\notag
\begin{split}
\|\partial_x^l\nabla\psi\|_{L^2}^2&\leq \||\partial_x|^{-s}\nabla\psi\|_{L^2}^{\frac{2}{s+l+1}}\|\partial_x^{l+1}\nabla\psi\|_{L^2}^{\frac{2(l+s)}{s+l+1}},\\
\|\partial_x^l\nabla^2\psi\|_{L^2}^2&\leq \||\partial_x|^{-s}D^{s+l+1}\nabla\psi\|_{L^2}^{\frac{2}{s+l+1}}\|\partial_x^{l+1}\nabla\psi\|_{L^2}^{\frac{2(l+s)}{s+l+1}}
\end{split}
\end{equation}
and similar interpolation inequalities hold for $\|\partial_x^lb\|_{L^2}^2$ and $\|\partial_x^l\nabla b\|_{L^2}^2$. Therefore we infer
\begin{equation}\label{interpo}
E_l(t)\leq C\left(E_{s,0}(t)+E_{s,l+s}(t)\right)^{\frac{1}{s+l+1}}D_l^{\frac{s+l}{s+l+1}}(t).
\end{equation}
In view of (\ref{assu-E}), (\ref{energy-ED}) and (\ref{interpo}) we obtain
\begin{equation}\notag
\frac{d}{dt}E_l(t)+c\varepsilon^{-\frac{2}{s+l}} E_l^{1+\frac{1}{s+l}}(t)\leq 0
\end{equation}
which along with (\ref{ED}) implies 
\begin{equation}\notag
\|\partial_x^lb\|_{L^2}+\|\partial_x^l\nabla b\|_{L^2}+\|\partial_x^l\nabla\psi\|_{L^2}+\|\partial_x^l\nabla^2\psi\|_{L^2}\leq C\varepsilon (1+t)^{-\frac{s+l}{2}}.
\end{equation}

\cbdu

\medskip

\begin{Lemma}\label{le-gg}
Let $s>\frac13$. Assume 
\begin{equation}\label{assu-Es}
E_{0,5}(t)+E_{s,-s}(t)+E_{s,5}(t)\leq c_1\varepsilon^2, \ \ t\in[0,T).
\end{equation}
The estimates
\begin{equation}\notag
g_1(t)\leq C\varepsilon, \ \ \int_0^t g_2(\tau)\, d\tau\leq C\varepsilon, \ \ \int_0^t g_3^2(\tau)\, d\tau\leq C\varepsilon^2, \ \ t\in[0,T) 
\end{equation}
hold for a constant $C>0$.
\end{Lemma}
\pf
It follows from the anisotropic Besov space definition that
\begin{equation}\notag
\begin{split}
\|b\|_{B^{\frac12,\frac32}}&\leq \sum_{k\leq j}\lambda_j^{\frac32}\lambda_k^{\frac12}\|\Delta_j\Delta_k^hb\|_{L^2}\\
&\leq \sum_{j}\lambda_j^{2}\|\Delta_j\Delta_k^hb\|_{L^2}\\
&\leq \sum_{j}\|\Delta_j\Delta_k^h\nabla b\|_{L^2}^{\frac12}\|\Delta_j\Delta_k^h\nabla^3 b\|_{L^2}^{\frac12}\\
&\leq C\left(\|\nabla b\|_{L^2}+\|\nabla^3 b\|_{L^2}\right),
\end{split}
\end{equation}
and similarly
\begin{equation}\notag
\begin{split}
\|\nabla b\|_{B^{\frac12,\frac32}}&\leq \sum_{k\leq j}\lambda_j^{\frac32}\lambda_k^{\frac12}\|\Delta_j\Delta_k^h\nabla b\|_{L^2}\\
&\leq C\left(\|\nabla^2 b\|_{L^2}+\|\nabla^4 b\|_{L^2}\right),
\end{split}
\end{equation}
\begin{equation}\notag
\begin{split}
\|\nabla\psi\|_{B^{\frac12,\frac32}}&\leq \sum_{k\leq j}\lambda_j^{\frac32}\lambda_k^{\frac12}\|\Delta_j\Delta_k^h\nabla\psi\|_{L^2}\\
&\leq \sum_{j}\lambda_j^{2}\|\Delta_j\Delta_k^h\nabla\psi\|_{L^2}\\
&\leq \sum_{j}\|\Delta_j\Delta_k^h\nabla^2\psi\|_{L^2}^{\frac12}\|\Delta_j\Delta_k^h\nabla^4\psi\|_{L^2}^{\frac12}\\
&\leq C\left(\|\nabla^2\psi\|_{L^2}+\|\nabla^4 \psi\|_{L^2}\right),
\end{split}
\end{equation}
\begin{equation}\notag
\begin{split}
\|\nabla\psi\|_{B^{\frac32,\frac12}}&\leq \sum_{k\leq j}\lambda_j^{\frac12}\lambda_k^{\frac32}\|\Delta_j\Delta_k^h\nabla\psi\|_{L^2}\\
&\leq \sum_{j}\lambda_j^{2}\|\Delta_j\Delta_k^h\nabla\psi\|_{L^2}\\
&\leq C\left(\|\nabla^2\psi\|_{L^2}+\|\nabla^4 \psi\|_{L^2}\right).
\end{split}
\end{equation}
Hence 
\begin{equation}\notag
g_1(t)\leq CE_{0,3}^{\frac12}(t)\leq C\varepsilon.
\end{equation}

In view of 
\begin{equation}\notag
\begin{split}
\|b\|_{B^{\frac32,\frac12}}&\leq \sum_{k\leq j}\lambda_j^{\frac12}\lambda_k^{\frac32}\|\Delta_j\Delta_k^hb\|_{L^2}\\
&\leq \sum_{j}\lambda_j^{\frac12}\|\Delta_jb\|_{L^2}^{\frac14}\|\Delta_j\partial_x^2b\|_{L^2}^{\frac34}\\
&\leq \sum_{j}\|\Delta_j\nabla b\|_{L^2}^{\frac14}\|\Delta_j\partial_x^2b\|_{L^2}^{\frac12}\|\Delta_j\partial_x^2\nabla b\|_{L^2}^{\frac14}\\
&=\sum_{j}\left(\|\Delta_j\nabla b\|_{L^2}^{\frac16}\|\Delta_j\partial_x^2b\|_{L^2}^{\frac12}\right)\left(\|\Delta_j\nabla b\|_{L^2}^{\frac{1}{12}}\|\Delta_j\partial_x^2\nabla b\|_{L^2}^{\frac14}\right)\\
&\leq C\|\nabla b\|_{L^2}^{\frac14}\|\partial_x^2b\|_{L^2}^{\frac34}+C\|\nabla b\|_{L^2}^{\frac{1}{4}}\|\partial_x^2\nabla b\|_{L^2}^{\frac34},
\end{split}
\end{equation}
we infer from (\ref{assu-Es}) and Lemma \ref{le-high-decay} 
\begin{equation}\notag
\begin{split}
\int_0^t g_2(\tau)\, d\tau=&\int_0^t \|b(\tau)\|_{B^{\frac32,\frac12}}\, d\tau\\
\leq&\ C\int_0^t\|\nabla b(\tau)\|_{L^2}^{\frac14}\|\partial_x^2b(\tau)\|_{L^2}^{\frac34} \, d\tau+C\int_0^t\|\nabla b(\tau)\|_{L^2}^{\frac14}\|\partial_x^2\nabla b(\tau)\|_{L^2}^{\frac34} \, d\tau\\
\leq&\ C\varepsilon^2 \left(\int_0^t\|\nabla b(\tau)\|_{L^2}^{2} \, d\tau \right)^{\frac18}\left(\int_0^t\|\partial_x^2 b(\tau)\|_{L^2}^{\frac67} \, d\tau \right)^{\frac78}\\
&+C\varepsilon^2 \left(\int_0^t\|\nabla b(\tau)\|_{L^2}^{2} \, d\tau \right)^{\frac18}\left(\int_0^t\|\partial_x^2 \nabla b(\tau)\|_{L^2}^{\frac67} \, d\tau \right)^{\frac78}\\
\leq&\ C\varepsilon^2\left(\int_0^t(1+\tau)^{-\frac{s+2}{2}\cdot\frac67} \, d\tau \right)^{\frac78}\\
\leq&\ C\varepsilon^2
\end{split}
\end{equation}
for $s>\frac13$.

On the other hand, due to
\begin{equation}\notag
\begin{split}
\|\psi\|_{B^{\frac32,\frac12}}&\leq \sum_{k\leq j}\lambda_j^{\frac12}\lambda_k^{\frac32}\|\Delta_j\Delta_k^h\psi\|_{L^2}\\
&\leq \sum_{j}\lambda_j^{\frac12}\|\Delta_j\partial_x\psi\|_{L^2}^{\frac12}\|\Delta_j\partial_x^2\psi\|_{L^2}^{\frac12}\\
&\leq \left(\||D_x|^{-s}\partial_x\psi\|_{L^2}+\|\nabla\partial_x\psi\|_{L^2} \right)^{\frac12}\|\nabla\partial_x^2\psi\|_{L^2}^{\frac12},
\end{split}
\end{equation}
it follows from assumption (\ref{assu-Es}) and Lemma \ref{le-high-decay} that
\begin{equation}\notag
\begin{split}
\int_0^t\|\psi(\tau)\|_{B^{\frac32,\frac12}}^2 \, d\tau\leq& \int_0^t\left(\||D_x|^{-s}\partial_x\psi(\tau)\|_{L^2}+\|\nabla\partial_x\psi(\tau)\|_{L^2} \right)\|\nabla\partial_x^2\psi(\tau)\|_{L^2} \, d\tau\\
\leq& \left(\sup_{\tau\in[0,t)}\||D_x|^{-s}\partial_x\psi(\tau)\|_{L^2}\right) \int_0^t\|\nabla\partial_x^2\psi(\tau)\|_{L^2}\,d\tau\\
&+\left( \int_0^t\|\nabla\partial_x\psi(\tau)\|_{L^2}^2\,d\tau\right)^{\frac12}\left( \int_0^t\|\nabla\partial_x^2\psi(\tau)\|_{L^2}^2\,d\tau\right)^{\frac12}\\
\leq&\ C\varepsilon \int_0^t\varepsilon (1+\tau)^{-\frac{s+2}{2}}\, d\tau+C\varepsilon \left( \int_0^t\varepsilon^2(1+\tau)^{-(s+2)}\,d\tau\right)^{\frac12}\\
\leq&\ C\varepsilon^2
\end{split}
\end{equation}
as $0<s<\frac12$.

Since for some $0<\theta<\frac{s}{4+2s}<\frac12$,
\begin{equation}\notag
\begin{split}
\|b\|_{B^{\frac12,\frac32}}&\leq \sum_{k\leq j}\lambda_j^{\frac32}\lambda_k^{\frac12}\|\Delta_j\Delta_k^hb\|_{L^2}\\
&\leq \sum_{j}\lambda_j^{\frac32}\|\Delta_jb\|_{L^2}^{\frac12}\|\Delta_j\partial_x b\|_{L^2}^{\frac12}\\
&\leq \sum_j\left(\|\Delta_j\partial_xb\|_{L^2}^{\theta}\|\Delta_j\nabla b\|_{L^2}^{\theta}  \right)\left(\|\Delta_jD^{\frac{2}{1-2\theta}}\partial_xb\|_{L^2}^{\frac12-\theta}\|\Delta_j\nabla b\|_{L^2}^{\frac12-\theta}  \right)\\
&\leq C\|\partial_xb\|_{L^2}^{\frac12}\|\nabla b\|_{L^2}^{\frac12}+C\|D^{\frac{2}{1-2\theta}}\partial_xb\|_{L^2}^{\frac12}\|\nabla b\|_{L^2}^{\frac12}\\
&\leq C\|\partial_xb\|_{L^2}^{\frac12}\|\nabla b\|_{L^2}^{\frac12}+C\|\nabla\partial_xb\|_{L^2}^{\frac12}\|\nabla^2 b\|_{L^2}^{\frac{1-6\theta}{2-4\theta}}\|\nabla^3 b\|_{L^2}^{\frac{2\theta}{1-2\theta}},
\end{split}
\end{equation}
we deduce from (\ref{assu-Es}) and Lemma \ref{le-high-decay} 
\begin{equation}\notag
\begin{split}
&\int_0^t\|b(\tau)\|_{B^{\frac12,\frac32}}^2 \, d\tau\\
\leq&\ C\int_0^t \|\partial_xb(\tau)\|_{L^2}\|\nabla b\|_{L^2}\, d\tau\\
&+C\int_0^t \|\nabla\partial_xb(\tau)\|_{L^2}\|\nabla^2 b(\tau)\|_{L^2}^{\frac{1-6\theta}{1-2\theta}}\|\nabla^3 b(\tau)\|_{L^2}^{\frac{4\theta}{1-2\theta}}\, d\tau\\
\leq&\ C\left(\int_0^t \|\partial_xb(\tau)\|_{L^2}^2\, d\tau\right)^{\frac12}\left(\int_0^t\|\nabla b\|_{L^2}^2\, d\tau\right)^{\frac12}\\
&+C\left(\int_0^t \|\partial_xb(\tau)\|_{L^2}^{\frac{2-4\theta}{1+2\theta}}\, d\tau\right)^{\frac{1+2\theta}{2-4\theta}}\left(\int_0^t\|\nabla^2 b\|_{L^2}^2\|\nabla^3 b\|_{L^2}^{\frac{8\theta}{1-6\theta}}\, d\tau\right)^{\frac{1-6\theta}{2-4\theta}}\\
\leq&\ C\varepsilon^2\left(\int_0^t (1+\tau)^{-(s+1)}\, d\tau\right)^{\frac12}\\
&+C\varepsilon^2\left(\int_0^t (1+\tau)^{-(s+1)\frac{1-2\theta}{1+2\theta}}\, d\tau\right)^{\frac{1+2\theta}{2-4\theta}}\left(\int_0^t\|\nabla^2 b\|_{L^2}^2\, d\tau\right)^{\frac{1-6\theta}{2-4\theta}}\\
\leq&\ C\varepsilon^2
\end{split}
\end{equation}
where in the last step we used the fact $(s+1)\frac{1-2\theta}{1+2\theta}>1$. Noting that $\|b\|_{B^{\frac32,\frac12}}\leq \|b\|_{B^{\frac12,\frac32}}$, we also have \[\int_0^t\|b(\tau)\|_{B^{\frac32,\frac12}}^2 \, d\tau\leq C\varepsilon^2.\]
Analogously for some $0<\theta<\frac12$, using interpolation
\begin{equation}\notag
\begin{split}
\|\nabla\psi\|_{B^{\frac12,\frac32}}&\leq \sum_{k\leq j}\lambda_j^{\frac32}\lambda_k^{\frac12}\|\Delta_j\Delta_k^h\nabla\psi\|_{L^2}\\
&\leq \sum_{j}\lambda_j^{\frac32}\|\Delta_j\nabla\psi\|_{L^2}^{\frac12}\|\Delta_j\partial_x \nabla\psi\|_{L^2}^{\frac12}\\
&\leq \sum_j\left(\|\Delta_j\partial_x\nabla\psi\|_{L^2}^{\theta}\|\Delta_j\nabla^2 \psi\|_{L^2}^{\theta}  \right)\left(\|\Delta_jD^{\frac{2}{1-2\theta}}\nabla\partial_x\psi\|_{L^2}^{\frac12-\theta}\|\Delta_j\nabla^2 \psi\|_{L^2}^{\frac12-\theta}  \right)\\
&\leq C\|\nabla\partial_x\psi\|_{L^2}^{\theta}\|\nabla^2 \psi\|_{L^2}^{\theta}\|\nabla^2\partial_x\psi\|_{L^2}^{\frac12-\theta}\|\nabla^3 \psi\|_{L^2}^{\frac{1-6\theta}{2}}\|\nabla^4 \psi\|_{L^2}^{2\theta},
\end{split}
\end{equation}
we thus obtain from (\ref{assu-Es}) and Lemma \ref{le-high-decay} 
\begin{equation}\notag
\begin{split}
&\int_0^t\|\nabla\psi(\tau)\|_{B^{\frac12,\frac32}}^4 \, d\tau\\
\leq&\ C\int_0^t \|\nabla\partial_x\psi\|_{L^2}^{4\theta}\|\nabla^2 \psi\|_{L^2}^{4\theta}\|\nabla^2\partial_x\psi\|_{L^2}^{2(1-2\theta)}\|\nabla^3 \psi\|_{L^2}^{2(1-6\theta)}\|\nabla^4 \psi\|_{L^2}^{8\theta}\,d\tau\\
\leq&\ C\varepsilon^2 \int_0^t (1+\tau)^{-\frac12(s+1)\cdot 4\theta}(1+\tau)^{-\frac12(s+1)\cdot 2(1-2\theta)}\, d\tau\\
\leq&\ C\varepsilon^2 \int_0^t (1+\tau)^{-(s+1)}\, d\tau\\
\leq&\ C\varepsilon^2.
\end{split}
\end{equation}
Again it follows from 
\begin{equation}\notag
\begin{split}
\|\nabla\psi\|_{B^{\frac32,\frac12}}&\leq \sum_{k\leq j}\lambda_j^{\frac12}\lambda_k^{\frac32}\|\Delta_j\Delta_k^h\nabla\psi\|_{L^2}\\
&\leq \sum_{j}\lambda_j^{\frac12}\|\Delta_j\partial_x\nabla\psi\|_{L^2}^{\frac12}\|\Delta_j\partial_x^2 \nabla\psi\|_{L^2}^{\frac12}\\
&\leq \sum_j\|\Delta_j\partial_x\nabla\psi\|_{L^2}^{\frac12}\|\Delta_j\partial_x^2 \nabla^2\psi\|_{L^2}^{\frac12}\\
&\leq C\|\partial_x\nabla\psi\|_{L^2}^{\frac12}\|\partial_x^2\nabla^2\psi\|_{L^2}^{\frac12},
\end{split}
\end{equation}
and Lemma \ref{le-high-decay} that
\begin{equation}\notag
\begin{split}
&\int_0^t\|\nabla\psi(\tau)\|_{B^{\frac32,\frac12}}^2 \, d\tau\\
\leq&\ C\int_0^t \|\partial_x\nabla\psi\|_{L^2}\|\partial_x^2\nabla^2\psi\|_{L^2}\,d\tau\\
\leq&\ C\left(\int_0^t \|\partial_x\nabla\psi\|_{L^2}^2\,d\tau\right)^{\frac12}\left(\int_0^t \|\partial_x^2\nabla\psi\|_{L^2}^2\,d\tau\right)^{\frac12}\\
\leq&\ C\varepsilon^2 \int_0^t (1+\tau)^{-(s+2)}\, d\tau\\
\leq&\ C\varepsilon^2.
\end{split}
\end{equation}
We also have
\begin{equation}\notag
\begin{split}
\|\psi\|_{B^{\frac12,\frac32}}&\leq \sum_{k\leq j}\lambda_j^{\frac32}\lambda_k^{\frac12}\|\Delta_j\Delta_k^h\psi\|_{L^2}\\
&\leq \sum_{j}\lambda_j^{\frac32}\|\Delta_j\partial_x\psi\|_{L^2}^{\frac12}\|\Delta_j\psi\|_{L^2}^{\frac12}\\
&\leq C\|\nabla\partial_x\psi\|_{L^2}^{\frac12}\left(\|\nabla\psi\|_{L^2} +\|\nabla^3\psi\|_{L^2}\right)^{\frac12}
\end{split}
\end{equation}
and hence
\begin{equation}\notag
\begin{split}
&\int_0^t\|\psi(\tau)\|_{B^{\frac12,\frac32}}^4 \, d\tau\\
\leq&\ C\int_0^t \|\nabla\partial_x\psi\|_{L^2}^{2}\left(\|\nabla\psi\|_{L^2} +\|\nabla^3\psi\|_{L^2}\right)^{2}\,d\tau\\
\leq&\ C\varepsilon^2 \int_0^t (1+\tau)^{-\frac12(s+1)\cdot 2}\, d\tau\\
\leq&\ C\varepsilon^2.
\end{split}
\end{equation}
Combining the analysis above gives
\[\int_0^t g_3^2(\tau)\, d\tau\leq C\varepsilon^2, \ \ t\in[0,T). \]

\cbdu

\medskip

\subsection{Finishing the proof of Theorem \ref{thm-main}}
Note the decay estimates in (ii) of Theorem \ref{thm-main} are already obtained in Lemma \ref{le-high-decay}. The uniqueness is trivial once we show existence in the Sobolev space as in (i). We are only left to complete the proof of global in time existence. First of all, the local well-posedness in $H^{s,-s}\cap H^{-s, 6} \cap H^6$ can be obtained on a short time interval $[0,T)$ through standard method. Also, for small enough initial data and small time $T>0$, the solution satisfies (\ref{assu-Es}).
We choose $\varepsilon>0$ small enough such that 
\[g_1(t)\leq C\varepsilon\leq c_1\]
from Lemma \ref{le-gg}, where $c_1$ is the constant in Lemma \ref{le-energy3}. It then follows from Lemma \ref{le-energy3} and Lemma \ref{le-gg} that
\begin{equation}\notag
\begin{split}
&E_{0,5}(t)+E_{s,-s}(t)+E_{s,5}(t)\\
\leq&\ C\left( E_{0,5}(t)+E_{s,-s}(t)+E_{s,5}(0)\right)\\
&+C\int_0^t\left(g_2(\tau)+g_3^2(\tau)\right)\left(E_{0,5}(\tau)+E_{s,-s}(\tau)+E_{s,5}(\tau) \right)\,d\tau\\
\leq&\ C\left( E_{0,5}(t)+E_{s,-s}(t)+E_{s,5}(0)\right)\\
&+C\varepsilon\sup_{\tau\in[0,t)}\left(E_{0,5}(\tau)+E_{s,-s}(\tau)+E_{s,5}(\tau) \right).
\end{split}
\end{equation}
Taking $\varepsilon$ small enough such that $C\varepsilon<\frac12$, we infer
\begin{equation}\notag
E_{0,5}(t)+E_{s,-s}(t)+E_{s,5}(t)\leq 2C\left( E_{0,5}(t)+E_{s,-s}(t)+E_{s,5}(0)\right)\leq 2C\varepsilon^2\leq \varepsilon.
\end{equation}
Thus the global existence in Theorem \ref{thm-main} follows from a continuous argument.

\cbdu


\section{Auxiliary estimates to obtain Lemma \ref{le-energy3}}
\label{sec-aux}

In this section, we provide detailed estimates for the integrals in (\ref{energy-high-1}).

\begin{Lemma}\label{le-high-aux12}
Let $0<s<\frac12$. We have
\begin{equation}\notag
\begin{split}
&\sum_{j,k}\lambda_j^{2s_1}\lambda_k^{-2s}|\tilde I_1+\tilde I_2+\tilde I_3+\tilde I_4|\\
\leq &\ \frac{1}{16} c\left(\|\nabla b\|_{H^{-s, s_1+1}}^2+\|\nabla \psi_x\|_{H^{-s, s_1}}^2\right)\\
&+ C\left(\|b\|_{B^{\frac32,\frac12}}+ \|\psi\|_{B^{\frac32,\frac12}}^2+\|b\|_{B^{\frac12,\frac32}}^2 +\|\psi\|_{B^{\frac12,\frac32}}^4\right)\\
&\cdot\left(\|\nabla \psi\|_{H^{-s, s_1}}^2+\|\nabla \psi\|_{H^{-s, s_1+1}}^2\right)
\end{split}
\end{equation}
for a constant $C>0$.
\end{Lemma}
\pf
Applying first integration by parts and exploring cancellations yields
\begin{equation}\notag
\begin{split}
\tilde I_1+\tilde I_2=&-\int_{\mathbb R^2}\Delta_j\Delta_k^hb_x \Delta_j\Delta_k^h(\psi_y\Delta \psi) \, dxdy-\int_{\mathbb R^2}\Delta_j\Delta_k^hb \Delta_j\Delta_k^h(\psi_{xy}\Delta \psi) \, dxdy\\
&+\int_{\mathbb R^2}\Delta_j\Delta_k^hb_y \Delta_j\Delta_k^h(\psi_x\Delta \psi) \, dxdy+\int_{\mathbb R^2}\Delta_j\Delta_k^hb \Delta_j\Delta_k^h(\psi_{xy}\Delta \psi) \, dxdy\\
=&-\int_{\mathbb R^2}\Delta_j\Delta_k^hb_x \Delta_j\Delta_k^h(\psi_y\Delta \psi) \, dxdy+\int_{\mathbb R^2}\Delta_j\Delta_k^hb_y \Delta_j\Delta_k^h(\psi_x\Delta \psi) \, dxdy\\
=& \int_{\mathbb R^2}\Delta_j\Delta_k^h\nabla b_x \Delta_j\Delta_k^h(\psi_y\nabla \psi) \, dxdy+\int_{\mathbb R^2}\Delta_j\Delta_k^hb_x \Delta_j\Delta_k^h(\nabla \psi_y\nabla \psi) \, dxdy\\
&-\int_{\mathbb R^2}\Delta_j\Delta_k^h\nabla b_y \Delta_j\Delta_k^h(\psi_x\nabla \psi) \, dxdy-\int_{\mathbb R^2}\Delta_j\Delta_k^hb_y \Delta_j\Delta_k^h(\nabla\psi_x\nabla \psi) \, dxdy\\
=&\int_{\mathbb R^2}\Delta_j\Delta_k^h\nabla b_x \Delta_j\Delta_k^h(\psi_y\nabla \psi) \, dxdy-\int_{\mathbb R^2}\Delta_j\Delta_k^h\nabla b_y \Delta_j\Delta_k^h(\psi_x\nabla \psi) \, dxdy
\end{split}
\end{equation}
since 
\begin{equation}\notag
\begin{split}
&\int_{\mathbb R^2}\Delta_j\Delta_k^hb_x \Delta_j\Delta_k^h(\nabla \psi_y\nabla \psi) \, dxdy-\int_{\mathbb R^2}\Delta_j\Delta_k^hb_y \Delta_j\Delta_k^h(\nabla\psi_x\nabla \psi) \, dxdy\\
=&-\frac12\int_{\mathbb R^2}\Delta_j\Delta_k^hb_{xy} \Delta_j\Delta_k^h(\nabla \psi\nabla \psi) \, dxdy+\frac12\int_{\mathbb R^2}\Delta_j\Delta_k^hb_{xy} \Delta_j\Delta_k^h(\nabla \psi\nabla \psi) \, dxdy\\
=&\ 0.
\end{split}
\end{equation}
We further explore cancellations by applying commutator and integration by parts to the integrals in $\tilde I_1+\tilde I_2$,
\begin{equation}\label{commu-12}
\begin{split}
&\tilde I_1+\tilde I_2\\
=&\int_{\mathbb R^2}\Delta_j\Delta_k^h\nabla b_x \psi_y\Delta_j\Delta_k^h\nabla \psi \, dxdy+\int_{\mathbb R^2}\Delta_j\Delta_k^h\nabla b_x [\Delta_j\Delta_k^h, \psi_y\nabla] \psi \, dxdy\\
&-\int_{\mathbb R^2}\Delta_j\Delta_k^h\nabla b_y \psi_x\Delta_j\Delta_k^h\nabla \psi \, dxdy-\int_{\mathbb R^2}\Delta_j\Delta_k^h\nabla b_y[\Delta_j\Delta_k^h, \psi_x\nabla] \psi \, dxdy\\
=&\int_{\mathbb R^2}\Delta_j\Delta_k^h\nabla b_x \psi_y\Delta_j\Delta_k^h\nabla \psi \, dxdy-\int_{\mathbb R^2}\Delta_j\Delta_k^h\nabla b [\Delta_j\Delta_k^h, \psi_y\nabla] \psi_x \, dxdy\\
&-\int_{\mathbb R^2}\Delta_j\Delta_k^h\nabla b_y \psi_x\Delta_j\Delta_k^h\nabla \psi \, dxdy+\int_{\mathbb R^2}\Delta_j\Delta_k^h\nabla b[\Delta_j\Delta_k^h, \psi_x\nabla] \psi_y \, dxdy
\end{split}
\end{equation}
since
\begin{equation}\notag
\begin{split}
&\int_{\mathbb R^2}\Delta_j\Delta_k^h\nabla b_x [\Delta_j\Delta_k^h, \psi_y\nabla] \psi \, dxdy-\int_{\mathbb R^2}\Delta_j\Delta_k^h\nabla b_y[\Delta_j\Delta_k^h, \psi_x\nabla] \psi \, dxdy\\
=&-\int_{\mathbb R^2}\Delta_j\Delta_k^h\nabla b [\Delta_j\Delta_k^h, \psi_y\nabla] \psi_x \, dxdy-\int_{\mathbb R^2}\Delta_j\Delta_k^h\nabla b [\Delta_j\Delta_k^h, \psi_{xy}\nabla] \psi \, dxdy\\
&+\int_{\mathbb R^2}\Delta_j\Delta_k^h\nabla b[\Delta_j\Delta_k^h, \psi_x\nabla] \psi_y \, dxdy+\int_{\mathbb R^2}\Delta_j\Delta_k^h\nabla b[\Delta_j\Delta_k^h, \psi_{yx}\nabla] \psi \, dxdy\\
=&-\int_{\mathbb R^2}\Delta_j\Delta_k^h\nabla b [\Delta_j\Delta_k^h, \psi_y\nabla] \psi_x \, dxdy+\int_{\mathbb R^2}\Delta_j\Delta_k^h\nabla b[\Delta_j\Delta_k^h, \psi_x\nabla] \psi_y \, dxdy.
\end{split}
\end{equation}

Now we rearrange the integrals in $\tilde I_3+\tilde I_4$,
\begin{equation}\notag
\begin{split}
\tilde I_3+\tilde I_4
=&\int_{\mathbb R^2}\Delta_j\Delta_k^h\nabla\psi \Delta_j\Delta_k^h\nabla(b_y \psi_x) \, dxdy-\int_{\mathbb R^2}\Delta_j\Delta_k^h\nabla\psi \Delta_j\Delta_k^h\nabla(b_x \psi_y) \, dxdy\\
=&\int_{\mathbb R^2}\Delta_j\Delta_k^h\nabla\psi \Delta_j\Delta_k^h(\nabla b_y \psi_x) \, dxdy+\int_{\mathbb R^2}\Delta_j\Delta_k^h\nabla\psi \Delta_j\Delta_k^h(b_y \nabla\psi_x) \, dxdy\\
&-\int_{\mathbb R^2}\Delta_j\Delta_k^h\nabla\psi \Delta_j\Delta_k^h(\nabla b_x \psi_y) \, dxdy
-\int_{\mathbb R^2}\Delta_j\Delta_k^h\nabla \psi \Delta_j\Delta_k^h(b_x\nabla\psi_y) \, dxdy.
\end{split}
\end{equation}
Applying commutator to the first and third integrals we have
\begin{equation}\notag
\begin{split}
&\int_{\mathbb R^2}\Delta_j\Delta_k^h\nabla\psi \Delta_j\Delta_k^h(\nabla b_y \psi_x) \, dxdy-\int_{\mathbb R^2}\Delta_j\Delta_k^h\nabla\psi \Delta_j\Delta_k^h(\nabla b_x \psi_y) \, dxdy\\
=&\int_{\mathbb R^2}\Delta_j\Delta_k^h\nabla\psi  \psi_x\Delta_j\Delta_k^h\nabla b_y \, dxdy+\int_{\mathbb R^2}\Delta_j\Delta_k^h\nabla\psi  [\Delta_j\Delta_k^h,\psi_x\nabla] b_y \, dxdy\\
&-\int_{\mathbb R^2}\Delta_j\Delta_k^h\nabla\psi \psi_y\Delta_j\Delta_k^h\nabla b_x \, dxdy
-\int_{\mathbb R^2}\Delta_j\Delta_k^h\nabla\psi [\Delta_j\Delta_k^h, \psi_y\nabla] b_x \, dxdy.
\end{split}
\end{equation}
Applying commutator to the second and forth integrals yields
\begin{equation}\notag
\begin{split}
&\int_{\mathbb R^2}\Delta_j\Delta_k^h\nabla\psi \Delta_j\Delta_k^h( b_y \nabla\psi_x) \, dxdy-\int_{\mathbb R^2}\Delta_j\Delta_k^h\nabla\psi \Delta_j\Delta_k^h( b_x \nabla\psi_y) \, dxdy\\
=&\int_{\mathbb R^2}\Delta_j\Delta_k^h\nabla\psi  b_y\Delta_j\Delta_k^h\nabla \psi_x \, dxdy+\int_{\mathbb R^2}\Delta_j\Delta_k^h\nabla\psi  [\Delta_j\Delta_k^h,b_y\nabla] \psi_x \, dxdy\\
&-\int_{\mathbb R^2}\Delta_j\Delta_k^h\nabla\psi b_x\Delta_j\Delta_k^h\nabla \psi_y \, dxdy
-\int_{\mathbb R^2}\Delta_j\Delta_k^h\nabla\psi [\Delta_j\Delta_k^h, b_x\nabla] \psi_y \, dxdy\\
=&\int_{\mathbb R^2}\Delta_j\Delta_k^h\nabla\psi  [\Delta_j\Delta_k^h,b_y\nabla] \psi_x \, dxdy-\int_{\mathbb R^2}\Delta_j\Delta_k^h\nabla\psi [\Delta_j\Delta_k^h, b_x\nabla] \psi_y \, dxdy
\end{split}
\end{equation}
since we have the cancellation
\begin{equation}\notag
\begin{split}
&\int_{\mathbb R^2}\Delta_j\Delta_k^h\nabla\psi  b_y\Delta_j\Delta_k^h\nabla \psi_x \, dxdy-\int_{\mathbb R^2}\Delta_j\Delta_k^h\nabla\psi b_x\Delta_j\Delta_k^h\nabla \psi_y \, dxdy\\
=&-\frac12 \int_{\mathbb R^2}  b_{xy}\left(\Delta_j\Delta_k^h\nabla \psi\right)^2 \, dxdy+\frac12 \int_{\mathbb R^2}  b_{yx}\left(\Delta_j\Delta_k^h\nabla \psi\right)^2 \, dxdy\\
=&\ 0.
\end{split}
\end{equation}
Therefore we infer
\begin{equation}\label{commu-34}
\begin{split}
&\tilde I_3+\tilde I_4\\
=&\int_{\mathbb R^2}\Delta_j\Delta_k^h\nabla\psi  \psi_x\Delta_j\Delta_k^h\nabla b_y \, dxdy+\int_{\mathbb R^2}\Delta_j\Delta_k^h\nabla\psi  [\Delta_j\Delta_k^h,\psi_x\nabla] b_y \, dxdy\\
&-\int_{\mathbb R^2}\Delta_j\Delta_k^h\nabla\psi \psi_y\Delta_j\Delta_k^h\nabla b_x \, dxdy
-\int_{\mathbb R^2}\Delta_j\Delta_k^h\nabla\psi [\Delta_j\Delta_k^h, \psi_y\nabla] b_x \, dxdy\\
&+\int_{\mathbb R^2}\Delta_j\Delta_k^h\nabla\psi  [\Delta_j\Delta_k^h,b_y\nabla] \psi_x \, dxdy-\int_{\mathbb R^2}\Delta_j\Delta_k^h\nabla\psi [\Delta_j\Delta_k^h, b_x\nabla] \psi_y \, dxdy.
\end{split}
\end{equation}
Observing cancellations among $\tilde I_1+\tilde I_2$ and $\tilde I_3+\tilde I_4$, it follows from (\ref{commu-12}) and (\ref{commu-34}) 
\begin{equation}\label{commu-1234}
\begin{split}
&\tilde I_1+\tilde I_2+\tilde I_3+\tilde I_4\\
=&-\int_{\mathbb R^2}\Delta_j\Delta_k^h\nabla b [\Delta_j\Delta_k^h, \psi_y\nabla] \psi_x \, dxdy
+\int_{\mathbb R^2}\Delta_j\Delta_k^h\nabla b[\Delta_j\Delta_k^h, \psi_x\nabla] \psi_y \, dxdy\\
&+\int_{\mathbb R^2}\Delta_j\Delta_k^h\nabla\psi  [\Delta_j\Delta_k^h,\psi_x\nabla] b_y \, dxdy
-\int_{\mathbb R^2}\Delta_j\Delta_k^h\nabla\psi [\Delta_j\Delta_k^h, \psi_y\nabla] b_x \, dxdy\\
&+\int_{\mathbb R^2}\Delta_j\Delta_k^h\nabla\psi  [\Delta_j\Delta_k^h,b_y\nabla] \psi_x \, dxdy
-\int_{\mathbb R^2}\Delta_j\Delta_k^h\nabla\psi [\Delta_j\Delta_k^h, b_x\nabla] \psi_y \, dxdy.
\end{split}
\end{equation} 
We estimate the integrals from (\ref{commu-1234}) in the following. Applying H\"older's inequality and the commutator estimate gives
\begin{equation}\notag
\begin{split}
&\sum_{j,k}\lambda_j^{2s_1}\lambda_k^{-2s}\left|\int_{\mathbb R^2}\Delta_j\Delta_k^h\nabla b [\Delta_j\Delta_k^h, \psi_y\nabla] \psi_x \, dxdy\right|\\
\leq&\sum_{j,k}\lambda_j^{2s_1}\lambda_k^{-2s} \|\Delta_j\Delta_k^h\nabla b\|_{L^2}\|\Delta_j\Delta_k^h\nabla \psi_y\|_{L^2}\|\psi_x\|_{L^\infty}\\
\leq&\ \|\psi_x\|_{B^{\frac12,\frac12}}\sum_{j,k} \|\Delta_j\Delta_k^h|D_x|^{-s}D^{s_1+1}\nabla b\|_{L^2}\|\Delta_j\Delta_k^h|D_x|^{-s}D^{s_1-1}\nabla \psi_y\|_{L^2}\\
\leq&\ \|\psi_x\|_{B^{\frac12,\frac12}}\|\nabla b\|_{H^{-s,s_1+1}}\|\psi_y\|_{H^{-s, s_1}}\\
\leq& \  \frac{1}{256} c\|\nabla b\|_{H^{-s, s_1+1}}^2+ C\|\psi_x\|_{B^{\frac12,\frac12}}^2\|\nabla \psi\|_{H^{-s, s_1}}^2,
\end{split}
\end{equation}
\begin{equation}\notag
\begin{split}
&\sum_{j,k}\lambda_j^{2s_1}\lambda_k^{-2s}\left|\int_{\mathbb R^2}\Delta_j\Delta_k^h\nabla b [\Delta_j\Delta_k^h, \psi_x\nabla] \psi_y \, dxdy\right|\\
\leq&\sum_{j,k}\lambda_j^{2s_1}\lambda_k^{-2s} \|\Delta_j\Delta_k^h\nabla b\|_{L^2}\|\Delta_j\Delta_k^h\nabla \psi_x\|_{L^2}\|\psi_y\|_{L^\infty}\\
\leq&\ \|\psi_y\|_{B^{\frac12,\frac12}}\sum_{j,k} \|\Delta_j\Delta_k^h|D_x|^{-s}D^{s_1}\nabla b\|_{L^2}\|\Delta_j\Delta_k^h|D_x|^{-s}D^{s_1}\nabla \psi_x\|_{L^2}\\
\leq&\ \|\psi_y\|_{B^{\frac12,\frac12}}\| b\|_{H^{-s,s_1+1}}\|\nabla\psi_x\|_{H^{-s, s_1}}\\
\leq& \  \frac{1}{256}  c\left(\|b\|_{H^{-s, s_1+1}}^2+\|\nabla \psi_x\|_{H^{-s, s_1}}^2\right)+ C\|\psi_y\|_{B^{\frac12,\frac12}}^4\|\nabla \psi_x\|_{H^{-s, s_1}}^2,
\end{split}
\end{equation}
\begin{equation}\notag
\begin{split}
&\sum_{j,k}\lambda_j^{2s_1}\lambda_k^{-2s}\left|\int_{\mathbb R^2}\Delta_j\Delta_k^h\nabla \psi [\Delta_j\Delta_k^h, \psi_x\nabla] b_y \, dxdy\right|\\
\leq&\sum_{j,k}\lambda_j^{2s_1}\lambda_k^{-2s} \|\Delta_j\Delta_k^h\nabla \psi\|_{L^2}\|\Delta_j\Delta_k^h\nabla \psi_x\|_{L^2}\|b_y\|_{L^\infty}\\
\leq&\ \|b_y\|_{B^{\frac12,\frac12}}\sum_{j,k} \|\Delta_j\Delta_k^h|D_x|^{-s}D^{s_1}\nabla \psi\|_{L^2}\|\Delta_j\Delta_k^h|D_x|^{-s}D^{s_1}\nabla \psi_x\|_{L^2}\\
\leq&\ \|b_y\|_{B^{\frac12,\frac12}}\| \nabla\psi\|_{H^{-s,s_1}}\|\nabla\psi_x\|_{H^{-s, s_1}}\\
\leq& \  \frac{1}{256}  c\|\nabla \psi_x\|_{H^{-s, s_1}}^2+ C\|b_y\|_{B^{\frac12,\frac12}}^2\|\nabla \psi\|_{H^{-s, s_1}}^2,
\end{split}
\end{equation}
\begin{equation}\notag
\begin{split}
&\sum_{j,k}\lambda_j^{2s_1}\lambda_k^{-2s}\left|\int_{\mathbb R^2}\Delta_j\Delta_k^h\nabla \psi [\Delta_j\Delta_k^h, \psi_y\nabla] b_x \, dxdy\right|\\
\leq&\sum_{j,k}\lambda_j^{2s_1}\lambda_k^{-2s} \|\Delta_j\Delta_k^h\nabla \psi\|_{L^2}\|\Delta_j\Delta_k^h\nabla \psi_y\|_{L^2}\|b_x\|_{L^\infty}\\
\leq&\ \|b_x\|_{B^{\frac12,\frac12}}\sum_{j,k} \|\Delta_j\Delta_k^h|D_x|^{-s}D^{s_1}\nabla \psi\|_{L^2}\|\Delta_j\Delta_k^h|D_x|^{-s}D^{s_1}\nabla \psi_y\|_{L^2}\\
\leq&\ \|b_x\|_{B^{\frac12,\frac12}}\| \nabla\psi\|_{H^{-s,s_1}}\|\nabla\psi_y\|_{H^{-s, s_1}}\\
\leq& \  \|b_x\|_{B^{\frac12,\frac12}}\left(\|\nabla \psi\|_{H^{-s, s_1}}^2+\|\nabla \psi\|_{H^{-s, s_1+1}}^2\right),
\end{split}
\end{equation}
\begin{equation}\notag
\begin{split}
&\sum_{j,k}\lambda_j^{2s_1}\lambda_k^{-2s}\left|\int_{\mathbb R^2}\Delta_j\Delta_k^h\nabla \psi [\Delta_j\Delta_k^h, b_y\nabla] \psi_x \, dxdy\right|\\
\leq&\sum_{j,k}\lambda_j^{2s_1}\lambda_k^{-2s} \|\Delta_j\Delta_k^h\nabla \psi\|_{L^2}\|\Delta_j\Delta_k^h\nabla b_y\|_{L^2}\|\psi_x\|_{L^\infty}\\
\leq&\ \|\psi_x\|_{B^{\frac12,\frac12}}\sum_{j,k} \|\Delta_j\Delta_k^h|D_x|^{-s}D^{s_1}\nabla \psi\|_{L^2}\|\Delta_j\Delta_k^h|D_x|^{-s}D^{s_1}\nabla b_y\|_{L^2}\\
\leq&\ \|\psi_x\|_{B^{\frac12,\frac12}}\| \nabla\psi\|_{H^{-s,s_1}}\|\nabla b_y\|_{H^{-s, s_1}}\\
\leq& \ \frac{1}{256}  c\|\nabla b\|_{H^{-s, s_1+1}}^2+C \|\psi_x\|_{B^{\frac12,\frac12}}^2\|\nabla \psi\|_{H^{-s, s_1}}^2.
\end{split}
\end{equation}
To estimate the last integral from (\ref{commu-1234}), we further manipulate it through integration by parts as 
\begin{equation}\notag
\begin{split}
&-\int_{\mathbb R^2}\Delta_j\Delta_k^h\nabla\psi [\Delta_j\Delta_k^h, b_x\nabla] \psi_y \, dxdy\\
=&\int_{\mathbb R^2}\Delta_j\Delta_k^h\nabla\psi_x [\Delta_j\Delta_k^h, b\nabla] \psi_y \, dxdy
+\int_{\mathbb R^2}\Delta_j\Delta_k^h\nabla\psi [\Delta_j\Delta_k^h, b\nabla] \psi_{xy} \, dxdy\\
=&\int_{\mathbb R^2}\Delta_j\Delta_k^h\nabla\psi_x [\Delta_j\Delta_k^h, b\nabla] \psi_y \, dxdy
-\int_{\mathbb R^2}\Delta_j\Delta_k^h\nabla\psi_y [\Delta_j\Delta_k^h, b\nabla] \psi_{x} \, dxdy\\
&-\int_{\mathbb R^2}\Delta_j\Delta_k^h\nabla\psi [\Delta_j\Delta_k^h, b_y\nabla] \psi_{x} \, dxdy.
\end{split}
\end{equation}
Note the integral from the last term 
\[\lambda_j^{2s_1}\lambda_k^{-2s}\int_{\mathbb R^2}\Delta_j\Delta_k^h\nabla\psi [\Delta_j\Delta_k^h, b_y\nabla] \psi_{x} \, dxdy\]
is already estimated. The other two terms are handled as
\begin{equation}\notag
\begin{split}
&\sum_{j,k}\lambda_j^{2s_1}\lambda_k^{-2s}\left|\int_{\mathbb R^2}\Delta_j\Delta_k^h\nabla \psi_x [\Delta_j\Delta_k^h, b\nabla] \psi_y \, dxdy\right|\\
\leq&\sum_{j,k}\lambda_j^{2s_1}\lambda_k^{-2s} \|\Delta_j\Delta_k^h\nabla \psi_x\|_{L^2}\|\Delta_j\Delta_k^h\nabla b\|_{L^2}\|\psi_y\|_{L^\infty}\\
\leq&\ \|\psi_y\|_{B^{\frac12,\frac12}}\sum_{j,k} \|\Delta_j\Delta_k^h|D_x|^{-s}D^{s_1}\nabla \psi_x\|_{L^2}\|\Delta_j\Delta_k^h|D_x|^{-s}D^{s_1}\nabla b\|_{L^2}\\
\leq&\ \|\psi_y\|_{B^{\frac12,\frac12}}\| \nabla\psi_x\|_{H^{-s,s_1}}\|\nabla b\|_{H^{-s, s_1}}\\
\leq& \ \frac{1}{256}  c\left(\|\nabla b\|_{H^{-s, s_1}}^2+\|\nabla \psi_x\|_{H^{-s, s_1}}^2\right)+C \|\psi_y\|_{B^{\frac12,\frac12}}^4\|\nabla \psi_x\|_{H^{-s, s_1}}^2,
\end{split}
\end{equation}
\begin{equation}\notag
\begin{split}
&\sum_{j,k}\lambda_j^{2s_1}\lambda_k^{-2s}\left|\int_{\mathbb R^2}\Delta_j\Delta_k^h\nabla \psi_y [\Delta_j\Delta_k^h, b\nabla] \psi_x \, dxdy\right|\\
\leq&\sum_{j,k}\lambda_j^{2s_1}\lambda_k^{-2s} \|\Delta_j\Delta_k^h\nabla \psi_y\|_{L^2}\|\Delta_j\Delta_k^h\nabla b\|_{L^2}\|\psi_x\|_{L^\infty}\\
\leq&\ \|\psi_x\|_{B^{\frac12,\frac12}}\sum_{j,k} \|\Delta_j\Delta_k^h|D_x|^{-s}D^{s_1-1}\nabla \psi_y\|_{L^2}\|\Delta_j\Delta_k^h|D_x|^{-s}D^{s_1+1}\nabla b\|_{L^2}\\
\leq&\ \|\psi_x\|_{B^{\frac12,\frac12}}\| \nabla\psi\|_{H^{-s,s_1}}\|\nabla b\|_{H^{-s, s_1+1}}\\
\leq& \ \frac{1}{256}  c\|\nabla b\|_{H^{-s, s_1+1}}^2+C \|\psi_x\|_{B^{\frac12,\frac12}}^2\|\nabla \psi\|_{H^{-s, s_1}}^2.
\end{split}
\end{equation}

The lemma follows from summarizing the estimates above.

\cbdu

\medskip

The highest nonlinear and most difficult terms are in $\tilde I_5+\tilde I_6$ and $\tilde I_7+\tilde I_8$. We need to fully explore the cancellations among them and hence handle them together as well.

\begin{Lemma}\label{le-high-aux5678}
Let $0<s<\frac12$ and $s_1>-\frac12$. The estimate
\begin{equation}\notag
\begin{split}
&\sum_{j,k}\lambda_j^{2s_1}\lambda_k^{-2s}|\tilde I_5+\tilde I_6+\tilde I_7+\tilde I_8|\\
\leq &\ \frac{1}{16} c\left( \|\nabla b\|_{H^{-s, s_1+1}}^2+\|\nabla \psi_x\|_{H^{-s, s_1}}^2\right)+C  \|\nabla b\|_{B^{\frac12,\frac32}}\|\nabla \psi\|_{H^{-s, s_1+1}}^2\\
&+ C\left(\|\nabla\psi\|_{B^{\frac32,\frac12}}^2+\|\nabla\psi\|_{B^{\frac12,\frac32}}^4\right)\|\nabla \psi\|_{H^{-s, s_1+1}}^2
\end{split}
\end{equation}
holds for a constant $C>0$.
\end{Lemma}
\pf
As usual, we start with integration by parts
\begin{equation}\notag
\begin{split}
\tilde I_5+\tilde I_6=&-\int_{\mathbb R^2}\Delta_j\Delta_k^h\nabla b_x \Delta_j\Delta_k^h\nabla(\psi_y\Delta \psi) \, dxdy-\int_{\mathbb R^2}\Delta_j\Delta_k^h\nabla b \Delta_j\Delta_k^h\nabla(\psi_{xy}\Delta \psi) \, dxdy\\
&+\int_{\mathbb R^2}\Delta_j\Delta_k^h\nabla b_y \Delta_j\Delta_k^h\nabla(\psi_x\Delta \psi) \, dxdy+\int_{\mathbb R^2}\Delta_j\Delta_k^h\nabla b \Delta_j\Delta_k^h\nabla(\psi_{xy}\Delta \psi) \, dxdy\\
=&-\int_{\mathbb R^2}\Delta_j\Delta_k^h\nabla b_x \Delta_j\Delta_k^h\nabla (\psi_y\Delta \psi) \, dxdy+\int_{\mathbb R^2}\Delta_j\Delta_k^h\nabla b_y \Delta_j\Delta_k^h\nabla(\psi_x\Delta \psi) \, dxdy\\
=& \int_{\mathbb R^2}\Delta_j\Delta_k^h\nabla^2 b_x \Delta_j\Delta_k^h\nabla(\psi_y\nabla \psi) \, dxdy+\int_{\mathbb R^2}\Delta_j\Delta_k^h\nabla b_x \Delta_j\Delta_k^h\nabla (\nabla \psi_y\nabla \psi) \, dxdy\\
&-\int_{\mathbb R^2}\Delta_j\Delta_k^h\nabla^2 b_y \Delta_j\Delta_k^h\nabla(\psi_x\nabla \psi) \, dxdy-\int_{\mathbb R^2}\Delta_j\Delta_k^h\nabla b_y \Delta_j\Delta_k^h\nabla(\nabla\psi_x\nabla \psi) \, dxdy\\
=&\int_{\mathbb R^2}\Delta_j\Delta_k^h\nabla^2 b_x \Delta_j\Delta_k^h\nabla(\psi_y\nabla \psi) \, dxdy-\int_{\mathbb R^2}\Delta_j\Delta_k^h\nabla^2 b_y \Delta_j\Delta_k^h\nabla(\psi_x\nabla \psi) \, dxdy
\end{split}
\end{equation}
since 
\begin{equation}\notag
\begin{split}
&\int_{\mathbb R^2}\Delta_j\Delta_k^h\nabla b_x \Delta_j\Delta_k^h\nabla (\nabla \psi_y\nabla \psi) \, dxdy-\int_{\mathbb R^2}\Delta_j\Delta_k^h\nabla b_y \Delta_j\Delta_k^h\nabla (\nabla\psi_x\nabla \psi) \, dxdy\\
=&-\frac12\int_{\mathbb R^2}\Delta_j\Delta_k^h\nabla b_{xy} \Delta_j\Delta_k^h\nabla (\nabla \psi\nabla \psi) \, dxdy+\frac12\int_{\mathbb R^2}\Delta_j\Delta_k^h\nabla b_{xy} \Delta_j\Delta_k^h\nabla (\nabla \psi\nabla \psi) \, dxdy\\
=&\ 0.
\end{split}
\end{equation}
We further use integration by parts and obtain
\begin{equation}\notag
\begin{split}
\tilde I_5+\tilde I_6
=&\int_{\mathbb R^2}\Delta_j\Delta_k^h\nabla^2 b_x \Delta_j\Delta_k^h(\nabla\psi_y\nabla \psi) \, dxdy
+\int_{\mathbb R^2}\Delta_j\Delta_k^h\nabla^2 b_x \Delta_j\Delta_k^h(\psi_y\nabla^2 \psi) \, dxdy\\
&-\int_{\mathbb R^2}\Delta_j\Delta_k^h\nabla^2 b_y \Delta_j\Delta_k^h(\nabla\psi_x\nabla \psi) \, dxdy-\int_{\mathbb R^2}\Delta_j\Delta_k^h\nabla^2 b_y \Delta_j\Delta_k^h(\psi_x\nabla^2 \psi) \, dxdy\\
=&\int_{\mathbb R^2}\Delta_j\Delta_k^h\nabla^2 b_x \Delta_j\Delta_k^h(\psi_y\nabla^2 \psi) \, dxdy-
\int_{\mathbb R^2}\Delta_j\Delta_k^h\nabla^2 b_y \Delta_j\Delta_k^h(\psi_x\nabla^2 \psi) \, dxdy\\
\end{split}
\end{equation}
thanks to
\begin{equation}\notag
\begin{split}
&\int_{\mathbb R^2}\Delta_j\Delta_k^h\nabla^2 b_x \Delta_j\Delta_k^h (\nabla \psi_y\nabla \psi) \, dxdy-\int_{\mathbb R^2}\Delta_j\Delta_k^h\nabla^2 b_y \Delta_j\Delta_k^h (\nabla\psi_x\nabla \psi) \, dxdy\\
=&-\frac12\int_{\mathbb R^2}\Delta_j\Delta_k^h\nabla^2 b_{xy} \Delta_j\Delta_k^h (\nabla \psi\nabla \psi) \, dxdy+\frac12\int_{\mathbb R^2}\Delta_j\Delta_k^h\nabla^2 b_{xy} \Delta_j\Delta_k^h (\nabla \psi\nabla \psi) \, dxdy\\
=&\ 0.
\end{split}
\end{equation}
Further applying commutators to the integrals in $\tilde I_5+\tilde I_6$ yields
\begin{equation}\label{high-56}
\begin{split}
\tilde I_5+\tilde I_6
=&\int_{\mathbb R^2}\Delta_j\Delta_k^h\nabla^2 b_x \psi_y\Delta_j\Delta_k^h\nabla^2 \psi \, dxdy\\
&+\int_{\mathbb R^2}\Delta_j\Delta_k^h\nabla^2 b_x [\Delta_j\Delta_k^h,\psi_y\nabla]\nabla \psi \, dxdy\\
&-\int_{\mathbb R^2}\Delta_j\Delta_k^h\nabla^2 b_y \psi_x\Delta_j\Delta_k^h\nabla^2 \psi \, dxdy\\
&-\int_{\mathbb R^2}\Delta_j\Delta_k^h\nabla^2 b_y [\Delta_j\Delta_k^h,\psi_x\nabla]\nabla \psi \, dxdy.
\end{split}
\end{equation}

Now we rearrange the integrals in $\tilde I_7+\tilde I_8$,
\begin{equation}\notag
\begin{split}
\tilde I_7+\tilde I_8=&\int_{\mathbb R^2}\Delta_j\Delta_k^h\nabla^2\psi \Delta_j\Delta_k^h(\nabla^2b_y \psi_x) \, dxdy+\int_{\mathbb R^2}\Delta_j\Delta_k^h\nabla^2\psi \Delta_j\Delta_k^h(b_y \nabla^2\psi_x) \, dxdy\\
&+2\int_{\mathbb R^2}\Delta_j\Delta_k^h\nabla^2\psi \Delta_j\Delta_k^h(\nabla b_y \nabla\psi_x) \, dxdy
-\int_{\mathbb R^2}\Delta_j\Delta_k^h\nabla^2\psi \Delta_j\Delta_k^h(\nabla^2b_x \psi_y) \, dxdy\\
&-\int_{\mathbb R^2}\Delta_j\Delta_k^h\nabla^2 \psi \Delta_j\Delta_k^h(b_x\nabla^2\psi_y) \, dxdy-2\int_{\mathbb R^2}\Delta_j\Delta_k^h\nabla^2\psi \Delta_j\Delta_k^h(\nabla b_x \nabla\psi_y) \, dxdy.
\end{split}
\end{equation}
Applying commutator to the first and forth integrals we have
\begin{equation}\notag
\begin{split}
&\int_{\mathbb R^2}\Delta_j\Delta_k^h\nabla^2\psi \Delta_j\Delta_k^h(\nabla^2b_y \psi_x) \, dxdy-\int_{\mathbb R^2}\Delta_j\Delta_k^h\nabla^2\psi \Delta_j\Delta_k^h(\nabla^2b_x \psi_y) \, dxdy\\
=&\int_{\mathbb R^2}\Delta_j\Delta_k^h\nabla^2\psi  \psi_x\Delta_j\Delta_k^h\nabla^2b_y \, dxdy+\int_{\mathbb R^2}\Delta_j\Delta_k^h\nabla^2\psi  [\Delta_j\Delta_k^h,\psi_x\nabla]\nabla b_y \, dxdy\\
&-\int_{\mathbb R^2}\Delta_j\Delta_k^h\nabla^2\psi \psi_y\Delta_j\Delta_k^h\nabla^2b_x \, dxdy
-\int_{\mathbb R^2}\Delta_j\Delta_k^h\nabla^2\psi [\Delta_j\Delta_k^h, \psi_y\nabla]\nabla b_x \, dxdy.
\end{split}
\end{equation}
Applying first commutator and then integration by parts to the second and fifth integrals in  $\tilde I_7+\tilde I_8$ gives
\begin{equation}\notag
\begin{split}
&\int_{\mathbb R^2}\Delta_j\Delta_k^h\nabla^2\psi \Delta_j\Delta_k^h(b_y \nabla^2\psi_x) \, dxdy-\int_{\mathbb R^2}\Delta_j\Delta_k^h\nabla^2 \psi \Delta_j\Delta_k^h(b_x\nabla^2\psi_y) \, dxdy\\
=&\int_{\mathbb R^2}\Delta_j\Delta_k^h\nabla^2\psi b_y\Delta_j\Delta_k^h\nabla^2\psi_x \, dxdy
+\int_{\mathbb R^2}\Delta_j\Delta_k^h\nabla^2\psi [\Delta_j\Delta_k^h, b_y\nabla]\nabla\psi_x \, dxdy\\
&-\int_{\mathbb R^2}\Delta_j\Delta_k^h\nabla^2 \psi b_x\Delta_j\Delta_k^h\nabla^2\psi_y \, dxdy
-\int_{\mathbb R^2}\Delta_j\Delta_k^h\nabla^2 \psi [\Delta_j\Delta_k^h, b_x\nabla]\nabla\psi_y \, dxdy\\
=&\int_{\mathbb R^2}\Delta_j\Delta_k^h\nabla^2\psi [\Delta_j\Delta_k^h, b_y\nabla]\nabla\psi_x \, dxdy
-\int_{\mathbb R^2}\Delta_j\Delta_k^h\nabla^2 \psi [\Delta_j\Delta_k^h, b_x\nabla]\nabla\psi_y \, dxdy
\end{split}
\end{equation}
since it follows from integration by parts that
\begin{equation}\notag
\begin{split}
&\int_{\mathbb R^2}\Delta_j\Delta_k^h\nabla^2\psi b_y\Delta_j\Delta_k^h\nabla^2\psi_x \, dxdy
-\int_{\mathbb R^2}\Delta_j\Delta_k^h\nabla^2 \psi b_x\Delta_j\Delta_k^h\nabla^2\psi_y \, dxdy\\
=&-\frac12 \int_{\mathbb R^2}b_{xy}(\Delta_j\Delta_k^h\nabla^2\psi)^2 \, dxdy+\frac12 \int_{\mathbb R^2}b_{xy}(\Delta_j\Delta_k^h\nabla^2\psi)^2 \, dxdy\\
=&\ 0.
\end{split}
\end{equation}
Using integration by parts to the third and sixth terms in $\tilde I_7+\tilde I_8$ we get
\begin{equation}\notag
\begin{split}
&2\int_{\mathbb R^2}\Delta_j\Delta_k^h\nabla^2\psi \Delta_j\Delta_k^h(\nabla b_y \nabla\psi_x) \, dxdy
-2\int_{\mathbb R^2}\Delta_j\Delta_k^h\nabla^2\psi \Delta_j\Delta_k^h(\nabla b_x \nabla\psi_y) \, dxdy\\
=&-2\int_{\mathbb R^2}\Delta_j\Delta_k^h\nabla^2\psi \Delta_j\Delta_k^h(\nabla b \nabla\psi_{xy}) \, dxdy
-2\int_{\mathbb R^2}\Delta_j\Delta_k^h\nabla^2\psi_y \Delta_j\Delta_k^h(\nabla b \nabla\psi_x) \, dxdy\\
&+2\int_{\mathbb R^2}\Delta_j\Delta_k^h\nabla^2\psi \Delta_j\Delta_k^h(\nabla b \nabla\psi_{xy}) \, dxdy
+2\int_{\mathbb R^2}\Delta_j\Delta_k^h\nabla^2\psi_x \Delta_j\Delta_k^h(\nabla b \nabla\psi_y) \, dxdy\\
=&-2\int_{\mathbb R^2}\Delta_j\Delta_k^h\nabla^2\psi_y \Delta_j\Delta_k^h(\nabla b \nabla\psi_x) \, dxdy
+2\int_{\mathbb R^2}\Delta_j\Delta_k^h\nabla^2\psi_x \Delta_j\Delta_k^h(\nabla b \nabla\psi_y) \, dxdy.
\end{split}
\end{equation}
Combining the manipulations above gives
\begin{equation}\label{high-78}
\begin{split}
\tilde I_7+\tilde I_8
=&\int_{\mathbb R^2}\Delta_j\Delta_k^h\nabla^2\psi  \psi_x\Delta_j\Delta_k^h\nabla^2b_y \, dxdy\\
&+\int_{\mathbb R^2}\Delta_j\Delta_k^h\nabla^2\psi  [\Delta_j\Delta_k^h,\psi_x\nabla]\nabla b_y \, dxdy\\
&-\int_{\mathbb R^2}\Delta_j\Delta_k^h\nabla^2\psi \psi_y\Delta_j\Delta_k^h\nabla^2b_x \, dxdy\\
&-\int_{\mathbb R^2}\Delta_j\Delta_k^h\nabla^2\psi [\Delta_j\Delta_k^h, \psi_y\nabla]\nabla b_x \, dxdy\\
&+\int_{\mathbb R^2}\Delta_j\Delta_k^h\nabla^2\psi [\Delta_j\Delta_k^h, b_y\nabla]\nabla\psi_x \, dxdy\\
&-\int_{\mathbb R^2}\Delta_j\Delta_k^h\nabla^2 \psi [\Delta_j\Delta_k^h, b_x\nabla]\nabla\psi_y \, dxdy\\
&-2\int_{\mathbb R^2}\Delta_j\Delta_k^h\nabla^2\psi_y \Delta_j\Delta_k^h(\nabla b \nabla\psi_x) \, dxdy\\
&+2\int_{\mathbb R^2}\Delta_j\Delta_k^h\nabla^2\psi_x \Delta_j\Delta_k^h(\nabla b \nabla\psi_y) \, dxdy.
\end{split}
\end{equation}

Observing the cancellations among (\ref{high-56}) and (\ref{high-78}), we obtain
\begin{equation}\label{high-5678}
\begin{split}
\tilde I_5+\tilde I_6+\tilde I_7+\tilde I_8
=&\int_{\mathbb R^2}\Delta_j\Delta_k^h\nabla^2 b_x [\Delta_j\Delta_k^h,\psi_y\nabla]\nabla \psi \, dxdy\\
&-\int_{\mathbb R^2}\Delta_j\Delta_k^h\nabla^2 b_y [\Delta_j\Delta_k^h,\psi_x\nabla]\nabla \psi \, dxdy\\
&+\int_{\mathbb R^2}\Delta_j\Delta_k^h\nabla^2\psi  [\Delta_j\Delta_k^h,\psi_x\nabla]\nabla b_y \, dxdy\\
&-\int_{\mathbb R^2}\Delta_j\Delta_k^h\nabla^2\psi [\Delta_j\Delta_k^h, \psi_y\nabla]\nabla b_x \, dxdy\\
&+\int_{\mathbb R^2}\Delta_j\Delta_k^h\nabla^2\psi [\Delta_j\Delta_k^h, b_y\nabla]\nabla\psi_x \, dxdy\\
&-\int_{\mathbb R^2}\Delta_j\Delta_k^h\nabla^2 \psi [\Delta_j\Delta_k^h, b_x\nabla]\nabla\psi_y \, dxdy\\
&-2\int_{\mathbb R^2}\Delta_j\Delta_k^h\nabla^2\psi_y \Delta_j\Delta_k^h(\nabla b \nabla\psi_x) \, dxdy\\
&+2\int_{\mathbb R^2}\Delta_j\Delta_k^h\nabla^2\psi_x \Delta_j\Delta_k^h(\nabla b \nabla\psi_y) \, dxdy.
\end{split}
\end{equation}

Further applying integration parts to the first two terms in (\ref{high-5678}) and exploring cancellations yields
\begin{equation}\notag
\begin{split}
&\int_{\mathbb R^2}\Delta_j\Delta_k^h\nabla^2 b_x [\Delta_j\Delta_k^h,\psi_y\nabla]\nabla \psi \, dxdy-\int_{\mathbb R^2}\Delta_j\Delta_k^h\nabla^2 b_y [\Delta_j\Delta_k^h,\psi_x\nabla]\nabla \psi \, dxdy\\
=&-\int_{\mathbb R^2}\Delta_j\Delta_k^h\nabla^2 b [\Delta_j\Delta_k^h,\psi_y\nabla]\nabla \psi_x \, dxdy+\int_{\mathbb R^2}\Delta_j\Delta_k^h\nabla^2 b [\Delta_j\Delta_k^h,\psi_x\nabla]\nabla \psi_y \, dxdy.
\end{split}
\end{equation}
Applying commutator estimate yields
\begin{equation}\notag
\begin{split}
&\lambda_j^{2s_1}\lambda_k^{-2s}\left|\int_{\mathbb R^2}\Delta_j\Delta_k^h\nabla^2 b [\Delta_j\Delta_k^h,\psi_y\nabla]\nabla \psi_x \, dxdy\right|\\
\leq&\ \lambda_j^{2s_1}\lambda_k^{-2s}\|\nabla\psi_x\|_{L^\infty}\|\Delta_j\Delta_k^h\nabla^2 b\|_{L^2}\|\Delta_j\Delta_k^h\nabla \psi_y\|_{L^2}\\
\leq&\ \|\nabla\psi_x\|_{B^{\frac12,\frac12}}\|\Delta_j\Delta_k^h|D_x|^{-s}D^{s_1}\nabla^2 b\|_{L^2}\|\Delta_j\Delta_k^h|D_x|^{-s}D^{s_1}\nabla \psi_y\|_{L^2}\\
\leq&\ \|\nabla\psi_x\|_{B^{\frac12,\frac12}}\|\nabla b\|_{H^{-s,s_1+1}}\|\nabla \psi\|_{H^{-s, s_1+1}}\\
\leq&\ \frac{1}{256} c\|\nabla b\|_{H^{-s,s_1+1}}^2+C\|\nabla\psi_x\|_{B^{\frac12,\frac12}}^2\|\nabla \psi\|_{H^{-s, s_1+1}}^2
\end{split}
\end{equation}
and similarly 
\begin{equation}\notag
\begin{split}
&\lambda_j^{2s_1}\lambda_k^{-2s}\left|\int_{\mathbb R^2}\Delta_j\Delta_k^h\nabla^2 b [\Delta_j\Delta_k^h,\psi_x\nabla]\nabla \psi_y \, dxdy\right|\\
\leq&\ \|\nabla\psi_y\|_{B^{\frac12,\frac12}}\|\nabla b\|_{H^{-s,s_1+1}}\|\nabla \psi_x\|_{H^{-s, s_1}}\\
\leq&\ \frac{1}{256}c\left(\|\nabla b\|_{H^{-s,s_1+1}}^2+\|\nabla \psi_x\|_{H^{-s, s_1}}^2\right)+C\|\nabla\psi_y\|_{B^{\frac12,\frac12}}^4\|\nabla \psi_x\|_{H^{-s, s_1}}^2.
\end{split}
\end{equation}
The third term of (\ref{high-5678}) is estimated as
\begin{equation}\notag
\begin{split}
&\lambda_j^{2s_1}\lambda_k^{-2s}\left|\int_{\mathbb R^2}\Delta_j\Delta_k^h\nabla^2 \psi [\Delta_j\Delta_k^h,\psi_x\nabla]\nabla b_y \, dxdy\right|\\
\leq&\ \lambda_j^{2s_1}\lambda_k^{-2s}\|\nabla b_y\|_{L^\infty}\|\Delta_j\Delta_k^h\nabla^2 \psi\|_{L^2}\|\Delta_j\Delta_k^h\nabla \psi_x\|_{L^2}\\
\leq&\ \|\nabla b\|_{B^{\frac12,\frac12}}\|\Delta_j\Delta_k^h|D_x|^{-s}D^{s_1}\nabla^2 \psi\|_{L^2}\|\Delta_j\Delta_k^h|D_x|^{-s}D^{s_1}\nabla \psi_x\|_{L^2}\\
\leq&\ \|\nabla b_y\|_{B^{\frac12,\frac12}}\|\nabla \psi\|_{H^{-s, s_1+1}}^2
\end{split}
\end{equation}
and the forth term analogously
\begin{equation}\notag
\begin{split}
&\lambda_j^{2s_1}\lambda_k^{-2s}\left|\int_{\mathbb R^2}\Delta_j\Delta_k^h\nabla^2 \psi [\Delta_j\Delta_k^h,\psi_y\nabla]\nabla b_x \, dxdy\right|\\
\leq&\ \|\nabla b_x\|_{B^{\frac12,\frac12}}\|\nabla \psi\|_{H^{-s, s_1+1}}^2.
\end{split}
\end{equation}
The fifth term of (\ref{high-5678}) is estimated as
\begin{equation}\notag
\begin{split}
&\lambda_j^{2s_1}\lambda_k^{-2s}\left|\int_{\mathbb R^2}\Delta_j\Delta_k^h\nabla^2 \psi [\Delta_j\Delta_k^h,b_y\nabla]\nabla \psi_x \, dxdy\right|\\
\leq&\ \lambda_j^{2s_1}\lambda_k^{-2s}\|\nabla \psi_x\|_{L^\infty}\|\Delta_j\Delta_k^h\nabla^2 \psi\|_{L^2}\|\Delta_j\Delta_k^h\nabla b_y\|_{L^2}\\
\leq&\ \|\nabla \psi_x\|_{B^{\frac12,\frac12}}\|\Delta_j\Delta_k^h|D_x|^{-s}D^{s_1}\nabla^2 \psi\|_{L^2}\|\Delta_j\Delta_k^h|D_x|^{-s}D^{s_1}\nabla b_y\|_{L^2}\\
\leq&\ \|\nabla \psi_x\|_{B^{\frac12,\frac12}}\|\nabla \psi\|_{H^{-s, s_1+1}}\|\nabla b\|_{H^{-s, s_1+1}}\\
\leq&\ \frac{1}{256}  c\|\nabla b\|_{H^{-s,s_1+1}}^2+C\|\nabla\psi_x\|_{B^{\frac12,\frac12}}^2\|\nabla \psi\|_{H^{-s, s_1+1}}^2
\end{split}
\end{equation}
and similarly
\begin{equation}\notag
\begin{split}
&\lambda_j^{2s_1}\lambda_k^{-2s}\left|\int_{\mathbb R^2}\Delta_j\Delta_k^h\nabla^2 \psi [\Delta_j\Delta_k^h,b_x\nabla]\nabla \psi_y \, dxdy\right|\\
\leq&\ \frac{1}{256}  c\left(\|\nabla b\|_{H^{-s,s_1+1}}^2+\|\nabla \psi_x\|_{H^{-s, s_1}}^2\right)+C\|\nabla\psi_y\|_{B^{\frac12,\frac12}}^4\|\nabla \psi_x\|_{H^{-s, s_1}}^2.
\end{split}
\end{equation}

We proceed to estimate the last two integrals in (\ref{high-5678}). First, we need to rearrange the integrals to further explore cancellations and the benefit of commutators,
\begin{equation}\label{high-5678-last}
\begin{split}
&-\int_{\mathbb R^2}\Delta_j\Delta_k^h\nabla^2\psi_y \Delta_j\Delta_k^h(\nabla b \nabla\psi_x) \, dxdy
+\int_{\mathbb R^2}\Delta_j\Delta_k^h\nabla^2\psi_x \Delta_j\Delta_k^h(\nabla b \nabla\psi_y) \, dxdy\\
=&\int_{\mathbb R^2}\Delta_j\Delta_k^h\nabla\psi_y \Delta_j\Delta_k^h(\nabla^2 b \nabla\psi_x) \, dxdy+\int_{\mathbb R^2}\Delta_j\Delta_k^h\nabla\psi_y \Delta_j\Delta_k^h(\nabla b \nabla^2\psi_x) \, dxdy\\
&-\int_{\mathbb R^2}\Delta_j\Delta_k^h\nabla\psi_x \Delta_j\Delta_k^h(\nabla^2 b \nabla\psi_y) \, dxdy
-\int_{\mathbb R^2}\Delta_j\Delta_k^h\nabla\psi_x \Delta_j\Delta_k^h(\nabla b \nabla^2\psi_y) \, dxdy\\
=&\int_{\mathbb R^2}\Delta_j\Delta_k^h\nabla\psi_y \Delta_j\Delta_k^h(\nabla^2 b \nabla\psi_x) \, dxdy
-\int_{\mathbb R^2}\Delta_j\Delta_k^h\nabla\psi_x \Delta_j\Delta_k^h(\nabla^2 b \nabla\psi_y) \, dxdy\\
&+\int_{\mathbb R^2}\Delta_j\Delta_k^h\nabla\psi_y \nabla b \Delta_j\Delta_k^h\nabla^2\psi_x \, dxdy
+\int_{\mathbb R^2}\Delta_j\Delta_k^h\nabla\psi_y[\Delta_j\Delta_k^h, \nabla b\nabla] \nabla\psi_x \, dxdy\\
&-\int_{\mathbb R^2}\Delta_j\Delta_k^h\nabla\psi_x \nabla b\Delta_j\Delta_k^h\nabla^2\psi_y \, dxdy
-\int_{\mathbb R^2}\Delta_j\Delta_k^h\nabla\psi_x [\Delta_j\Delta_k^h, \nabla b\nabla]\nabla\psi_y \, dxdy\\
=&\int_{\mathbb R^2}\Delta_j\Delta_k^h\nabla\psi_y \Delta_j\Delta_k^h(\nabla^2 b \nabla\psi_x) \, dxdy
-\int_{\mathbb R^2}\Delta_j\Delta_k^h\nabla\psi_x \Delta_j\Delta_k^h(\nabla^2 b \nabla\psi_y) \, dxdy\\
&+\int_{\mathbb R^2}\Delta_j\Delta_k^h\nabla\psi_y[\Delta_j\Delta_k^h, \nabla b\nabla] \nabla\psi_x \, dxdy
-\int_{\mathbb R^2}\Delta_j\Delta_k^h\nabla\psi_x [\Delta_j\Delta_k^h, \nabla b\nabla]\nabla\psi_y \, dxdy\\
&-\int_{\mathbb R^2}\Delta_j\Delta_k^h\nabla\psi_y \nabla b_x \Delta_j\Delta_k^h\nabla^2\psi \, dxdy
+\int_{\mathbb R^2}\Delta_j\Delta_k^h\nabla\psi_x \nabla b_y\Delta_j\Delta_k^h\nabla^2\psi \, dxdy
\end{split}
\end{equation}
in view of the cancellations 
\begin{equation}\notag
\begin{split}
&\int_{\mathbb R^2}\Delta_j\Delta_k^h\nabla\psi_y \nabla b \Delta_j\Delta_k^h\nabla^2\psi_x \, dxdy
-\int_{\mathbb R^2}\Delta_j\Delta_k^h\nabla\psi_x \nabla b\Delta_j\Delta_k^h\nabla^2\psi_y \, dxdy\\
=&-\int_{\mathbb R^2}\Delta_j\Delta_k^h\nabla\psi_{xy} \nabla b \Delta_j\Delta_k^h\nabla^2\psi \, dxdy
-\int_{\mathbb R^2}\Delta_j\Delta_k^h\nabla\psi_y \nabla b_x \Delta_j\Delta_k^h\nabla^2\psi \, dxdy\\
&+\int_{\mathbb R^2}\Delta_j\Delta_k^h\nabla\psi_{xy} \nabla b\Delta_j\Delta_k^h\nabla^2\psi \, dxdy
+\int_{\mathbb R^2}\Delta_j\Delta_k^h\nabla\psi_x \nabla b_y\Delta_j\Delta_k^h\nabla^2\psi \, dxdy\\
=&-\int_{\mathbb R^2}\Delta_j\Delta_k^h\nabla\psi_y \nabla b_x \Delta_j\Delta_k^h\nabla^2\psi \, dxdy
+\int_{\mathbb R^2}\Delta_j\Delta_k^h\nabla\psi_x \nabla b_y\Delta_j\Delta_k^h\nabla^2\psi \, dxdy.
\end{split}
\end{equation}

Observing the structures, the first two terms of (\ref{high-5678-last}) can be estimated similarly, the middle two terms can be estimated similarly, and so are the last two terms. Therefore, we only present the estimates of the first, third and fifth terms in (\ref{high-5678-last}).
The first integral in (\ref{high-5678-last}) can be further decomposed as
\begin{equation}\notag
\begin{split}
&\lambda_j^{2s_1}\lambda_k^{-2s}\int_{\mathbb R^2}\Delta_j\Delta_k^h\nabla\psi_y \Delta_j\Delta_k^h(\nabla^2 b \nabla\psi_x) \, dxdy\\
=&\sum_{m=1}^9 \lambda_j^{2s_1}\lambda_k^{-2s} \int_{\mathbb R^2}\Delta_j\Delta_k^h \nabla\psi_y \mathcal B_{j,k}^m(\nabla^2 b,\nabla\psi_x) \, dxdy.
\end{split}
\end{equation}
The nine terms are estimated as
\begin{equation}\notag
\begin{split}
&\lambda_j^{2s_1}\lambda_k^{-2s}\left| \int_{\mathbb R^2}\Delta_j\Delta_k^h \nabla\psi_y \mathcal B_{j,k}^1(\nabla^2 b, \nabla\psi_x) \, dxdy\right|\\
\leq&\ \lambda_j^{2s_1}\lambda_k^{-2s} \sum_{|j-j'|\leq 2, |k-k'|\leq 2}\int_{\mathbb R^2}\left|\Delta_j\Delta_k^h \nabla\psi_y \Delta_j\Delta_k^h\left(S_{j'-1}S_{k'-1}^h\nabla^2 b\Delta_{j'}\Delta_{k'}^h \nabla\psi_x \right) \right|\, dxdy\\
\leq &\ \lambda_j^{2s_1}\lambda_k^{-2s} \sum_{|j-j'|\leq 2, |k-k'|\leq 2} \|\Delta_j\Delta_k^h \nabla\psi_y\|_{L^2}\|S_{j'-1}S_{k'-1}^h\nabla^2 b\|_{L^\infty}\|\Delta_{j'}\Delta_{k'}^h \nabla\psi_x\|_{L^2}\\
\lesssim&\sum_{|j-j'|\leq 2, |k-k'|\leq 2} \|\nabla^2b\|_{B^{\frac12,\frac12}} \|\Delta_j\Delta_k^h|D_x|^{-s}D^{s_1}\nabla\psi_y\|_{L^2}\|\Delta_{j'}\Delta_{k'}^h|D_x|^{-s}D^{s_1}\nabla\psi_x\|_{L^2}\\
\lesssim &\  \|\nabla^2 b\|_{B^{\frac12,\frac12}}\|\nabla \psi\|_{H^{-s, s_1+1}}^2,
\end{split}
\end{equation}

\begin{equation}\notag
\begin{split}
&\lambda_j^{2s_1}\lambda_k^{-2s}\left| \int_{\mathbb R^2}\Delta_j\Delta_k^h\nabla\psi_y \mathcal B_{j,k}^2(\nabla^2b, \nabla\psi_x) \, dxdy\right|\\
\leq&\ \lambda_j^{2s_1}\lambda_k^{-2s}\sum_{|j-j'|\leq 2, |k-k'|\leq 2}\int_{\mathbb R^2}\left|\Delta_j\Delta_k^h\nabla\psi_y \Delta_j\Delta_k^h\left(\Delta_{j'}S_{k'-1}^h\nabla^2bS_{j'-1}\Delta_{k'}^h\nabla\psi_x \right) \right|\, dxdy\\
\leq &\ \lambda_j^{2s_1}\lambda_k^{-2s} \sum_{|j-j'|\leq 2, |k-k'|\leq 2} \|\Delta_j\Delta_k^h\nabla\psi_y\|_{L^2}\|\Delta_{j'}\nabla^2b\|_{L_x^\infty L_y^2}\|\Delta_{k'}^h\nabla\psi_x\|_{L_x^2 L_y^\infty}\\
\lesssim&\sum_{|j-j'|\leq 2, |k-k'|\leq 2} \left(\lambda_{j'}^{\frac12}\|\Delta_{j'}\nabla^2b\|_{L_x^\infty L_y^2}\right) \|\Delta_j\Delta_k^h|D_x|^{-s}D^{s_1-1}\nabla \psi_y\|_{L^2}\\
&\cdot\|\Delta_{k'}^h|D_x|^{-s}D^{s_1+\frac12}\psi_x\|_{L_x^2 L_y^\infty}\\
\lesssim&\ \|\nabla^2b\|_{B^{\frac12,\frac12}} \|\nabla \psi\|_{H^{-s, s_1+1}}^2,
\end{split}
\end{equation}

\begin{equation}\notag
\begin{split}
&\lambda_j^{2s_1}\lambda_k^{-2s}\left| \int_{\mathbb R^2}\Delta_j\Delta_k^h\nabla\psi_y \mathcal B_{j,k}^3(\nabla^2b, \nabla\psi_x) \, dxdy\right|\\
\leq&\ \lambda_j^{2s_1}\lambda_k^{-2s} \sum_{j'\geq j-2, |k-k'|\leq 2}\int_{\mathbb R^2}\left|\Delta_j\Delta_k^h\nabla\psi_y \Delta_j\Delta_k^h\left(\widetilde{\Delta}_{j'}S_{k'-1}^h\nabla^2b\Delta_{j'}\Delta_{k'}^h\nabla\psi_x \right) \right|\, dxdy\\
\leq &\ \lambda_j^{2s_1}\lambda_k^{-2s} \sum_{j'\geq j- 2, |k-k'|\leq 2} \|\Delta_j\Delta_k^h\nabla\psi_y\|_{L^\infty}\|\widetilde{\Delta}_{j'}S_{k'-1}^h\nabla^2b\|_{L^2}\|\Delta_{j'}\Delta_{k'}^h\nabla\psi_x\|_{L^2}\\
\lesssim&\sum_{j'\geq j- 2, |k-k'|\leq 2} \|\nabla\psi_y\|_{B^{\frac12,\frac12}} \|\widetilde{\Delta}_{j'}S_{k'-1}^h|D_x|^{-s}D^{s_1}\nabla^2b\|_{L^2}\|\Delta_{j'}\Delta_{k'}^h|D_x|^{-s}D^{s_1}\nabla\psi_x\|_{L^2}\\
\lesssim&\ \|\nabla\psi_y\|_{B^{\frac12,\frac12}}\|\nabla b\|_{H^{-s, s_1+1}}\|\nabla \psi_x\|_{H^{-s, s_1}}\\
\leq&\ \frac{1}{256} c\left(\|\nabla b\|_{H^{-s,s_1+1}}^2+\|\nabla \psi_x\|_{H^{-s, s_1}}^2\right)+C\|\nabla\psi_y\|_{B^{\frac12,\frac12}}^4\|\nabla \psi_x\|_{H^{-s, s_1}}^2,
\end{split}
\end{equation}

\begin{equation}\notag
\begin{split}
&\lambda_j^{2s_1}\lambda_k^{-2s}\left| \int_{\mathbb R^2}\Delta_j\Delta_k^h\nabla\psi_y \mathcal B_{j,k}^4(\nabla^2b, \nabla\psi_x) \, dxdy\right|\\
\leq&\ \lambda_j^{2s_1}\lambda_k^{-2s} \sum_{|j-j'|\leq 2, |k-k'|\leq 2}\int_{\mathbb R^2}\left|\Delta_j\Delta_k^h\nabla\psi_y \Delta_j\Delta_k^h\left(S_{j'-1}\Delta_{k'}^h\nabla^2b\Delta_{j'}S_{k'-1}^h\nabla\psi_x \right) \right|\, dxdy\\
\leq &\ \lambda_j^{2s_1}\lambda_k^{-2s} \sum_{|j-j'|\leq 2, |k-k'|\leq 2} \|\Delta_j\Delta_k^h\nabla\psi_y\|_{L^2}\|\Delta_{k'}^h\nabla^2b\|_{L^\infty}\|\Delta_{j'}S_{k'-1}^h\nabla\psi_x\|_{L^2}\\
\lesssim&\sum_{|j-j'|\leq 2, |k-k'|\leq 2} \|\nabla^2b\|_{B^{\frac12,\frac12}} \|\Delta_j\Delta_k^h|D_x|^{-s} D^{s_1}\nabla\psi_y\|_{L^2}\|\Delta_{j'}S_{k'-1}^h|D_x|^{-s}D^{s_1}\nabla\psi_x\|_{L^2}\\
\lesssim&\ \|\nabla^2b\|_{B^{\frac12,\frac12}} \|\nabla \psi\|_{H^{-s, s_1+1}}^2,
\end{split}
\end{equation}

\begin{equation}\notag
\begin{split}
&\lambda_j^{2s_1}\lambda_k^{-2s}\left| \int_{\mathbb R^2}\Delta_j\Delta_k^h\nabla\psi_y \mathcal B_{j,k}^5(\nabla^2b, \nabla\psi_x) \, dxdy\right|\\
\leq&\ \lambda_j^{2s_1}\lambda_k^{-2s} \sum_{|j-j'|\leq 2, |k-k'|\leq 2}\int_{\mathbb R^2}\left|\Delta_j\Delta_k^h\nabla\psi_y \Delta_j\Delta_k^h\left(\Delta_{j'}\Delta_{k'}^h\nabla^2b S_{j'-1}S_{k'-1}^h\nabla\psi_x \right) \right|\, dxdy\\
\leq &\ \lambda_j^{2s_1}\lambda_k^{-2s} \sum_{|j-j'|\leq 2, |k-k'|\leq 2} \|\Delta_j\Delta_k^h \nabla\psi_y\|_{L^2}\|\Delta_{j'}\Delta_{k'}^h\nabla^2b \|_{L^2}\|S_{j'-1}S_{k'-1}^h\nabla\psi_x\|_{L^\infty}\\
\lesssim&\sum_{|j-j'|\leq 2, |k-k'|\leq 2} \|\nabla\psi_x\|_{L^{\infty}} \|\Delta_j\Delta_k^h|D_x|^{-s}D^{s_1}\nabla\psi_y\|_{L^2}\|\Delta_{j'}\Delta_{k'}^h|D_x|^{-s}D^{s_1}\nabla^2b \|_{L^2}\\
\lesssim&\ \|\nabla\psi_x\|_{B^{\frac12,\frac12}}\|\nabla b\|_{H^{-s, s_1+1}}\|\nabla \psi\|_{H^{-s, s_1+1}}\\
\leq&\ \frac{1}{256}  c \|\nabla b\|_{H^{-s, s_1+1}}^2+\|\nabla\psi_x\|_{B^{\frac12,\frac12}}^2\|\nabla \psi\|_{H^{-s, s_1+1}}^2,
\end{split}
\end{equation}

\begin{equation}\notag
\begin{split}
&\lambda_j^{2s_1}\lambda_k^{-2s}\left| \int_{\mathbb R^2}\Delta_j\Delta_k^h\nabla\psi_y \mathcal B_{j,k}^6(\nabla^2b, \nabla\psi_x) \, dxdy\right|\\
\leq&\ \lambda_j^{2s_1}\lambda_k^{-2s} \sum_{j'\geq j-2, |k-k'|\leq 2}\int_{\mathbb R^2}\left|\Delta_j\Delta_k^h\nabla\psi_y \Delta_j\Delta_k^h\left(\widetilde{\Delta}_{j'}\Delta_{k'}^h\nabla^2b\Delta_{j'}S_{k'-1}^h\nabla\psi_x \right) \right|\, dxdy\\
\leq &\ \lambda_j^{2s_1}\lambda_k^{-2s} \sum_{j'\geq j- 2, |k-k'|\leq 2} \|\Delta_j\Delta_k^h\nabla\psi_y\|_{L_x^2L_y^\infty}\|\widetilde{\Delta}_{j'}\Delta_{k'}^h\nabla^2b\|_{L^2}\|\Delta_{j'}S_{k'-1}^h\nabla\psi_x\|_{L_x^\infty L_y^2}\\
\lesssim&\sum_{j'\geq j- 2, |k-k'|\leq 2}\lambda_{j'-j}^{-s_1-\frac12} \left(\lambda_{j'}^{\frac12}\|\Delta_{j'}S_{k'-1}^h\nabla\psi_x\|_{L_x^\infty L_y^2}\right)\\
&\cdot \|\Delta_j\Delta_k^h|D_x|^{-s}D^{s_1-\frac12}\nabla\psi_y\|_{L_x^2L_y^\infty}\|\widetilde{\Delta}_{j'}\Delta_{k'}^h|D_x|^{-s}D^{s_1}\nabla^2b\|_{L^2}\\
\lesssim&\ \|\nabla\psi_x\|_{B^{\frac12,\frac12}}\|\nabla b\|_{H^{-s, s_1+1}}\|\nabla \psi\|_{H^{-s, s_1+1}}\\
\leq&\ \frac{1}{256}  c \|\nabla b\|_{H^{-s, s_1+1}}^2+\|\nabla\psi_x\|_{B^{\frac12,\frac12}}^2\|\nabla \psi\|_{H^{-s, s_1+1}}^2
\end{split}
\end{equation}
for $s_1>-\frac12$,

\begin{equation}\notag
\begin{split}
&\lambda_j^{2s_1}\lambda_k^{-2s}\left| \int_{\mathbb R^2}\Delta_j\Delta_k^h\nabla\psi_y \mathcal B_{j,k}^7(\nabla^2b, \nabla\psi_x) \, dxdy\right|\\
\leq&\ \lambda_j^{2s_1}\lambda_k^{-2s} \sum_{|j-j'|\leq 2, k'\geq k- 2}\int_{\mathbb R^2}\left|\Delta_j\Delta_k^h\nabla\psi_y \Delta_j\Delta_k^h\left(S_{j'-1}\widetilde{\Delta}_{k'}^h\nabla^2b\Delta_{j'}\Delta_{k'}^h\nabla\psi_x \right) \right|\, dxdy\\
\leq &\ \lambda_j^{2s_1}\lambda_k^{-2s} \sum_{|j-j'|\leq 2, k'\geq k- 2} \|\Delta_j\Delta_k^h\nabla\psi_y\|_{L_x^\infty L_y^2}\|\widetilde{\Delta}_{k'}^h\nabla^2b\|_{L_x^2L_y^\infty}\|\Delta_{j'}\Delta_{k'}^h\nabla\psi_x\|_{L^2}\\
\lesssim&\sum_{|j-j'|\leq 2, k'\geq k- 2} \lambda_{k'-k}^{s-\frac12}\left(\lambda_{k'}^{\frac12} \|\widetilde{\Delta}_{k'}^h\nabla^2b\|_{L_x^2L_y^\infty}\right)\|\Delta_j\Delta_k^h|D_x|^{-s-\frac12}D^{s_1}\nabla\psi_y\|_{L_x^\infty L_y^2}\\
&\cdot\|\Delta_{j'}\Delta_{k'}^h|D_x|^{-s}D^{s_1}\nabla\psi_x\|_{L^2}\\
\lesssim&\ \|\nabla^2b\|_{B^{\frac12,\frac12}} \|\nabla \psi\|_{H^{-s, s_1+1}}^2
\end{split}
\end{equation}
for $s<\frac12$,

\begin{equation}\notag
\begin{split}
&\lambda_j^{2s_1}\lambda_k^{-2s}\left| \int_{\mathbb R^2}\Delta_j\Delta_k^h\nabla\psi_y \mathcal B_{j,k}^8(\nabla^2b, \nabla\psi_x) \, dxdy\right|\\
\leq&\ \lambda_j^{2s_1}\lambda_k^{-2s} \sum_{|j-j'|\leq 2, k'\geq k- 2}\int_{\mathbb R^2}\left|\Delta_j\Delta_k^h\nabla\psi_y \Delta_j\Delta_k^h\left(\Delta_{j'}\widetilde{\Delta}_{k'}^h\nabla^2b S_{j'-1}\Delta_{k'}^h\nabla\psi_x \right) \right|\, dxdy\\
\leq &\ \lambda_j^{2s_1}\lambda_k^{-2s} \sum_{|j-j'|\leq 2, k'\geq k- 2} \|\Delta_j\Delta_k^h\nabla\psi_y\|_{L_x^\infty L_y^2}\|\Delta_{j'}\widetilde{\Delta}_{k'}^h\nabla^2b\|_{L^2}\|\Delta_{k'}^h\nabla\psi_x\|_{L_x^2L_y^\infty}\\
\lesssim&\sum_{|j-j'|\leq 2, k'\geq k- 2} \lambda_{k'-k}^{s-\frac12}\left(\lambda_{k'}^{\frac12}\|\Delta_{k'}^h\nabla\psi_x\|_{L_x^2L_y^\infty}\right) \|\Delta_j\Delta_k^h|D_x|^{-s-\frac12}D^{s_1}\nabla\psi_y\|_{L_x^\infty L_y^2}\\
&\cdot\|\Delta_{j'}\widetilde{\Delta}_{k'}^h|D_x|^{-s}D^{s_1}\nabla^2b\|_{L^2}\\
\lesssim&\ \|\nabla\psi_x\|_{B^{\frac12,\frac12}}\|\nabla b\|_{H^{-s, s_1+1}}\|\nabla \psi\|_{H^{-s, s_1+1}}\\
\leq&\ \frac{1}{256} c \|\nabla b\|_{H^{-s, s_1+1}}^2+\|\nabla\psi_x\|_{B^{\frac12,\frac12}}^2\|\nabla \psi\|_{H^{-s, s_1+1}}^2
\end{split}
\end{equation}
for $s<\frac12$,

\begin{equation}\notag
\begin{split}
&\lambda_j^{2s_1}\lambda_k^{-2s}\left| \int_{\mathbb R^2}\Delta_j\Delta_k^h\nabla\psi_y \mathcal B_{j,k}^9(\nabla^2b, \nabla\psi_x) \, dxdy\right|\\
\leq&\ \lambda_j^{2s_1}\lambda_k^{-2s} \sum_{j'\geq j- 2, k'\geq k- 2}\int_{\mathbb R^2}\left|\Delta_j\Delta_k^h\nabla\psi_y \Delta_j\Delta_k^h\left(\widetilde{\Delta}_{j'}\widetilde{\Delta}_{k'}^h\nabla^2b \Delta_{j'}\Delta_{k'}^h\nabla\psi_x \right) \right|\, dxdy\\
\leq &\ \lambda_j^{2s_1}\lambda_k^{-2s} \sum_{j'\geq j- 2, k'\geq k- 2} \|\Delta_j\Delta_k^h\nabla\psi_y\|_{L^\infty}\|\widetilde{\Delta}_{j'}\widetilde{\Delta}_{k'}^h\nabla^2b\|_{L^2}\|\Delta_{j'}\Delta_{k'}^h\nabla\psi_x\|_{L^2}\\
\lesssim&\sum_{j'\geq j- 2, k'\geq k- 2}\lambda_{j'-j}^{-s_1-\frac12} \lambda_{k'-k}^{s-\frac12}\left(\lambda_{j'}^{\frac12}\lambda_{k'}^{\frac12}\|\Delta_{j'}\Delta_{k'}^h\nabla\psi_x\|_{L^2} \right)\\
&\cdot \|\Delta_j\Delta_k^h|D_x|^{-s-\frac12}D^{s_1-\frac12}\nabla\psi_y\|_{L^\infty}\|\widetilde{\Delta}_{j'}\widetilde{\Delta}_{k'}^h|D_x|^{-s}D^{s_1}\nabla^2b\|_{L^2}\\
\lesssim&\ \|\nabla\psi_x\|_{B^{\frac12,\frac12}}\|\nabla b\|_{H^{-s, s_1+1}}\|\nabla \psi\|_{H^{-s, s_1+1}}\\
\leq&\ \frac{1}{256}  c \|\nabla b\|_{H^{-s, s_1+1}}^2+\|\nabla\psi_x\|_{B^{\frac12,\frac12}}^2\|\nabla \psi\|_{H^{-s, s_1+1}}^2
\end{split}
\end{equation}
for $s<\frac12$ and $s_1>-\frac12$.

Applying commutator estimate to the third and fifth integrals in (\ref{high-5678-last}) gives
\begin{equation}\notag
\begin{split}
&\lambda_j^{2s_1}\lambda_k^{-2s}\left|\int_{\mathbb R^2}\Delta_j\Delta_k^h\nabla\psi_y[\Delta_j\Delta_k^h, \nabla b\nabla] \nabla\psi_x \, dxdy\right|\\
\leq&\ \lambda_j^{2s_1}\lambda_k^{-2s}\|\Delta_j\Delta_k^h\nabla\psi_y\|_{L^2}\|\Delta_j\Delta_k^h\nabla^2b\|_{L^2} \|\nabla\psi_x\|_{L^\infty}\\
\leq&\ \|\nabla\psi_x\|_{B^{\frac12,\frac12}}\|\Delta_j\Delta_k^h|D_x|^{-s}D^{s_1}\nabla\psi_y\|_{L^2}\|\Delta_j\Delta_k^h|D_x|^{-s}D^{s_1}\nabla^2b\|_{L^2}\\
\lesssim&\ \|\nabla\psi_x\|_{B^{\frac12,\frac12}}\|\nabla b\|_{H^{-s, s_1+1}}\|\nabla \psi\|_{H^{-s, s_1+1}}\\
\leq&\ \frac{1}{256} c \|\nabla b\|_{H^{-s, s_1+1}}^2+\|\nabla\psi_x\|_{B^{\frac12,\frac12}}^2\|\nabla \psi\|_{H^{-s, s_1+1}}^2,
\end{split}
\end{equation}
\begin{equation}\notag
\begin{split}
&\lambda_j^{2s_1}\lambda_k^{-2s}\left|\int_{\mathbb R^2}\Delta_j\Delta_k^h\nabla\psi_y\nabla b_x\Delta_j\Delta_k^h\nabla^2\psi \, dxdy\right|\\
\leq&\ \lambda_j^{2s_1}\lambda_k^{-2s}\|\Delta_j\Delta_k^h\nabla\psi_y\|_{L^2}\|\Delta_j\Delta_k^h\nabla^2\psi\|_{L^2} \|\nabla b_x\|_{L^\infty}\\
\leq&\ \|\nabla b_x\|_{B^{\frac12,\frac12}}\|\Delta_j\Delta_k^h|D_x|^{-s}D^{s_1}\nabla\psi_y\|_{L^2}\|\Delta_j\Delta_k^h|D_x|^{-s}D^{s_1}\nabla^2\psi\|_{L^2}\\
\lesssim&\ \|\nabla b_x\|_{B^{\frac12,\frac12}}\|\nabla \psi\|_{H^{-s, s_1+1}}^2.
\end{split}
\end{equation}

The proof of the lemma is complete.
\cbdu

\medskip

\begin{Lemma}\label{le-high-aux910}
Let $0<s<\frac12$ and $s_1>-\frac12$. We have 
\begin{equation}\notag
\begin{split}
&\sum_{j,k}\lambda_j^{2s_1}\lambda_k^{-2s}|\tilde I_9+\tilde I_{10}|\\
\leq &\ C \left(\|\nabla\psi\|_{B^{\frac32,\frac12}}+\|\nabla\psi\|_{B^{\frac12,\frac32}} \right)\left( \|\nabla\psi\|_{H^{-s,s_1}}^2+ \|\nabla\psi\|_{H^{-s,s_1+1}}^2\right)\\
&+C \left(\|\psi\|_{B^{\frac32,\frac12}}^2 +\|\psi\|_{B^{\frac12,\frac32}}^4\right) \|\nabla\psi\|_{H^{-s,s_1+1}}^2+\frac{1}{16} c\|\nabla\psi_x\|_{H^{-s,s_1}}^2
\end{split}
\end{equation}
for a constant $C>0$.
\end{Lemma}
\pf
First, applying integration by parts and observing cancellations gives 
\begin{equation}\notag
\begin{split}
&\tilde I_9+\tilde I_{10}\\
=&\ \varepsilon_1 \int_{\mathbb R^2}\Delta_j\Delta_k^h\psi_x \Delta_j\Delta_k^h(\psi_y \Delta\psi_x) \, dxdy
-\varepsilon_1 \int_{\mathbb R^2}\Delta_j\Delta_k^h\psi_x \Delta_j\Delta_k^h(\psi_x \Delta\psi_y) \, dxdy\\
=&-\varepsilon_1 \int_{\mathbb R^2}\Delta_j\Delta_k^h\nabla\psi_x \Delta_j\Delta_k^h(\psi_y \nabla\psi_x) \, dxdy
-\varepsilon_1 \int_{\mathbb R^2}\Delta_j\Delta_k^h\psi_x \Delta_j\Delta_k^h(\nabla\psi_y \nabla\psi_x) \, dxdy\\
&+\varepsilon_1 \int_{\mathbb R^2}\Delta_j\Delta_k^h\nabla\psi_x \Delta_j\Delta_k^h(\psi_x \nabla\psi_y) \, dxdy
+\varepsilon_1 \int_{\mathbb R^2}\Delta_j\Delta_k^h\psi_x \Delta_j\Delta_k^h(\nabla\psi_x \nabla\psi_y) \, dxdy\\
=&-\varepsilon_1 \int_{\mathbb R^2}\Delta_j\Delta_k^h\nabla\psi_x \Delta_j\Delta_k^h(\psi_y \nabla\psi_x) \, dxdy
+\varepsilon_1 \int_{\mathbb R^2}\Delta_j\Delta_k^h\nabla\psi_x \Delta_j\Delta_k^h (\psi_x\nabla\psi_y) \, dxdy\\
\end{split}
\end{equation}
We only need to show the estimates for one integral and the other one can be analyzed in an analogous way.
The first integral in $\tilde I_9+\tilde I_{10}$ is decomposed as
\begin{equation}\notag
\begin{split}
&\lambda_j^{2s_1} \lambda_k^{-2s}\int_{\mathbb R^2}\Delta_j\Delta_k^h\nabla\psi_x \Delta_j\Delta_k^h(\nabla \psi_x \psi_y) \, dxdy\\
=&\sum_{m=1}^9  \lambda_j^{2s_1} \lambda_k^{-2s}\int_{\mathbb R^2}\Delta_j\Delta_k^h \nabla\psi_x \mathcal B_{j,k}^m(\nabla \psi_x,  \psi_y) \, dxdy.
\end{split}
\end{equation}
For $0<s<\frac12$ and $s_1>-\frac12$, the estimates of the nine terms are given as
\begin{equation}\notag
\begin{split}
&\lambda_j^{2s_1} \lambda_k^{-2s}\left| \int_{\mathbb R^2}\Delta_j\Delta_k^h \nabla\psi_x \mathcal B_{j,k}^1(\nabla \psi_x, \psi_y) \, dxdy\right|\\
\leq&\ \lambda_j^{2s_1} \lambda_k^{-2s} \sum_{|j-j'|\leq 2, |k-k'|\leq 2}\int_{\mathbb R^2}\left|\Delta_j\Delta_k^h \nabla\psi_x \Delta_j\Delta_k^h\left(S_{j'-1}S_{k'-1}^h\nabla \psi_x\Delta_{j'}\Delta_{k'}^h\psi_y \right) \right|\, dxdy\\
\leq &\ \lambda_j^{2s_1} \lambda_k^{-2s} \sum_{|j-j'|\leq 2, |k-k'|\leq 2} \|\Delta_j\Delta_k^h \nabla\psi_x\|_{L^2}\|S_{j'-1}S_{k'-1}^h\nabla \psi_x\|_{L^\infty}\|\Delta_{j'}\Delta_{k'}^h\psi_y\|_{L^2}\\
\lesssim&\sum_{|j-j'|\leq 2, |k-k'|\leq 2} \|\nabla\psi_x\|_{B^{\frac12,\frac12}} \|\Delta_j\Delta_k^h |D_x|^{-s}D^{s_1}\nabla\psi_x\|_{L^2}\|\Delta_{j'}\Delta_{k'}^h|D_x|^{-s}D^{s_1}\psi_y\|_{L^2}\\
\lesssim&\ \|\nabla \psi_x\|_{B^{\frac12,\frac12}} \|\nabla\psi\|_{H^{-s,s_1}} \|\nabla\psi\|_{H^{-s,s_1+1}},
\end{split}
\end{equation}
\begin{equation}\notag
\begin{split}
&\lambda_j^{2s_1} \lambda_k^{-2s}\left| \int_{\mathbb R^2}\Delta_j\Delta_k^h\nabla\psi_x \mathcal B_{j,k}^2(\nabla \psi_x, \psi_y) \, dxdy\right|\\
\leq&\ \lambda_j^{2s_1} \lambda_k^{-2s}\sum_{|j-j'|\leq 2, |k-k'|\leq 2}\int_{\mathbb R^2}\left|\Delta_j\Delta_k^h\nabla\psi_x \Delta_j\Delta_k^h\left(\Delta_{j'}S_{k'-1}^h\nabla \psi_xS_{j'-1}\Delta_{k'}^h\psi_y \right) \right|\, dxdy\\
\leq &\ \lambda_j^{2s_1} \lambda_k^{-2s} \sum_{|j-j'|\leq 2, |k-k'|\leq 2} \|\Delta_j\Delta_k^h\nabla\psi_x\|_{L^2}\|\Delta_{j'}\nabla \psi_x\|_{L_x^\infty L_y^2}\|\Delta_{k'}^h\psi_y\|_{L_x^2 L_y^\infty}\\
\lesssim&\sum_{|j-j'|\leq 2, |k-k'|\leq 2} \|\nabla\psi_x\|_{B^{\frac12,\frac12}}\|\Delta_j\Delta_k^h |D_x|^{-s}D^{s_1}\nabla\psi_x\|_{L^2}\|\Delta_{k'}^h |D_x|^{-s}D^{s_1-\frac12}\psi_y\|_{L_x^2 L_y^\infty}\\
\lesssim&\ \|\nabla \psi_x\|_{B^{\frac12,\frac12}} \|\nabla\psi\|_{H^{-s,s_1}} \|\nabla\psi\|_{H^{-s,s_1+1}},
\end{split}
\end{equation}
\begin{equation}\notag
\begin{split}
&\lambda_j^{2s_1} \lambda_k^{-2s}\left| \int_{\mathbb R^2}\Delta_j\Delta_k^h \nabla\psi_x \mathcal B_{j,k}^3(\nabla \psi_x, \psi_y) \, dxdy\right|\\
\leq&\ \lambda_j^{2s_1} \lambda_k^{-2s} \sum_{j'\geq j-2, |k-k'|\leq 2}\int_{\mathbb R^2}\left|\Delta_j\Delta_k^h \nabla\psi_x \Delta_j\Delta_k^h\left(\widetilde{\Delta}_{j'}S_{k'-1}^h\nabla \psi_x\Delta_{j'}\Delta_{k'}^h\psi_y \right) \right|\, dxdy\\
\leq &\ \lambda_j^{2s_1} \lambda_k^{-2s} \sum_{j'\geq j- 2, |k-k'|\leq 2} \|\Delta_j\Delta_k^h\nabla\psi_x\|_{L_x^2L_y^\infty}\|\widetilde{\Delta}_{j'}S_{k'-1}^h\nabla \psi_x\|_{L_x^\infty L_y^2}\|\Delta_{j'}\Delta_{k'}^h\psi_y\|_{L^2}\\
\lesssim&\sum_{j'\geq j- 2, |k-k'|\leq 2} \lambda_{j'-j}^{-s_1-\frac12}\left(\lambda_{j'}^{\frac12}\|\widetilde{\Delta}_{j'}S_{k'-1}^h\nabla \psi_x\|_{L_x^\infty L_y^2}\right)\\
 &\cdot\|\Delta_j\Delta_k^h |D_x|^{-s}D^{s_1-\frac12}\nabla\psi_x\|_{L_x^2L_y^\infty} \|\Delta_{j'}\Delta_{k'}^h |D_x|^{-s}D^{s_1}\psi_y\|_{L^2}\\
\lesssim&\ \|\nabla \psi_x\|_{B^{\frac12,\frac12}} \|\nabla\psi\|_{H^{-s,s_1}} \|\nabla\psi\|_{H^{-s,s_1+1}},
\end{split}
\end{equation}

\begin{equation}\notag
\begin{split}
&\lambda_j^{2s_1} \lambda_k^{-2s}\left| \int_{\mathbb R^2}\Delta_j\Delta_k^h \nabla\psi_x \mathcal B_{j,k}^4(\nabla \psi_x,  \psi_y) \, dxdy\right|\\
\leq&\ \lambda_j^{2s_1} \lambda_k^{-2s} \sum_{|j-j'|\leq 2, |k-k'|\leq 2}\int_{\mathbb R^2}\left|\Delta_j\Delta_k^h \nabla\psi_x \Delta_j\Delta_k^h\left(S_{j'-1}\Delta_{k'}^h\nabla \psi_x\Delta_{j'}S_{k'-1}^h\psi_y \right) \right|\, dxdy\\
\leq &\ \lambda_j^{2s_1} \lambda_k^{-2s} \sum_{|j-j'|\leq 2, |k-k'|\leq 2} \|\Delta_j\Delta_k^h \nabla\psi_x\|_{L^2}\|\Delta_{k'}^h\nabla \psi_x\|_{L^\infty}\|\Delta_{j'}S_{k'-1}^h\psi_y\|_{L^2}\\
\lesssim &\ \sum_{|j-j'|\leq 2, |k-k'|\leq 2} \|\nabla \psi_x\|_{B^{\frac12,\frac12}}\|\Delta_j\Delta_k^h|D_x|^{-s}D^{s_1}\nabla\psi_x\|_{L^2}\|\Delta_{j'}S_{k'-1}^h|D_x|^{-s}D^{s_1}\psi_y\|_{L^2}\\
\lesssim&\ \|\nabla\psi_x\|_{B^{\frac12,\frac12}} \|\nabla\psi\|_{H^{-s,s_1}} \|\nabla\psi\|_{H^{-s,s_1+1}},
\end{split}
\end{equation}

\begin{equation}\notag
\begin{split}
&\lambda_j^{2s_1} \lambda_k^{-2s}\left| \int_{\mathbb R^2}\Delta_j\Delta_k^h\nabla\psi_x \mathcal B_{j,k}^5(\nabla \psi_x, \psi_y) \, dxdy\right|\\
\leq&\ \lambda_j^{2s_1} \lambda_k^{-2s} \sum_{|j-j'|\leq 2, |k-k'|\leq 2}\int_{\mathbb R^2}\left|\Delta_j\Delta_k^h\nabla\psi_x \Delta_j\Delta_k^h\left(\Delta_{j'}\Delta_{k'}^h\nabla \psi_x S_{j'-1}S_{k'-1}^h\psi_y \right) \right|\, dxdy\\
\leq &\ \lambda_j^{2s_1} \lambda_k^{-2s} \sum_{|j-j'|\leq 2, |k-k'|\leq 2} \|\Delta_j\Delta_k^h \nabla\psi_x\|_{L^2}\|\Delta_{j'}\Delta_{k'}^h\nabla \psi_x \|_{L^2}\|S_{j'-1}S_{k'-1}^h\psi_y\|_{L^\infty}\\
\lesssim&\sum_{|j-j'|\leq 2, |k-k'|\leq 2} \|\psi_y\|_{L^\infty} \|\Delta_j\Delta_k^h|D_x|^{-s}D^{s_1}\nabla\psi_x\|_{L^2}\|\Delta_{j'}\Delta_{k'}^h|D_x|^{-s}D^{s_1}\nabla \psi_x \|_{L^2}\\
\lesssim&\ \|\psi_y\|_{B^{\frac12,\frac12}} \|\nabla\psi_x\|_{H^{-s,s_1}}^2\\
\leq&\ \frac{1}{64}  c\|\nabla\psi_x\|_{H^{-s,s_1}}^2+C\|\psi_y\|_{B^{\frac12,\frac12}}^4 \|\nabla\psi_x\|_{H^{-s,s_1}}^2,
\end{split}
\end{equation}

\begin{equation}\notag
\begin{split}
&\lambda_j^{2s_1} \lambda_k^{-2s}\left| \int_{\mathbb R^2}\Delta_j\Delta_k^h\nabla\psi_x \mathcal B_{j,k}^6(\nabla \psi_x, \psi_y) \, dxdy\right|\\
\leq&\ \lambda_j^{2s_1} \lambda_k^{-2s} \sum_{j'\geq j-2, |k-k'|\leq 2}\int_{\mathbb R^2}\left|\Delta_j\Delta_k^h\nabla\psi_x \Delta_j\Delta_k^h\left(\widetilde{\Delta}_{j'}\Delta_{k'}^h\nabla \psi_x\Delta_{j'}S_{k'-1}^h\psi_y \right) \right|\, dxdy\\
\leq &\ \lambda_j^{2s_1} \lambda_k^{-2s} \sum_{j'\geq j- 2, |k-k'|\leq 2} \|\Delta_j\Delta_k^h\nabla\psi_x\|_{L_x^2L_y^\infty}\|\widetilde{\Delta}_{j'}\Delta_{k'}^h\nabla \psi_x\|_{L^2}\|\Delta_{j'}S_{k'-1}^h\psi_y\|_{L_x^\infty L_y^2}\\
\lesssim&\sum_{j'\geq j- 2, |k-k'|\leq 2}\lambda_{j'-j}^{-s_1-\frac12}\left(\lambda_{j'}^{\frac12}\|\Delta_{j'}S_{k'-1}^h\psi_y\|_{L_x^\infty L_y^2}\right)\\
&\cdot \|\Delta_j\Delta_k^h|D_x|^{-s}D^{s_1-\frac12}\nabla\psi_x\|_{L_x^2L_y^\infty}\|\widetilde{\Delta}_{j'}\Delta_{k'}^h|D_x|^{-s}D^{s_1}\nabla \psi_x\|_{L^2}\\
\lesssim&\ \|\psi_y\|_{B^{\frac12,\frac12}} \|\nabla\psi_x\|_{H^{-s,s_1}}^2\\
\leq&\ \frac{1}{64} c\|\nabla\psi_x\|_{H^{-s,s_1}}^2+C\|\psi_y\|_{B^{\frac12,\frac12}}^4 \|\nabla\psi_x\|_{H^{-s,s_1}}^2,
\end{split}
\end{equation}

\begin{equation}\notag
\begin{split}
&\lambda_j^{2s_1} \lambda_k^{-2s}\left| \int_{\mathbb R^2}\Delta_j\Delta_k^h\nabla\psi_x \mathcal B_{j,k}^7(\nabla \psi_x, \psi_y) \, dxdy\right|\\
\leq&\ \lambda_j^{2s_1} \lambda_k^{-2s} \sum_{|j-j'|\leq 2, k'\geq k- 2}\int_{\mathbb R^2}\left|\Delta_j\Delta_k^h\nabla\psi_x \Delta_j\Delta_k^h\left(S_{j'-1}\widetilde{\Delta}_{k'}^h\nabla \psi_x\Delta_{j'}\Delta_{k'}^h\psi_y \right) \right|\, dxdy\\
\leq &\ \lambda_j^{2s_1} \lambda_k^{-2s} \sum_{|j-j'|\leq 2, k'\geq k- 2} \|\Delta_j\Delta_k^h\nabla\psi_x\|_{L_x^\infty L_y^2}\|S_{j'-1}\widetilde{\Delta}_{k'}^h\nabla \psi_x\|_{L_x^2L_y^\infty}\|\Delta_{j'}\Delta_{k'}^h\psi_y\|_{L^2}\\
\lesssim&\sum_{|j-j'|\leq 2, k'\geq k- 2} \lambda_{k'-k}^{s-\frac12}\left(\lambda_{k'}^{\frac12}\|S_{j'-1}\widetilde{\Delta}_{k'}^h\nabla \psi_x\|_{L_x^2L_y^\infty} \right)\\
&\cdot \|\Delta_j\Delta_k^h|D_x|^{-s-\frac12}D^{s_1-1}\nabla\psi_x\|_{L_x^\infty L_y^2}\|\Delta_{j'}\Delta_{k'}^h|D_x|^{-s}D^{s_1+1}\psi_y\|_{L^2}\\
\lesssim&\ \|\nabla\psi_x\|_{B^{\frac12,\frac12}} \|\nabla\psi\|_{H^{-s,s_1}} \|\nabla\psi\|_{H^{-s,s_1+1}},
\end{split}
\end{equation}

\begin{equation}\notag
\begin{split}
&\lambda_j^{2s_1} \lambda_k^{-2s}\left| \int_{\mathbb R^2}\Delta_j\Delta_k^h\nabla\psi_x \mathcal B_{j,k}^8(\nabla \psi_x, \psi_y) \, dxdy\right|\\
\leq&\ \lambda_j^{2s_1} \lambda_k^{-2s} \sum_{|j-j'|\leq 2, k'\geq k- 2}\int_{\mathbb R^2}\left|\Delta_j\Delta_k^h\nabla\psi_x \Delta_j\Delta_k^h\left(\Delta_{j'}\widetilde{\Delta}_{k'}^h\nabla \psi_x S_{j'-1}\Delta_{k'}^h\psi_y \right) \right|\, dxdy\\
\leq &\ \lambda_j^{2s_1} \lambda_k^{-2s} \sum_{|j-j'|\leq 2, k'\geq k- 2} \|\Delta_j\Delta_k^h\nabla\psi_x\|_{L_x^\infty L_y^2}\|\Delta_{j'}\widetilde{\Delta}_{k'}^h\nabla \psi_x\|_{L^2}\|S_{j'-1}\Delta_{k'}^h\psi_y\|_{L_x^2L_y^\infty}\\
\lesssim&\sum_{|j-j'|\leq 2, k'\geq k- 2} \lambda_{k'-k}^{s-\frac12}\left(\lambda_{k'}^{\frac12} \|S_{j'-1}\Delta_{k'}^h\psi_y\|_{L_x^2L_y^\infty}\right)\\
&\cdot\|\Delta_j\Delta_k^h|D_x|^{-s-\frac12}D^{s_1}\nabla\psi_x\|_{L_x^\infty L_y^2}\|\Delta_{j'}\widetilde{\Delta}_{k'}^h|D_x|^{-s}D^{s_1}\nabla \psi_x\|_{L^2}\\
\lesssim&\ \|\psi_y\|_{B^{\frac12,\frac12}} \|\nabla\psi_x\|_{H^{-s,s_1}}^2\\
\leq&\ \frac{1}{64} c\|\nabla\psi_x\|_{H^{-s,s_1}}^2+C\|\psi_y\|_{B^{\frac12,\frac12}}^4 \|\nabla\psi_x\|_{H^{-s,s_1}}^2,
\end{split}
\end{equation}

\begin{equation}\notag
\begin{split}
&\lambda_j^{2s_1} \lambda_k^{-2s}\left| \int_{\mathbb R^2}\Delta_j\Delta_k^h\nabla\psi_x\mathcal B_{j,k}^9(\nabla \psi_x,\psi_y) \, dxdy\right|\\
\leq&\ \lambda_j^{2s_1} \lambda_k^{-2s} \sum_{j'\geq j- 2, k'\geq k- 2}\int_{\mathbb R^2}\left|\Delta_j\Delta_k^h\nabla\psi_x \Delta_j\Delta_k^h\left(\widetilde{\Delta}_{j'}\widetilde{\Delta}_{k'}^h\nabla \psi_x \Delta_{j'}\Delta_{k'}^h\psi_y \right) \right|\, dxdy\\
\leq &\ \lambda_j^{2s_1} \lambda_k^{-2s} \sum_{j'\geq j- 2, k'\geq k- 2} \|\Delta_j\Delta_k^h\nabla\psi_x\|_{L^\infty}\|\widetilde{\Delta}_{j'}\widetilde{\Delta}_{k'}^h\nabla \psi_x\|_{L^2}\|\Delta_{j'}\Delta_{k'}^h\psi_y\|_{L^2}\\
\lesssim&\sum_{j'\geq j- 2, k'\geq k- 2}\lambda_{j'-j}^{-s_1-\frac12} \lambda_{k'-k}^{s-\frac12}\left( \lambda_{j'}^{\frac12}\lambda_{k'}^{\frac12}\|\widetilde{\Delta}_{j'}\widetilde{\Delta}_{k'}^h\nabla \psi_x\|_{L^2}\right)\\
&\cdot \|\Delta_j\Delta_k^h|D_x|^{-s-\frac12}D^{s_1-\frac12}\nabla\psi_x\|_{L^\infty}\|\Delta_{j'}\Delta_{k'}^h|D_x|^{-s}D^{s_1}\psi_y\|_{L^2}\\
\lesssim&\ \|\nabla\psi_x\|_{B^{\frac12,\frac12}} \|\nabla\psi\|_{H^{-s,s_1}} \|\nabla\psi\|_{H^{-s,s_1+1}}.
\end{split}
\end{equation}

The proof of the lemma is complete.
\cbdu

\medskip

\begin{Lemma}\label{le-high-aux11}
Let $0<s<\frac12$ and $s_1>-\frac12$. We have 
\begin{equation}\notag
\begin{split}
&\sum_{j,k}\lambda_j^{2s_1}\lambda_k^{-2s}|\tilde I_{11}+\tilde I_{12}|\\
\leq &\ \frac{1}{16} c \|b\|_{H^{-s,s_1+1}}^2+C\left( \|b\|_{B^{\frac32,\frac12}}^2 +\|b\|_{B^{\frac12,\frac32}}^2 \right) \|\nabla\psi\|_{H^{-s,s_1}}^2\\
&+C \left(\|\psi\|_{B^{\frac32,\frac12}}^2+\|\psi\|_{B^{\frac12,\frac32}}^4 \right) \|b\|_{H^{-s,s_1+1}}^2
\end{split}
\end{equation}
for a constant $C>0$.
\end{Lemma}
\pf
Again we first use integration by parts and cancellations,
\begin{equation}\notag
\begin{split}
&\tilde I_{11}+\tilde I_{12}\\
=&\ \varepsilon_1 \int_{\mathbb R^2}\Delta_j\Delta_k^hb \Delta_j\Delta_k^h\partial_x(b_y\psi_x) \, dxdy
-\varepsilon_1 \int_{\mathbb R^2}\Delta_j\Delta_k^hb \Delta_j\Delta_k^h\partial_x(b_x \psi_y) \, dxdy\\
=&- \varepsilon_1 \int_{\mathbb R^2}\Delta_j\Delta_k^hb_x \Delta_j\Delta_k^h(b_y\psi_x) \, dxdy
+\varepsilon_1 \int_{\mathbb R^2}\Delta_j\Delta_k^hb_x \Delta_j\Delta_k^h(b_x \psi_y) \, dxdy.
\end{split}
\end{equation}
The first integral in $\tilde I_{11}+\tilde I_{12}$ can be further decomposed as:
\begin{equation}\notag
\begin{split}
&\lambda_j^{2s_1} \lambda_k^{-2s}\int_{\mathbb R^2}\Delta_j\Delta_k^hb_x \Delta_j\Delta_k^h(b_y \psi_x) \, dxdy\\
=&\sum_{m=1}^9\lambda_j^{2s_1} \lambda_k^{-2s} \int_{\mathbb R^2}\Delta_j\Delta_k^h b_x \mathcal B_{j,k}^m(b_y,  \psi_x) \, dxdy.
\end{split}
\end{equation}
We estimate each of the terms in the sum, for $0<s<\frac12$ and $s_1>-\frac12$
\begin{equation}\notag
\begin{split}
&\lambda_j^{2s_1} \lambda_k^{-2s}\left| \int_{\mathbb R^2}\Delta_j\Delta_k^h b_x \mathcal B_{j,k}^1(b_y, \psi_x) \, dxdy\right|\\
\leq&\ \lambda_j^{2s_1} \lambda_k^{-2s} \sum_{|j-j'|\leq 2, |k-k'|\leq 2}\int_{\mathbb R^2}\left|\Delta_j\Delta_k^h b_x \Delta_j\Delta_k^h\left(S_{j'-1}S_{k'-1}^hb_y\Delta_{j'}\Delta_{k'}^h\psi_x \right) \right|\, dxdy\\
\leq &\ \lambda_j^{2s_1} \lambda_k^{-2s} \sum_{|j-j'|\leq 2, |k-k'|\leq 2} \|\Delta_j\Delta_k^hb_x\|_{L^2}\|S_{j'-1}S_{k'-1}^hb_y\|_{L^\infty}\|\Delta_{j'}\Delta_{k'}^h\psi_x\|_{L^2}\\
\lesssim&\sum_{|j-j'|\leq 2, |k-k'|\leq 2} \|b_y\|_{B^{\frac12,\frac12}} \|\Delta_j\Delta_k^h|D_x|^{-s}D^{s_1}b_x\|_{L^2}\|\Delta_{j'}\Delta_{k'}^h|D_x|^{-s}D^{s_1}\psi_x\|_{L^2}\\
\lesssim&\  \|b_y\|_{B^{\frac12,\frac12}} \|\nabla b\|_{H^{-s,s_1}}\|\nabla \psi\|_{H^{-s,s_1}}\\
\leq&\  \frac{1}{256} c\|\nabla b\|_{H^{-s,s_1}}^2+C\|b_y\|_{B^{\frac12,\frac12}}^2 \|\nabla \psi\|_{H^{-s,s_1}}^2,
\end{split}
\end{equation}
\begin{equation}\notag
\begin{split}
&\lambda_j^{2s_1} \lambda_k^{-2s}\left| \int_{\mathbb R^2}\Delta_j\Delta_k^hb_x \mathcal B_{j,k}^2(b_y, \psi_x) \, dxdy\right|\\
\leq&\ \lambda_j^{2s_1} \lambda_k^{-2s} \sum_{|j-j'|\leq 2, |k-k'|\leq 2}\int_{\mathbb R^2}\left|\Delta_j\Delta_k^h b_x \Delta_j\Delta_k^h\left(\Delta_{j'}S_{k'-1}^hb_yS_{j'-1}\Delta_{k'}^h\psi_x \right) \right|\, dxdy\\
\leq &\ \lambda_j^{2s_1} \lambda_k^{-2s} \sum_{|j-j'|\leq 2, |k-k'|\leq 2} \|\Delta_j\Delta_k^hb_x\|_{L^2}\|\Delta_{j'}b_y\|_{L_x^\infty L_y^2}\|\Delta_{k'}^h\psi_x\|_{L_x^2L_y^\infty}\\
\lesssim&\sum_{|j-j'|\leq 2, |k-k'|\leq 2}  \left( \lambda_{j'}^{\frac12}\|\Delta_{j'}b_y\|_{L_x^\infty L_y^2}\right) \|\Delta_j\Delta_k^h|D_x|^{-s}D^{s_1}b_x\|_{L^2}\|\Delta_{k'}^h|D_x|^{-s}D^{s_1-\frac12}\psi_x\|_{L_x^2L_y^\infty}\\
\lesssim&\  \|b_y\|_{B^{\frac12,\frac12}} \|\nabla b\|_{H^{-s,s_1}}\|\nabla \psi\|_{H^{-s,s_1}}\\
\leq&\  \frac{1}{256} c\|\nabla b\|_{H^{-s,s_1}}^2+C\|b_y\|_{B^{\frac12,\frac12}}^2 \|\nabla \psi\|_{H^{-s,s_1}}^2,
\end{split}
\end{equation}
\begin{equation}\notag
\begin{split}
&\lambda_j^{2s_1} \lambda_k^{-2s}\left| \int_{\mathbb R^2}\Delta_j\Delta_k^hb_x \mathcal B_{j,k}^3(b_y, \psi_x) \, dxdy\right|\\
\leq&\ \lambda_j^{2s_1} \lambda_k^{-2s} \sum_{j'\geq j-2, |k-k'|\leq 2}\int_{\mathbb R^2}\left|\Delta_j\Delta_k^hb_x \Delta_j\Delta_k^h\left(\widetilde{\Delta}_{j'}S_{k'-1}^hb_y\Delta_{j'}\Delta_{k'}^h\psi_x \right) \right|\, dxdy\\
\leq &\ \lambda_j^{2s_1} \lambda_k^{-2s} \sum_{j'\geq j- 2, |k-k'|\leq 2} \|\Delta_j\Delta_k^hb_x\|_{L_x^2L_y^\infty}\|\widetilde{\Delta}_{j'}S_{k'-1}^hb_y\|_{L_x^\infty L_y^2}\|\Delta_{j'}\Delta_{k'}^h\psi_x\|_{L^2}\\
\lesssim&\sum_{j'\geq j- 2, |k-k'|\leq 2} \lambda_{j'-j}^{-s_1-\frac12}\left(\lambda_{j'}^{\frac12}\|\widetilde{\Delta}_{j'}S_{k'-1}^hb_y\|_{L_x^\infty L_y^2}\right)\\
&\cdot\|\Delta_j\Delta_k^h|D_x|^{-s}D^{s_1-\frac12}b_x\|_{L_x^2L_y^\infty}\|\Delta_{j'}\Delta_{k'}^h|D_x|^{-s}D^{s_1}\psi_x\|_{L^2}\\
\lesssim&\  \|b_y\|_{B^{\frac12,\frac12}} \|\nabla b\|_{H^{-s,s_1}}\|\nabla \psi\|_{H^{-s,s_1}}\\
\leq&\  \frac{1}{256} c\|\nabla b\|_{H^{-s,s_1}}^2+C\|b_y\|_{B^{\frac12,\frac12}}^2 \|\nabla \psi\|_{H^{-s,s_1}}^2,
\end{split}
\end{equation}
\begin{equation}\notag
\begin{split}
&\lambda_j^{2s_1} \lambda_k^{-2s}\left| \int_{\mathbb R^2}\Delta_j\Delta_k^hb_x \mathcal B_{j,k}^4(b_y, \psi_x) \, dxdy\right|\\
\leq&\ \lambda_j^{2s_1} \lambda_k^{-2s} \sum_{|j-j'|\leq 2, |k-k'|\leq 2}\int_{\mathbb R^2}\left|\Delta_j\Delta_k^hb_x \Delta_j\Delta_k^h\left(S_{j'-1}\Delta_{k'}^hb_y\Delta_{j'}S_{k'-1}^h\psi_x \right) \right|\, dxdy\\
\leq &\ \lambda_j^{2s_1} \lambda_k^{-2s}\sum_{|j-j'|\leq 2, |k-k'|\leq 2} \|\Delta_j\Delta_k^hb_x\|_{L^2}\|\Delta_{k'}^hb_y\|_{L^\infty}\|\Delta_{j'}S_{k'-1}^h\psi_x\|_{L^2}\\
\lesssim&\sum_{|j-j'|\leq 2, |k-k'|\leq 2} \|b_y\|_{B^{\frac12,\frac12}}\|\Delta_j\Delta_k^h|D_x|^{-s}D^{s_1}b_x\|_{L^2}\|\Delta_{j'}S_{k'-1}^h|D_x|^{-s}D^{s_1}\psi_x\|_{L^2}\\
\lesssim&\  \|b_y\|_{B^{\frac12,\frac12}} \|\nabla b\|_{H^{-s,s_1}}\|\nabla \psi\|_{H^{-s,s_1}}\\
\leq&\  \frac{1}{256} c\|\nabla b\|_{H^{-s,s_1}}^2+C\|b_y\|_{B^{\frac12,\frac12}}^2 \|\nabla \psi\|_{H^{-s,s_1}}^2,
\end{split}
\end{equation}
\begin{equation}\notag
\begin{split}
&\lambda_j^{2s_1} \lambda_k^{-2s}\left| \int_{\mathbb R^2}\Delta_j\Delta_k^hb_x \mathcal B_{j,k}^5(b_y, \psi_x) \, dxdy\right|\\
\leq&\ \lambda_j^{2s_1} \lambda_k^{-2s} \sum_{|j-j'|\leq 2, |k-k'|\leq 2}\int_{\mathbb R^2}\left|\Delta_j\Delta_k^hb_x \Delta_j\Delta_k^h\left(\Delta_{j'}\Delta_{k'}^hb_y S_{j'-1}S_{k'-1}^h\psi_x \right) \right|\, dxdy\\
\leq &\ \lambda_j^{2s_1} \lambda_k^{-2s} \sum_{|j-j'|\leq 2, |k-k'|\leq 2} \|\Delta_j\Delta_k^hb_x\|_{L^2}\|\Delta_{j'}\Delta_{k'}^hb_y \|_{L^2}\|S_{j'-1}S_{k'-1}^h\psi_x\|_{L^\infty}\\
\lesssim&\sum_{|j-j'|\leq 2, |k-k'|\leq 2} \|\psi_x\|_{B^{\frac12,\frac12}} \|\Delta_j\Delta_k^h |D_x|^{-s}D^{s_1}b_x\|_{L^2}\|\Delta_{j'}\Delta_{k'}^h|D_x|^{-s}D^{s_1}b_y \|_{L^2}\\
\lesssim&\  \|\psi_x\|_{B^{\frac12,\frac12}} \|\nabla b\|_{H^{-s,s_1}}^2\\
\leq&\ \frac{1}{256} c\|\nabla b\|_{H^{-s,s_1}}^2+C\|\psi_x\|_{B^{\frac12,\frac12}}^2 \|\nabla b\|_{H^{-s,s_1}}^2,
\end{split}
\end{equation}
\begin{equation}\notag
\begin{split}
&\lambda_j^{2s_1} \lambda_k^{-2s}\left| \int_{\mathbb R^2}\Delta_j\Delta_k^hb_x \mathcal B_{j,k}^6(b_y, \psi_x) \, dxdy\right|\\
\leq&\ \lambda_j^{2s_1} \lambda_k^{-2s} \sum_{j'\geq j-2, |k-k'|\leq 2}\int_{\mathbb R^2}\left|\Delta_j\Delta_k^hb_x \Delta_j\Delta_k^h\left(\widetilde{\Delta}_{j'}\Delta_{k'}^hb_y\Delta_{j'}S_{k'-1}^h\psi_x \right) \right|\, dxdy\\
\leq &\ \lambda_j^{2s_1} \lambda_k^{-2s}\sum_{j'\geq j- 2, |k-k'|\leq 2} \|\Delta_j\Delta_k^hb_x\|_{L_x^2L_y^\infty}\|\widetilde{\Delta}_{j'}\Delta_{k'}^hb_y\|_{L^2}\|\Delta_{j'}S_{k'-1}^h\psi_x\|_{L_x^\infty L_y^2}\\
\lesssim&\sum_{j'\geq j- 2, |k-k'|\leq 2}\lambda_{j'-j}^{-s_1-\frac12} \left(\lambda_{j'}^{\frac12}\|\Delta_{j'}S_{k'-1}^h\psi_x\|_{L_x^\infty L_y^2}\right)\\
&\cdot\|\Delta_j\Delta_k^h|D_x|^{-s}D^{s_1-\frac12}b_x\|_{L_x^2L_y^\infty}\|\widetilde{\Delta}_{j'}\Delta_{k'}^h|D_x|^{-s}D^{s_1}b_y\|_{L^2}\\
\lesssim&\  \|\psi_x\|_{B^{\frac12,\frac12}} \|\nabla b\|_{H^{-s,s_1}}^2\\
\leq&\  \frac{1}{256} c\|\nabla b\|_{H^{-s,s_1}}^2+C\|\psi_x\|_{B^{\frac12,\frac12}}^2 \|\nabla b\|_{H^{-s,s_1}}^2,
\end{split}
\end{equation}
\begin{equation}\notag
\begin{split}
&\lambda_j^{2s_1} \lambda_k^{-2s}\left| \int_{\mathbb R^2}\Delta_j\Delta_k^hb_x \mathcal B_{j,k}^7(b_y, \psi_x) \, dxdy\right|\\
\leq&\ \lambda_j^{2s_1} \lambda_k^{-2s} \sum_{|j-j'|\leq 2, k'\geq k- 2}\int_{\mathbb R^2}\left|\Delta_j\Delta_k^hb_x \Delta_j\Delta_k^h\left(S_{j'-1}\widetilde{\Delta}_{k'}^hb_y\Delta_{j'}\Delta_{k'}^h\psi_x \right) \right|\, dxdy\\
\leq &\ \lambda_j^{2s_1} \lambda_k^{-2s}\sum_{|j-j'|\leq 2, k'\geq k- 2} \|\Delta_j\Delta_k^hb_x\|_{L_x^\infty L_y^2}\|\widetilde{\Delta}_{k'}^hb_y\|_{L_x^2L_y^\infty}\|\Delta_{j'}\Delta_{k'}^h\psi_x\|_{L^2}\\
\lesssim&\sum_{|j-j'|\leq 2, k'\geq k- 2} \lambda_{k'-k}^{s-\frac12}\left( \lambda_{k'}^{\frac12}\|\widetilde{\Delta}_{k'}^hb_y\|_{L_x^2L_y^\infty}\right)\\
&\cdot\|\Delta_j\Delta_k^h|D_x|^{-s-\frac12}D^{s_1}b_x\|_{L_x^\infty L_y^2}\|\Delta_{j'}\Delta_{k'}^h|D_x|^{-s}D^{s_1}\psi_x\|_{L^2}\\
\lesssim&\  \|b_y\|_{B^{\frac12,\frac12}} \|\nabla b\|_{H^{-s,s_1}}\|\nabla \psi\|_{H^{-s,s_1}}\\
\leq&\  \frac{1}{256} c\|\nabla b\|_{H^{-s,s_1}}^2+C\|b_y\|_{B^{\frac12,\frac12}}^2 \|\nabla \psi\|_{H^{-s,s_1}}^2,
\end{split}
\end{equation}
\begin{equation}\notag
\begin{split}
&\lambda_j^{2s_1} \lambda_k^{-2s}\left| \int_{\mathbb R^2}\Delta_j\Delta_k^hb_x \mathcal B_{j,k}^8(b_y, \psi_x) \, dxdy\right|\\
\leq&\ \lambda_j^{2s_1} \lambda_k^{-2s} \sum_{|j-j'|\leq 2, k'\geq k- 2}\int_{\mathbb R^2}\left|\Delta_j\Delta_k^hb_x \Delta_j\Delta_k^h\left(\Delta_{j'}\widetilde{\Delta}_{k'}^hb_y S_{j'-1}\Delta_{k'}^h\psi_x \right) \right|\, dxdy\\
\leq &\ \lambda_j^{2s_1} \lambda_k^{-2s} \sum_{|j-j'|\leq 2, k'\geq k- 2} \|\Delta_j\Delta_k^hb_x\|_{L_x^\infty L_y^2}\|\Delta_{j'}\widetilde{\Delta}_{k'}^hb_y\|_{L^2}\|\Delta_{k'}^h\psi_x\|_{L_x^2L_y^\infty}\\
\lesssim&\sum_{|j-j'|\leq 2, k'\geq k- 2} \lambda_{k'-k}^{s-\frac12}\left(\lambda_{k'}^{\frac12}\|\Delta_{k'}^h\psi_x\|_{L_x^2L_y^\infty}\right)\\
&\cdot \|\Delta_j\Delta_k^h|D_x|^{-s-\frac12}D^{s_1}b_x\|_{L_x^\infty L_y^2}\|\Delta_{j'}\widetilde{\Delta}_{k'}^h|D_x|^{-s}D^{s_1}b_y\|_{L^2}\\
\lesssim&\  \|\psi_x\|_{B^{\frac12,\frac12}} \|\nabla b\|_{H^{-s,s_1}}^2\\
\leq&\  \frac{1}{256} c\|\nabla b\|_{H^{-s,s_1}}^2+C\|\psi_x\|_{B^{\frac12,\frac12}}^2 \|\nabla b\|_{H^{-s,s_1}}^2,
\end{split}
\end{equation}
\begin{equation}\notag
\begin{split}
&\lambda_j^{2s_1} \lambda_k^{-2s}\left| \int_{\mathbb R^2}\Delta_j\Delta_k^hb_x \mathcal B_{j,k}^9(b_y, \psi_x) \, dxdy\right|\\
\leq&\ \lambda_j^{2s_1} \lambda_k^{-2s} \sum_{j'\geq j- 2, k'\geq k- 2}\int_{\mathbb R^2}\left|\Delta_j\Delta_k^hb_x \Delta_j\Delta_k^h\left(\widetilde{\Delta}_{j'}\widetilde{\Delta}_{k'}^hb_y \Delta_{j'}\Delta_{k'}^h\psi_x \right) \right|\, dxdy\\
\leq &\ \lambda_j^{2s_1} \lambda_k^{-2s} \sum_{j'\geq j- 2, k'\geq k- 2} \|\Delta_j\Delta_k^hb_x\|_{L^\infty}\|\widetilde{\Delta}_{j'}\widetilde{\Delta}_{k'}^hb_y\|_{L^2}\|\Delta_{j'}\Delta_{k'}^h\psi_x\|_{L^2}\\
\lesssim&\sum_{j'\geq j- 2, k'\geq k- 2} \lambda_{j'-j}^{-s_1-\frac12}\lambda_{k'-k}^{s-\frac12}\left(\lambda_{j'}^{\frac12} \lambda_{k'}^{\frac12}\|\widetilde{\Delta}_{j'}\widetilde{\Delta}_{k'}^hb_y\|_{L^2}\right)\\
&\cdot \|\Delta_j\Delta_k^h|D_x|^{-s-\frac12}D^{s_1-\frac12}b_x\|_{L^\infty}\|\Delta_{j'}\Delta_{k'}^h|D_x|^{-s}D^{s_1}\psi_x\|_{L^2}\\
\lesssim&\ \|b_y\|_{B^{\frac12,\frac12}} \|\nabla b\|_{H^{-s,s_1}} \|\nabla\psi\|_{H^{-s,s_1}}\\
\leq&\  \frac{1}{256} c\|\nabla b\|_{H^{-s,s_1}}^2+C\|b_y\|_{B^{\frac12,\frac12}}^2 \|\nabla \psi\|_{H^{-s,s_1}}^2.
\end{split}
\end{equation}
We conclude the proof.


\cbdu

\bigskip

\section{Appendix: Estimates to obtain Lemma \ref{le-energy2}}
\label{sec-app}
We provide estimates for the integrals in (\ref{basic-energy}) in this section.  To limit the length of the paper, we only present estimates for the case $l=1$ and show details for some major terms. 

\begin{Lemma}\label{le-basic-aux12}
We have 
\begin{equation}\notag
\sum_{j,k}\lambda_k^{2}|\tilde I_{1}+\tilde I_{2}|
\leq  C \left(\|\nabla\psi\|_{H^{0,-s}}+\|\nabla^2\psi\|_{L^2} \right)^{\frac12}\|\nabla^2\partial_x\psi\|_{L^2}^{\frac12} D_1(t)
\end{equation}
for a constant $C>0$.
\end{Lemma}
\pf
Integration by parts leads to
\begin{equation}\notag
\begin{split}
\tilde I_1+\tilde I_2=&-\int_{\mathbb R^2}\Delta_j\Delta_k^hb_x \Delta_j\Delta_k^h(\psi_y\Delta \psi) \, dxdy-\int_{\mathbb R^2}\Delta_j\Delta_k^hb \Delta_j\Delta_k^h(\psi_{xy}\Delta \psi) \, dxdy\\
&+\int_{\mathbb R^2}\Delta_j\Delta_k^hb_y \Delta_j\Delta_k^h(\psi_x\Delta \psi) \, dxdy+\int_{\mathbb R^2}\Delta_j\Delta_k^hb \Delta_j\Delta_k^h(\psi_{xy}\Delta \psi) \, dxdy\\
=&-\int_{\mathbb R^2}\Delta_j\Delta_k^hb_x \Delta_j\Delta_k^h(\psi_y\Delta \psi) \, dxdy+\int_{\mathbb R^2}\Delta_j\Delta_k^hb_y \Delta_j\Delta_k^h(\psi_x\Delta \psi) \, dxdy\\
=& \int_{\mathbb R^2}\Delta_j\Delta_k^h\nabla b_x \Delta_j\Delta_k^h(\psi_y\nabla \psi) \, dxdy+\int_{\mathbb R^2}\Delta_j\Delta_k^hb_x \Delta_j\Delta_k^h(\nabla \psi_y\nabla \psi) \, dxdy\\
&-\int_{\mathbb R^2}\Delta_j\Delta_k^h\nabla b_y \Delta_j\Delta_k^h(\psi_x\nabla \psi) \, dxdy-\int_{\mathbb R^2}\Delta_j\Delta_k^hb_y \Delta_j\Delta_k^h(\nabla\psi_x\nabla \psi) \, dxdy\\
=&\int_{\mathbb R^2}\Delta_j\Delta_k^h\nabla b_x \Delta_j\Delta_k^h(\psi_y\nabla \psi) \, dxdy-\int_{\mathbb R^2}\Delta_j\Delta_k^h\nabla b_y \Delta_j\Delta_k^h(\psi_x\nabla \psi) \, dxdy.
\end{split}
\end{equation}
Using Bony's paraproduct introduced in Section \ref{sec-pre} gives
\begin{equation}\notag
\lambda_k^2 \int_{\mathbb R^2}\Delta_j\Delta_k^h\nabla b_x \Delta_j\Delta_k^h(\psi_y\nabla \psi) \, dxdy
=\sum_{m=1}^9 \lambda_k^2 \int_{\mathbb R^2}\Delta_j\Delta_k^h\nabla b_x \mathcal B_{j,k}^m(\psi_y, \nabla \psi) \, dxdy.
\end{equation}
We have the estimate
\begin{equation}\notag
\begin{split}
&\sum_{j,k}\lambda_k^2\left| \int_{\mathbb R^2}\Delta_j\Delta_k^h\nabla b_x \mathcal B_{j,k}^1(\psi_y, \nabla \psi) \, dxdy\right|\\
\leq&\sum_{j,k}\lambda_k^2 \sum_{|j-j'|\leq 2, |k-k'|\leq 2}\int_{\mathbb R^2}\left|\Delta_j\Delta_k^h\nabla b_x \Delta_j\Delta_k^h\left(S_{j'-1}S_{k'-1}^h\psi_y\Delta_{j'}\Delta_{k'}^h\nabla\psi \right) \right|\, dxdy\\
\leq & \sum_{j,k}\lambda_k^2 \sum_{|j-j'|\leq 2, |k-k'|\leq 2} \|\Delta_j\Delta_k^h\nabla b_x\|_{L^2}\|S_{j'-1}S_{k'-1}^h\psi_y\|_{L^\infty}\|\Delta_{j'}\Delta_{k'}^h\nabla\psi\|_{L^2}\\
\lesssim&\sum_{j,k}\sum_{|j-j'|\leq 2, |k-k'|\leq 2} \|\psi_y\|_{B^{\frac12,\frac12}} \|\Delta_j\Delta_k^h\nabla b_x\|_{L^2}\|\Delta_{j'}\Delta_{k'}^h\partial_x^2\nabla\psi\|_{L^2}\\
\lesssim&  \left(\|\psi_y\|_{H^{0,-s}}+\|\nabla\psi_y\|_{L^2} \right)^{\frac12}\|\nabla\partial_x\psi_y\|_{L^2}^{\frac12}\\
&\cdot \sum_{j,k}\sum_{|j-j'|\leq 2, |k-k'|\leq 2} \|\Delta_j\Delta_k^h\nabla b_x\|_{L^2}\|\Delta_{j'}\Delta_{k'}^h\partial_x^2\nabla\psi\|_{L^2}\\
\lesssim&  \left(\|\psi_y\|_{H^{0,-s}}+\|\nabla\psi_y\|_{L^2} \right)^{\frac12}\|\nabla\partial_x\psi_y\|_{L^2}^{\frac12} D_1(t).
\end{split}
\end{equation}
Similarly we have
\begin{equation}\notag
\begin{split}
&\sum_{j,k}\lambda_k^2\left| \int_{\mathbb R^2}\Delta_j\Delta_k^h\nabla b_x \mathcal B_{j,k}^2(\psi_y, \nabla \psi) \, dxdy\right|\\
\leq&\sum_{j,k}\lambda_k^2 \sum_{|j-j'|\leq 2, |k-k'|\leq 2}\int_{\mathbb R^2}\left|\Delta_j\Delta_k^h\nabla b_x \Delta_j\Delta_k^h\left(\Delta_{j'}S_{k'-1}^h\psi_yS_{j'-1}\Delta_{k'}^h\nabla\psi \right) \right|\, dxdy\\
\leq &\sum_{j,k}\lambda_k^2 \sum_{|j-j'|\leq 2, |k-k'|\leq 2} \|\Delta_j\Delta_k^h\nabla b_x\|_{L^2}\|\Delta_{j'}S_{k'-1}^h\psi_y\|_{L^\infty}\|S_{j'-1}\Delta_{k'}^h\nabla\psi\|_{L^2}\\
\lesssim&  \left(\|\psi_y\|_{H^{0,-s}}+\|\nabla\psi_y\|_{L^2} \right)^{\frac12}\|\nabla\partial_x\psi_y\|_{L^2}^{\frac12} D_1(t),
\end{split}
\end{equation}

\begin{equation}\notag
\begin{split}
&\sum_{j,k}\lambda_k^2\left| \int_{\mathbb R^2}\Delta_j\Delta_k^h\nabla b_x \mathcal B_{j,k}^3(\psi_y, \nabla \psi) \, dxdy\right|\\
\leq&\sum_{j,k} \lambda_k^2 \sum_{j'\geq j-2, |k-k'|\leq 2}\int_{\mathbb R^2}\left|\Delta_j\Delta_k^h\nabla b_x \Delta_j\Delta_k^h\left(\widetilde{\Delta}_{j'}S_{k'-1}^h\psi_y\Delta_{j'}\Delta_{k'}^h\nabla\psi \right) \right|\, dxdy\\
\leq &\sum_{j,k} \lambda_k^2 \sum_{j'\geq j- 2, |k-k'|\leq 2} \|\Delta_j\Delta_k^h\nabla b_x\|_{L^2}\|\widetilde{\Delta}_{j'}S_{k'-1}^h\psi_y\|_{L^\infty}\|\Delta_{j'}\Delta_{k'}^h\nabla\psi\|_{L^2}\\
\lesssim&  \left(\|\psi_y\|_{H^{0,-s}}+\|\nabla\psi_y\|_{L^2} \right)^{\frac12}\|\nabla\partial_x\psi_y\|_{L^2}^{\frac12} D_1(t),
\end{split}
\end{equation}

\begin{equation}\notag
\begin{split}
&\sum_{j,k}\lambda_k^2\left| \int_{\mathbb R^2}\Delta_j\Delta_k^h\nabla b_x \mathcal B_{j,k}^4(\psi_y, \nabla \psi) \, dxdy\right|\\
\leq&\sum_{j,k} \lambda_k^2 \sum_{|j-j'|\leq 2, |k-k'|\leq 2}\int_{\mathbb R^2}\left|\Delta_j\Delta_k^h\nabla b_x \Delta_j\Delta_k^h\left(S_{j'-1}\Delta_{k'}^h\psi_y\Delta_{j'}S_{k'-1}^h\nabla\psi \right) \right|\, dxdy\\
\leq &\sum_{j,k} \lambda_k^2 \sum_{|j-j'|\leq 2, |k-k'|\leq 2} \|\Delta_j\Delta_k^h\nabla b_x\|_{L^2}\|S_{j'-1}\Delta_{k'}^h\psi_y\|_{L^2}\|\Delta_{j'}S_{k'-1}^h\nabla\psi\|_{L^\infty}\\
\lesssim&\sum_{j,k}\sum_{|j-j'|\leq 2, |k-k'|\leq 2} \|\nabla \psi\|_{B^{\frac12,\frac12}} \|\Delta_j\Delta_k^h\nabla b_x\|_{L^2}\|S_{j'-1}\Delta_{k'}^h\partial_x^2\psi_y\|_{L^2}\\
\lesssim& \left(\|\nabla\psi\|_{H^{0,-s}}+\|\nabla^2\psi\|_{L^2} \right)^{\frac12}\|\nabla^2\partial_x\psi\|_{L^2}^{\frac12} D_1(t),
\end{split}
\end{equation}

\begin{equation}\notag
\begin{split}
&\sum_{j,k}\lambda_k^2\left| \int_{\mathbb R^2}\Delta_j\Delta_k^h\nabla b_x \mathcal B_{j,k}^5(\psi_y, \nabla \psi) \, dxdy\right|\\
\leq&\sum_{j,k} \lambda_k^2 \sum_{|j-j'|\leq 2, |k-k'|\leq 2}\int_{\mathbb R^2}\left|\Delta_j\Delta_k^h\nabla b_x \Delta_j\Delta_k^h\left(\Delta_{j'}\Delta_{k'}^h\psi_y S_{j'-1}S_{k'-1}^h\nabla\psi \right) \right|\, dxdy\\
\leq &\sum_{j,k} \lambda_k^2 \sum_{|j-j'|\leq 2, |k-k'|\leq 2} \|\Delta_j\Delta_k^h\nabla b_x\|_{L^2}\|\Delta_{j'}\Delta_{k'}^h\psi_y\|_{L^2}\|S_{j'-1}S_{k'-1}^h\nabla\psi\|_{L^\infty}\\
\lesssim& \left(\|\nabla\psi\|_{H^{0,-s}}+\|\nabla^2\psi\|_{L^2} \right)^{\frac12}\|\nabla^2\partial_x\psi\|_{L^2}^{\frac12} D_1(t),
\end{split}
\end{equation}

\begin{equation}\notag
\begin{split}
&\sum_{j,k}\lambda_k^2\left| \int_{\mathbb R^2}\Delta_j\Delta_k^h\nabla b_x \mathcal B_{j,k}^6(\psi_y, \nabla \psi) \, dxdy\right|\\
\leq&\sum_{j,k}\lambda_k^2 \sum_{j'\geq j-2, |k-k'|\leq 2}\int_{\mathbb R^2}\left|\Delta_j\Delta_k^h\nabla b_x \Delta_j\Delta_k^h\left(\widetilde{\Delta}_{j'}\Delta_{k'}^h\psi_y\Delta_{j'}S_{k'-1}^h\nabla\psi \right) \right|\, dxdy\\
\leq &\sum_{j,k} \lambda_k^2 \sum_{j'\geq j- 2, |k-k'|\leq 2} \|\Delta_j\Delta_k^h\nabla b_x\|_{L^2}\|\widetilde{\Delta}_{j'}\Delta_{k'}^h\psi_y\|_{L^2}\|\Delta_{j'}S_{k'-1}^h\nabla\psi\|_{L^\infty}\\
\lesssim&  \left(\|\nabla\psi\|_{H^{0,-s}}+\|\nabla^2\psi\|_{L^2} \right)^{\frac12}\|\nabla^2\partial_x\psi\|_{L^2}^{\frac12} D_1(t),
\end{split}
\end{equation}

\begin{equation}\notag
\begin{split}
&\sum_{j,k}\lambda_k^2\left| \int_{\mathbb R^2}\Delta_j\Delta_k^h\nabla b_x \mathcal B_{j,k}^7(\psi_y, \nabla \psi) \, dxdy\right|\\
\leq&\sum_{j,k} \lambda_k^2 \sum_{|j-j'|\leq 2, k'\geq k- 2}\int_{\mathbb R^2}\left|\Delta_j\Delta_k^h\nabla b_x \Delta_j\Delta_k^h\left(S_{j'-1}\widetilde{\Delta}_{k'}^h\psi_y\Delta_{j'}\Delta_{k'}^h\nabla\psi \right) \right|\, dxdy\\
\leq &\sum_{j,k} \lambda_k^2 \sum_{|j-j'|\leq 2, k'\geq k- 2} \|\Delta_j\Delta_k^h\nabla b_x\|_{L^2}\|S_{j'-1}\widetilde{\Delta}_{k'}^h\psi_y\|_{L^\infty}\|\Delta_{j'}\Delta_{k'}^h\nabla\psi\|_{L^2}\\
\lesssim&\sum_{j,k}\sum_{|j-j'|\leq 2, k'\geq k- 2} \lambda_{k-k'}^2\|\psi_y\|_{B^{\frac12,\frac12}} \|\Delta_j\Delta_k^h\nabla b_x\|_{L^2}\|\Delta_{j'}\Delta_{k'}^h\partial_x^2\nabla\psi\|_{L^2}\\
\lesssim&  \left(\|\psi_y\|_{H^{0,-s}}+\|\nabla\psi_y\|_{L^2} \right)^{\frac12}\|\nabla\partial_x\psi_y\|_{L^2}^{\frac12} D_1(t),
\end{split}
\end{equation}

\begin{equation}\notag
\begin{split}
&\sum_{j,k}\lambda_k^2\left| \int_{\mathbb R^2}\Delta_j\Delta_k^h\nabla b_x \mathcal B_{j,k}^8(\psi_y, \nabla \psi) \, dxdy\right|\\
\leq&\sum_{j,k} \lambda_k^2 \sum_{|j-j'|\leq 2, k'\geq k- 2}\int_{\mathbb R^2}\left|\Delta_j\Delta_k^h\nabla b_x \Delta_j\Delta_k^h\left(\Delta_{j'}\widetilde{\Delta}_{k'}^h\psi_y S_{j'-1}\Delta_{k'}^h\nabla\psi \right) \right|\, dxdy\\
\leq &\sum_{j,k} \lambda_k^2 \sum_{|j-j'|\leq 2, k'\geq k- 2} \|\Delta_j\Delta_k^h\nabla b_x\|_{L^2}\|\Delta_{j'}\widetilde{\Delta}_{k'}^h\psi_y\|_{L^2}\|S_{j'-1}\Delta_{k'}^h\nabla\psi\|_{L^\infty}\\
\lesssim&\sum_{j,k}\sum_{|j-j'|\leq 2, k'\geq k- 2} \lambda_{k-k'}^2\|\nabla\psi\|_{B^{\frac12,\frac12}} \|\Delta_j\Delta_k^h\nabla b_x\|_{L^2}\|\widetilde{\Delta}_{j'}\Delta_{k'}^h\partial_x^2\psi_y\|_{L^2}\\
\lesssim&   \left(\|\nabla\psi\|_{H^{0,-s}}+\|\nabla^2\psi\|_{L^2} \right)^{\frac12}\|\nabla^2\partial_x\psi\|_{L^2}^{\frac12} D_1(t),
\end{split}
\end{equation}

\begin{equation}\notag
\begin{split}
&\sum_{j,k}\lambda_k^2\left| \int_{\mathbb R^2}\Delta_j\Delta_k^h\nabla b_x \mathcal B_{j,k}^9(\psi_y, \nabla \psi) \, dxdy\right|\\
\leq&\sum_{j,k} \lambda_k^2 \sum_{j'\geq j- 2, k'\geq k- 2}\int_{\mathbb R^2}\left|\Delta_j\Delta_k^h\nabla b_x \Delta_j\Delta_k^h\left(\widetilde{\Delta}_{j'}\widetilde{\Delta}_{k'}^h\psi_y \Delta_{j'}\Delta_{k'}^h\nabla\psi \right) \right|\, dxdy\\
\leq &\sum_{j,k} \lambda_k^2 \sum_{j'\geq j- 2, k'\geq k- 2} \|\Delta_j\Delta_k^h\nabla b_x\|_{L^2}\|\widetilde{\Delta}_{j'}\widetilde{\Delta}_{k'}^h\psi_y\|_{L^\infty}\|\Delta_{j'}\Delta_{k'}^h\nabla\psi\|_{L^2}\\
\lesssim&\sum_{j,k}\sum_{j'\geq j- 2, k'\geq k- 2} \lambda_{k-k'}^2\|\psi_y\|_{B^{\frac12,\frac12}} \|\Delta_j\Delta_k^h\nabla b_x\|_{L^2}\|{\Delta}_{j'}{\Delta}_{k'}^h\partial_x^2\nabla\psi\|_{L^2}\\
\lesssim&   \left(\|\psi_y\|_{H^{0,-s}}+\|\nabla\psi_y\|_{L^2} \right)^{\frac12}\|\nabla\partial_x\psi_y\|_{L^2}^{\frac12} D_1(t).
\end{split}
\end{equation}
The second integral $\lambda_k^2 \int_{\mathbb R^2}\Delta_j\Delta_k^h\nabla b_y \Delta_j\Delta_k^h(\psi_x\nabla \psi) \, dxdy$ can be handled in an analogous way. Summarizing the estimates above leads to the conclusion of the lemma.

\cbdu

\begin{Lemma}\label{le-basic-aux34}
The estimate
\begin{equation}\notag
\begin{split}
&\sum_{j,k}\lambda_k^{2}|\tilde I_{3}+\tilde I_{4}|\\
\leq&\  C \left(\|\nabla\psi\|_{H^{0,-s}}+\|\nabla^2\psi\|_{L^2} \right)^{\frac12}\|\nabla^2\partial_x\psi\|_{L^2}^{\frac12} D_1(t)\\
&+C\left(\|\nabla^2b\|_{L^2}+\|\nabla^2\psi\|_{L^2}+\|\nabla^2\psi_x\|_{L^2}+\|\nabla^2b\|_{L^2}^{\frac14}\|\nabla\psi\|_{L^2}^{\frac34}\right) D_1(t)
\end{split}
\end{equation}
holds for a constant $C>0$.
\end{Lemma}
\pf
Again we first apply integration by parts
\begin{equation}\notag
\begin{split}
\tilde I_3+\tilde I_4=&-\int_{\mathbb R^2}\Delta_j\Delta_k^h\nabla\psi_y \Delta_j\Delta_k^h\nabla(b \psi_x) \, dxdy-\int_{\mathbb R^2}\Delta_j\Delta_k^h\nabla \psi \Delta_j\Delta_k^h\nabla(b\psi_{xy}) \, dxdy\\
&+\int_{\mathbb R^2}\Delta_j\Delta_k^h\nabla\psi_x \Delta_j\Delta_k^h\nabla(b \psi_y) \, dxdy+\int_{\mathbb R^2}\Delta_j\Delta_k^h\nabla \psi \Delta_j\Delta_k^h\nabla(b\psi_{xy}) \, dxdy\\
=&-\int_{\mathbb R^2}\Delta_j\Delta_k^h\nabla\psi_y \Delta_j\Delta_k^h\nabla(b \psi_x) \, dxdy+\int_{\mathbb R^2}\Delta_j\Delta_k^h\nabla\psi_x \Delta_j\Delta_k^h\nabla(b \psi_y) \, dxdy\\
=&-\int_{\mathbb R^2}\Delta_j\Delta_k^h\nabla\psi_y \Delta_j\Delta_k^h(\nabla b \psi_x) \, dxdy-\int_{\mathbb R^2}\Delta_j\Delta_k^h\nabla\psi_y \Delta_j\Delta_k^h(b \nabla\psi_x) \, dxdy\\
&+\int_{\mathbb R^2}\Delta_j\Delta_k^h\nabla\psi_x \Delta_j\Delta_k^h(\nabla b \psi_y) \, dxdy+\int_{\mathbb R^2}\Delta_j\Delta_k^h\nabla\psi_x \Delta_j\Delta_k^h(b \nabla\psi_y) \, dxdy.
\end{split}
\end{equation}
To continue, applying integration by parts to the first integral, commutator to the second and fourth integrals and observing cancellations among them, we obtain
\begin{equation}\notag
\begin{split}
\tilde I_3+\tilde I_4
=&\int_{\mathbb R^2}\Delta_j\Delta_k^h\psi_y \Delta_j\Delta_k^h(\nabla^2 b \psi_x) \, dxdy+\int_{\mathbb R^2}\Delta_j\Delta_k^h\psi_y \Delta_j\Delta_k^h(\nabla b \nabla\psi_x) \, dxdy\\
&-\int_{\mathbb R^2}\Delta_j\Delta_k^h\nabla\psi_y b \Delta_j\Delta_k^h\nabla\psi_x \, dxdy-\int_{\mathbb R^2}\Delta_j\Delta_k^h\nabla\psi_y [\Delta_j\Delta_k^h, b\nabla]\partial_x\psi \, dxdy\\
&+\int_{\mathbb R^2}\Delta_j\Delta_k^h\nabla\psi_x \Delta_j\Delta_k^h(\nabla b \psi_y) \, dxdy+\int_{\mathbb R^2}\Delta_j\Delta_k^h\nabla\psi_x b\Delta_j\Delta_k^h \nabla\psi_y \, dxdy\\
&+\int_{\mathbb R^2}\Delta_j\Delta_k^h\nabla\psi_x [\Delta_j\Delta_k^h,b\partial_y] \nabla\psi \, dxdy\\
=&\int_{\mathbb R^2}\Delta_j\Delta_k^h\psi_y \Delta_j\Delta_k^h(\nabla^2 b \psi_x) \, dxdy+\int_{\mathbb R^2}\Delta_j\Delta_k^h\psi_y \Delta_j\Delta_k^h(\nabla b \nabla\psi_x) \, dxdy\\
&-\int_{\mathbb R^2}\Delta_j\Delta_k^h\nabla\psi_y [\Delta_j\Delta_k^h, b\nabla]\psi_x \, dxdy
+\int_{\mathbb R^2}\Delta_j\Delta_k^h\nabla\psi_x \Delta_j\Delta_k^h(\nabla b \psi_y) \, dxdy\\
&+\int_{\mathbb R^2}\Delta_j\Delta_k^h\nabla\psi_x [\Delta_j\Delta_k^h,b\partial_y] \nabla\psi \, dxdy.
\end{split}
\end{equation}
For the third term on the right hand side, we further use integration by parts,
\begin{equation}\notag
-\int_{\mathbb R^2}\Delta_j\Delta_k^h\nabla\psi_y [\Delta_j\Delta_k^h, b\nabla]\psi_x \, dxdy=\int_{\mathbb R^2}\Delta_j\Delta_k^h\psi_y \nabla\left([\Delta_j\Delta_k^h, b\nabla]\psi_x\right) \, dxdy
\end{equation}

The first integral in $\tilde I_3+\tilde I_4$ can be further decomposed as:
\begin{equation}\notag
\lambda_k^2\int_{\mathbb R^2}\Delta_j\Delta_k^h\psi_y \Delta_j\Delta_k^h(\nabla^2b \psi_x) \, dxdy
=\sum_{m=1}^9 \lambda_k^2 \int_{\mathbb R^2}\Delta_j\Delta_k^h \psi_y \mathcal B_{j,k}^m(\nabla^2 b,  \psi_x) \, dxdy
\end{equation}
We have the estimate
\begin{equation}\notag
\begin{split}
&\sum_{j,k}\lambda_k^2\left| \int_{\mathbb R^2}\Delta_j\Delta_k^h \psi_y \mathcal B_{j,k}^1(\nabla^2 b, \psi_x) \, dxdy\right|\\
\leq &\sum_{j,k}\lambda_k^2 \sum_{|j-j'|\leq 2, |k-k'|\leq 2} \|\Delta_j\Delta_k^h \psi_y\|_{L^\infty}\|S_{j'-1}S_{k'-1}^h\nabla^2 b\|_{L^2}\|\Delta_{j'}\Delta_{k'}^h\psi_x\|_{L^2}\\
\lesssim&  \left(\|\psi_y\|_{H^{0,-s}}+\|\nabla\psi_y\|_{L^2} \right)^{\frac12}\|\nabla\partial_x\psi_y\|_{L^2}^{\frac12} D_1(t),
\end{split}
\end{equation}
\begin{equation}\notag
\begin{split}
&\sum_{j,k}\lambda_k^2\left| \int_{\mathbb R^2}\Delta_j\Delta_k^h\psi_y \mathcal B_{j,k}^2(\nabla^2b, \psi_x) \, dxdy\right|\\
\leq &\sum_{j,k} \lambda_k^2 \sum_{|j-j'|\leq 2, |k-k'|\leq 2} \|\Delta_j\Delta_k^h\psi_y\|_{L^\infty}\|\Delta_{j'}S_{k'-1}^h\nabla^2b\|_{L^2}\|S_{j'-1}\Delta_{k'}^h\psi_x\|_{L^2}\\
\lesssim&  \left(\|\psi_y\|_{H^{0,-s}}+\|\nabla\psi_y\|_{L^2} \right)^{\frac12}\|\nabla\partial_x\psi_y\|_{L^2}^{\frac12} D_1(t),
\end{split}
\end{equation}

\begin{equation}\notag
\begin{split}
&\sum_{j,k}\lambda_k^2\left| \int_{\mathbb R^2}\Delta_j\Delta_k^h\psi_y \mathcal B_{j,k}^3(\nabla^2b, \psi_x) \, dxdy\right|\\
\leq &\sum_{j,k} \lambda_k^2 \sum_{j'\geq j- 2, |k-k'|\leq 2} \|\Delta_j\Delta_k^h\psi_y\|_{L^\infty}\|\widetilde{\Delta}_{j'}S_{k'-1}^h\nabla^2b\|_{L^2}\|\Delta_{j'}\Delta_{k'}^h\psi_x\|_{L^2}\\
\lesssim&  \left(\|\psi_y\|_{H^{0,-s}}+\|\nabla\psi_y\|_{L^2} \right)^{\frac12}\|\nabla\partial_x\psi_y\|_{L^2}^{\frac12} D_1(t),
\end{split}
\end{equation}

\begin{equation}\notag
\begin{split}
&\sum_{j,k}\lambda_k^2\left| \int_{\mathbb R^2}\Delta_j\Delta_k^h\psi_y \mathcal B_{j,k}^4(\nabla^2b, \psi_x) \, dxdy\right|\\
\leq &\sum_{j,k} \lambda_k^2 \sum_{|j-j'|\leq 2, |k-k'|\leq 2} \|\Delta_j\Delta_k^h\psi_y\|_{L^2}\|S_{j'-1}\Delta_{k'}^h\nabla^2b\|_{L^2}\|\Delta_{j'}S_{k'-1}^h\psi_x\|_{L^\infty}\\
\lesssim& \left(\|\psi_x\|_{H^{0,-s}}+\|\nabla\psi_x\|_{L^2} \right)^{\frac12}\|\nabla\partial_x^2\psi\|_{L^2}^{\frac12} D_1(t),
\end{split}
\end{equation}

\begin{equation}\notag
\begin{split}
&\sum_{j,k}\lambda_k^2\left| \int_{\mathbb R^2}\Delta_j\Delta_k^h\psi_y \mathcal B_{j,k}^5(\nabla^2b, \psi_x) \, dxdy\right|\\
\leq &\sum_{j,k} \lambda_k^2 \sum_{|j-j'|\leq 2, |k-k'|\leq 2} \|\Delta_j\Delta_k^h \psi_y\|_{L^2}\|\Delta_{j'}\Delta_{k'}^h\nabla^2b \|_{L^2}\|S_{j'-1}S_{k'-1}^h\psi_x\|_{L^\infty}\\
\lesssim&  \left(\|\psi_x\|_{H^{0,-s}}+\|\nabla\psi_x\|_{L^2} \right)^{\frac12}\|\nabla\partial_x^2\psi\|_{L^2}^{\frac12} D_1(t),
\end{split}
\end{equation}

\begin{equation}\notag
\begin{split}
&\sum_{j,k}\lambda_k^2\left| \int_{\mathbb R^2}\Delta_j\Delta_k^h\psi_y \mathcal B_{j,k}^6(\nabla^2b, \psi_x) \, dxdy\right|\\
\leq &\sum_{j,k} \lambda_k^2 \sum_{j'\geq j- 2, |k-k'|\leq 2} \|\Delta_j\Delta_k^h\psi_y\|_{L^2}\|\widetilde{\Delta}_{j'}\Delta_{k'}^h\nabla^2b\|_{L^2}\|\Delta_{j'}S_{k'-1}^h\psi_x\|_{L^\infty}\\
\lesssim&  \left(\|\nabla\psi\|_{H^{0,-s}}+\|\nabla^2\psi\|_{L^2} \right)^{\frac12}\|\nabla\partial_x^2\psi\|_{L^2}^{\frac12} D_1(t),
\end{split}
\end{equation}

\begin{equation}\notag
\begin{split}
&\sum_{j,k}\lambda_k^2\left| \int_{\mathbb R^2}\Delta_j\Delta_k^h\psi_y \mathcal B_{j,k}^7(\nabla^2b, \psi_x) \, dxdy\right|\\
\leq &\sum_{j,k} \lambda_k^2 \sum_{|j-j'|\leq 2, k'\geq k- 2} \|\Delta_j\Delta_k^h\psi_y\|_{L^2}\|S_{j'-1}\widetilde{\Delta}_{k'}^h\nabla^2b\|_{L^2}\|\Delta_{j'}\Delta_{k'}^h\psi_x\|_{L^\infty}\\
\lesssim&\sum_{j,k}\sum_{|j-j'|\leq 2, k'\geq k- 2} \lambda_{k-k'}\|\psi_x\|_{B^{\frac12,\frac12}} \|\Delta_j\Delta_k^h\partial_x\psi_y\|_{L^2}\|S_{j'-1}\widetilde{\Delta}_{k'}^h\partial_x\nabla^2b\|_{L^2}\\
\lesssim&  \left(\|\psi_x\|_{H^{0,-s}}+\|\nabla\psi_x\|_{L^2} \right)^{\frac12}\|\nabla\partial_x^2\psi\|_{L^2}^{\frac12} D_1(t),
\end{split}
\end{equation}

\begin{equation}\notag
\begin{split}
&\sum_{j,k}\lambda_k^2\left| \int_{\mathbb R^2}\Delta_j\Delta_k^h\psi_y \mathcal B_{j,k}^8(\nabla^2b, \psi_x) \, dxdy\right|\\
\leq &\sum_{j,k} \lambda_k^2 \sum_{|j-j'|\leq 2, k'\geq k- 2} \|\Delta_j\Delta_k^h\psi_y\|_{L^2}\|\Delta_{j'}\widetilde{\Delta}_{k'}^h\nabla^2b\|_{L^2}\|S_{j'-1}\Delta_{k'}^h\psi_x\|_{L^\infty}\\
\lesssim&\sum_{j,k}\sum_{|j-j'|\leq 2, k'\geq k- 2} \lambda_{k-k'}\|\nabla\psi\|_{B^{\frac12,\frac12}} \|\Delta_j\Delta_k^h\partial_x\psi_y\|_{L^2}\|\Delta_{j'}\widetilde{\Delta}_{k'}^h\partial_x\nabla^2b\|_{L^2}\\
\lesssim&   \left(\|\nabla\psi\|_{H^{0,-s}}+\|\nabla^2\psi\|_{L^2} \right)^{\frac12}\|\nabla^2\partial_x\psi\|_{L^2}^{\frac12} D_1(t),
\end{split}
\end{equation}

\begin{equation}\notag
\begin{split}
&\sum_{j,k}\lambda_k^2\left| \int_{\mathbb R^2}\Delta_j\Delta_k^h\psi_y \mathcal B_{j,k}^9(\nabla^2b, \psi_x) \, dxdy\right|\\
\leq &\sum_{j,k}\lambda_k^2 \sum_{j'\geq j- 2, k'\geq k- 2} \|\Delta_j\Delta_k^h\psi_y\|_{L^2}\|\widetilde{\Delta}_{j'}\widetilde{\Delta}_{k'}^h\nabla^2b\|_{L^2}\|\Delta_{j'}\Delta_{k'}^h\psi_x\|_{L^\infty}\\
\lesssim&\sum_{j,k}\sum_{j'\geq j- 2, k'\geq k- 2} \lambda_{k-k'}\|\psi_x\|_{B^{\frac12,\frac12}}\|\Delta_j\Delta_k^h\partial_x\psi_y\|_{L^2}\|\widetilde{\Delta}_{j'}\widetilde{\Delta}_{k'}^h\partial_x\nabla^2b\|_{L^2}\\
\lesssim&   \left(\|\psi_y\|_{H^{0,-s}}+\|\nabla\psi_y\|_{L^2} \right)^{\frac12}\|\nabla\partial_x\psi_y\|_{L^2}^{\frac12} D_1(t).
\end{split}
\end{equation}
The second integral of $\tilde I_3+\tilde I_4$ can be dealt with similarly to yield
\begin{equation}\notag
\begin{split}
&\sum_{j,k}\lambda_k^2\left|\int_{\mathbb R^2}\Delta_j\Delta_k^h\psi_y \Delta_j\Delta_k^h(\nabla b \nabla\psi_x) \, dxdy\right|\\
\lesssim& \left(\|\nabla^2b\|_{L^2}+\|\nabla^2\psi_x\|_{L^2}\right) D_1(t). 
\end{split}
\end{equation}

Applying commutator estimate, we obtain for the third and fifth integrals in $\tilde I_3+\tilde I_4$:
\begin{equation}\notag
\begin{split}
&\sum_{j,k}\lambda_k^2\left|\int_{\mathbb R^2}\Delta_j\Delta_k^h\psi_y \nabla\left([\Delta_j\Delta_k^h, b]\nabla\psi_x\right) \, dxdy\right|\\
\leq&\sum_{j,k} \lambda_k^2\|\Delta_j\Delta_k^h\psi_y\|_{L^2}\|\Delta_j\Delta_k^h\nabla^2 b\|_{L^2}\|\psi_x\|_{L^\infty}\\
&+\sum_{j,k}\lambda_k^2\|\Delta_j\Delta_k^h\psi_y\|_{L^2}\|\Delta_j\Delta_k^h\nabla b\|_{L^\infty}\|\nabla\psi_x\|_{L^2}\\
\lesssim&\sum_{j,k} \lambda_k^2\|\Delta_j\Delta_k^h\psi_y\|_{L^2}\|\Delta_j\Delta_k^h\nabla^2 b\|_{L^2}\|\nabla\psi_x\|_{L^2}\\
\lesssim&\ \|\Delta_j\Delta_k^h\partial_x\psi_y\|_{L^2}\|\Delta_j\Delta_k^h\partial_x\nabla^2 b\|_{L^2}\|\nabla\psi_x\|_{L^2}\\
\lesssim&\ \|\nabla\psi_x\|_{L^2} D_1(t),
\end{split}
\end{equation}
\begin{equation}\notag
\begin{split}
&\sum_{j,k}\lambda_k^2\int_{\mathbb R^2}\left|\Delta_j\Delta_k^h\nabla\psi_x [\Delta_j\Delta_k^h,b]\partial_y \nabla\psi \right| \, dxdy\\
\leq&\sum_{j,k} \lambda_k^2 \|\Delta_j\Delta_k^h\nabla\psi_x\|_{L^2} \|\Delta_j\Delta_k^h \partial_yb\|_{L^2}\|\nabla \psi\|_{L^\infty}\\
\lesssim&\ \|\nabla^2\psi\|_{L^2} \|\Delta_j\Delta_k^h\partial_x\nabla\psi_x\|_{L^2} \|\Delta_j\Delta_k^h \partial_x\partial_yb\|_{L^2}\\
\lesssim&\ \|\nabla^2\psi\|_{L^2} D_1(t).
\end{split}
\end{equation}

We decompose the forth integral in $\tilde I_3+\tilde I_4$
\begin{equation}\notag
 \lambda_k^2\int_{\mathbb R^2}\Delta_j\Delta_k^h\nabla\psi_x \Delta_j\Delta_k^h(\nabla b \psi_y) \, dxdy=\sum_{m=1}^9 \lambda_k^2 \int_{\mathbb R^2}\Delta_j\Delta_k^h \nabla\psi_x \mathcal B_{j,k}^m(\nabla b,  \psi_y) \, dxdy
\end{equation}
It is easy to obtain
\begin{equation}\notag
\begin{split}
&\sum_{j,k}\lambda_k^2\left| \int_{\mathbb R^2}\Delta_j\Delta_k^h \nabla\psi_x \mathcal B_{j,k}^1(\nabla b, \psi_y) \, dxdy\right|\\
&+\sum_{j,k}\lambda_k^2\left| \int_{\mathbb R^2}\Delta_j\Delta_k^h\nabla\psi_x \mathcal B_{j,k}^2(\nabla b, \psi_y) \, dxdy\right|\\
&+\sum_{j,k}\lambda_k^2\left| \int_{\mathbb R^2}\Delta_j\Delta_k^h \nabla\psi_x \mathcal B_{j,k}^3(\nabla b, \psi_y) \, dxdy\right|\\
\lesssim&\ \|\nabla^2b\|_{L^2} D_1(t)
\end{split}
\end{equation}
and 
\begin{equation}\notag
\begin{split}
&\sum_{j,k}\lambda_k^2\left| \int_{\mathbb R^2}\Delta_j\Delta_k^h \nabla\psi_x \mathcal B_{j,k}^7(\nabla b, \psi_y) \, dxdy\right|\\
&+\sum_{j,k}\lambda_k^2\left| \int_{\mathbb R^2}\Delta_j\Delta_k^h\nabla\psi_x \mathcal B_{j,k}^8(\nabla b, \psi_y) \, dxdy\right|\\
&+\sum_{j,k}\lambda_k^2\left| \int_{\mathbb R^2}\Delta_j\Delta_k^h \nabla\psi_x \mathcal B_{j,k}^9(\nabla b, \psi_y) \, dxdy\right|\\
\lesssim&\ \|\nabla^2\psi\|_{L^2} D_1(t).
\end{split}
\end{equation}
The other terms are estimated as
\begin{equation}\notag
\begin{split}
&\sum_{j,k}\lambda_k^2\left| \int_{\mathbb R^2}\Delta_j\Delta_k^h \nabla\psi_x \mathcal B_{j,k}^4(\nabla b,  \psi_y) \, dxdy\right|\\
\leq &\sum_{j,k} \lambda_k^2 \sum_{|j-j'|\leq 2, |k-k'|\leq 2} \|\Delta_j\Delta_k^h \nabla\psi_x\|_{L^2}\|\Delta_{k'}^h\nabla b\|_{L_x^2L_y^\infty}\|\Delta_{j'}\psi_y\|_{L_{x}^\infty L_y^2}\\
\lesssim &\sum_{j,k} \sum_{|j-j'|\leq 2, |k-k'|\leq 2} \|\Delta_j\Delta_k^h \partial_x\nabla\psi_x\|_{L^2}\|\Delta_{k'}^h\partial_x\nabla b\|_{L^2}^{\frac12}\|\Delta_{k'}^h\partial_x\partial_y\nabla b\|_{L^2}^{\frac12}\\
&\cdot \|\Delta_{j'}\psi_y\|_{L^2}^{\frac34}  \|\Delta_{j'}\partial_x^2\psi_y\|_{L^2}^{\frac14}\\
\lesssim& \  \|\nabla^2 b\|_{L^2}^{\frac14} \|\nabla \psi\|_{L^2}^{\frac34}D_1(t),
\end{split}
\end{equation}
\begin{equation}\notag
\begin{split}
&\sum_{j,k}\lambda_k^2\left| \int_{\mathbb R^2}\Delta_j\Delta_k^h\nabla\psi_x \mathcal B_{j,k}^5(\nabla b, \psi_y) \, dxdy\right|\\
\leq &\sum_{j,k}\lambda_k^2 \sum_{|j-j'|\leq 2, |k-k'|\leq 2} \|\Delta_j\Delta_k^h \nabla\psi_x\|_{L^2}\|\Delta_{j'}\Delta_{k'}^h\nabla b \|_{L^2}\|S_{j'-1}S_{k'-1}^h\psi_y\|_{L^\infty}\\
\lesssim& \ \left(\|\psi_y\|_{H^{0,-s}}+\|\nabla\psi_y\|_{L^2} \right)^{\frac12}\|\nabla\partial_x\psi_y\|_{L^2}^{\frac12}D_1(t),
\end{split}
\end{equation}
\begin{equation}\notag
\begin{split}
&\sum_{j,k}\lambda_k^2\left| \int_{\mathbb R^2}\Delta_j\Delta_k^h\nabla\psi_x \mathcal B_{j,k}^6(\nabla b, \psi_y) \, dxdy\right|\\
\leq &\sum_{j,k} \lambda_k^2 \sum_{j'\geq j- 2, |k-k'|\leq 2} \|\Delta_j\Delta_k^h\nabla\psi_x\|_{L^2}\|\widetilde{\Delta}_{j'}\Delta_{k'}^h\nabla b\|_{L^2}\|\Delta_{j'}S_{k'-1}^h\psi_y\|_{L^\infty}\\
\lesssim&\ \left(\|\psi_y\|_{H^{0,-s}}+\|\nabla\psi_y\|_{L^2} \right)^{\frac12}\|\nabla\partial_x\psi_y\|_{L^2}^{\frac12} D_1(t).
\end{split}
\end{equation}

It completes the proof of the lemma.

\cbdu

\begin{Lemma}\label{le-basic-aux56}
The estimate
\begin{equation}\notag
\begin{split}
&\sum_{j,k}\lambda_k^{2}|\tilde I_{5}+\tilde I_{6}|\\
\leq&\  C \left(\|\nabla\psi_y\|_{H^{0,-s}}+\|\nabla^2\psi_y\|_{L^2} \right)^{\frac12}\|\nabla^2\partial_x\psi_y\|_{L^2}^{\frac12} D_1(t)\\
&+ C \left(\|\nabla\psi_x\|_{H^{0,-s}}+\|\nabla^2\psi_x\|_{L^2} \right)^{\frac12}\|\nabla^2\partial_x^2\psi\|_{L^2}^{\frac12} D_1(t)\\
&+C\left(\|\nabla^2\psi_x\|_{L^2}+\|\nabla^2\psi_y\|_{L^2}+\|\nabla\psi\|_{H^{-s,6}}+\|\nabla\psi\|_{H^{6}}\right) D_1(t)
\end{split}
\end{equation}
holds for a constant $C>0$.
\end{Lemma}
\pf
Following a few steps of integration by parts we deduce
\begin{equation}\notag
\begin{split}
\tilde I_5+\tilde I_6
=&\int_{\mathbb R^2}\Delta_j\Delta_k^h\nabla^2 b_x \Delta_j\Delta_k^h\nabla(\psi_y\nabla \psi) \, dxdy-\int_{\mathbb R^2}\Delta_j\Delta_k^h\nabla^2 b_y \Delta_j\Delta_k^h\nabla(\psi_x\nabla \psi) \, dxdy.
\end{split}
\end{equation}
We only present the estimate of the first integral in $\tilde I_5+\tilde I_6$. The second integral in $\tilde I_5+\tilde I_6$ can be analyzed in an analogous way.
Applying Bony's paraproduct decomposition gives
\begin{equation}\label{est-basic-bony56}
\begin{split}
&\lambda_k^2 \int_{\mathbb R^2}\Delta_j\Delta_k^h\nabla^2 b_x \Delta_j\Delta_k^h\nabla(\psi_y\nabla \psi) \, dxdy\\
=&\sum_{m=1}^9 \lambda_k^2 \int_{\mathbb R^2}\Delta_j\Delta_k^h\nabla^2 b_x  \mathcal B_{j,k}^m(\nabla\psi_y, \nabla \psi) \, dxdy\\
&+\sum_{m=1}^9 \lambda_k^2 \int_{\mathbb R^2}\Delta_j\Delta_k^h\nabla^2 b_x \mathcal B_{j,k}^m(\psi_y, \nabla^2 \psi) \, dxdy.
\end{split}
\end{equation}
It is trivial to have the estimates
\begin{equation}\notag
\begin{split}
&\sum_{j,k}\lambda_k^2\left| \int_{\mathbb R^2}\Delta_j\Delta_k^h\nabla^2 b_x \mathcal B_{j,k}^1(\nabla\psi_y, \nabla \psi) \, dxdy\right|\\
&+\sum_{j,k}\lambda_k^2\left| \int_{\mathbb R^2}\Delta_j\Delta_k^h\nabla^2 b_x \mathcal B_{j,k}^2(\nabla\psi_y, \nabla \psi) \, dxdy\right|\\
&+\sum_{j,k}\lambda_k^2\left| \int_{\mathbb R^2}\Delta_j\Delta_k^h\nabla^2 b_x \mathcal B_{j,k}^3(\nabla\psi_y, \nabla \psi) \, dxdy\right|\\
\lesssim&  \left(\|\nabla\psi_y\|_{H^{0,-s}}+\|\nabla^2\psi_y\|_{L^2} \right)^{\frac12}\|\nabla^2\partial_x\psi_y\|_{L^2}^{\frac12} D_1(t)
\end{split}
\end{equation}
and 
\begin{equation}\notag
\begin{split}
&\sum_{j,k}\lambda_k^2\left| \int_{\mathbb R^2}\Delta_j\Delta_k^h\nabla^2 b_x \mathcal B_{j,k}^7(\nabla\psi_y, \nabla \psi) \, dxdy\right|\\
&+\sum_{j,k}\lambda_k^2\left| \int_{\mathbb R^2}\Delta_j\Delta_k^h\nabla^2 b_x \mathcal B_{j,k}^9(\nabla\psi_y, \nabla \psi) \, dxdy\right|\\
\lesssim&\ \|\nabla^2\psi_y\|_{L^2} D_1(t).
\end{split}
\end{equation}
The other terms are estimated as
\begin{equation}\notag
\begin{split}
&\sum_{j,k}\lambda_k^2\left| \int_{\mathbb R^2}\Delta_j\Delta_k^h\nabla^2 b_x \mathcal B_{j,k}^4(\nabla\psi_y, \nabla \psi) \, dxdy\right|\\
\leq &\sum_{j,k} \lambda_k^2 \sum_{|j-j'|\leq 2, |k-k'|\leq 2} \|\Delta_j\Delta_k^h\nabla^2 b_x\|_{L^2}\|\Delta_{k'}^h\nabla\psi_y\|_{L^2_xL_y^\infty}\|\Delta_{j'}\nabla\psi\|_{L_x^\infty L_y^2}\\
\leq &\sum_{j,k} \sum_{|j-j'|\leq 2, |k-k'|\leq 2} \|\Delta_j\Delta_k^h\nabla^2 b_x\|_{L^2}\|\Delta_{k'}^h\partial_x^2\nabla\psi_y\|_{L^2_xL_y^\infty}\|\Delta_{j'}\nabla\psi\|_{L_x^\infty L_y^2}\\
\lesssim&\ \|\nabla^2 b_x\|_{L^2} \|\partial_x^2\nabla\psi\|_{L^2}^\theta \| |D_x|^{-\frac{3}{2(1-\theta)}} |D|^{\frac{3}{2(1-\theta)}}\nabla \psi\|_{L^2}^{1-\theta} \|\partial_x^2\nabla\psi\|_{L^2}^{1-\theta} \| |D_x|^{-(\frac{3}{2\theta}-2)} \nabla \psi\|_{L^2}^{\theta}\\
\lesssim& \ \|\nabla\psi\|_{H^{-\frac{3}{2(1-\theta)}, \frac{3}{2(1-\theta)}}}^{1-\theta}\|\nabla\psi\|_{H^{-\frac{3}{2\theta}+2, 0}}^{\theta} D_1(t)
\end{split}
\end{equation}
for $0<\theta<1$, and similarly
\begin{equation}\notag
\begin{split}
&\sum_{j,k}\lambda_k^2\left| \int_{\mathbb R^2}\Delta_j\Delta_k^h\nabla^2 b_x \mathcal B_{j,k}^5(\nabla\psi_y, \nabla \psi) \, dxdy\right|\\
\leq &\sum_{j,k} \sum_{|j-j'|\leq 2, |k-k'|\leq 2} \|\Delta_j\Delta_k^h\nabla^2 b_x\|_{L^2}\|\Delta_{j'}\Delta_{k'}^h\partial_x^2\nabla\psi_y\|_{L^2}\|\nabla\psi\|_{L^\infty}\\
\lesssim&\ \|\nabla^2 b_x\|_{L^2}\|\partial_x^2\nabla\psi\|_{L^2}^{\theta}\||D_x|^{-\frac{1+2\theta}{1-\theta}}D^{\frac{1}{1-\theta}}\nabla\psi\|_{L^2}^{1-\theta} \|\partial_x^2\nabla\psi\|_{L^2}^{1-\theta} \||D_x|^{-\frac{2-2\theta}{\theta}}D^{\frac{1}{2\theta}}\nabla\psi\|_{L^2}^{\theta} \\
\lesssim&\ \|\nabla\psi\|_{H^{-\frac{1+2\theta}{1-\theta}, \frac{1}{1-\theta}}}^{1-\theta} \|\nabla\psi\|_{H^{-\frac{2-2\theta}{\theta}, \frac{1}{2\theta}}}^{\theta}D_1(t),
\end{split}
\end{equation}
\begin{equation}\notag
\begin{split}
&\sum_{j,k}\lambda_k^2\left| \int_{\mathbb R^2}\Delta_j\Delta_k^h\nabla^2 b_x \mathcal B_{j,k}^6(\nabla\psi_y, \nabla \psi) \, dxdy\right|\\
\leq &\sum_{j,k} \lambda_k^2 \sum_{j'\geq j- 2, |k-k'|\leq 2} \|\Delta_j\Delta_k^h\nabla^2 b_x\|_{L^2}\|\widetilde{\Delta}_{j'}\Delta_{k'}^h\nabla\psi_y\|_{L^2}\|\Delta_{j'}S_{k'-1}^h\nabla\psi\|_{L^\infty}\\
\lesssim&\ \|\nabla\psi\|_{H^{-\frac{1+2\theta}{1-\theta}, \frac{1}{1-\theta}}}^{1-\theta} \|\nabla\psi\|_{H^{-\frac{2-2\theta}{\theta}, \frac{1}{2\theta}}}^{\theta} D_1(t),
\end{split}
\end{equation}
\begin{equation}\notag
\begin{split}
&\sum_{j,k}\lambda_k^2\left| \int_{\mathbb R^2}\Delta_j\Delta_k^h\nabla^2 b_x \mathcal B_{j,k}^8(\nabla\psi_y, \nabla \psi) \, dxdy\right|\\
\leq &\sum_{j,k} \sum_{|j-j'|\leq 2, k'\geq k- 2}  \lambda_{k-k'}^2\|\Delta_j\Delta_k^h\nabla^2 b_x\|_{L^2}\|\Delta_{j'}\widetilde{\Delta}_{k'}^h\partial_x\nabla\psi_y\|_{L^\infty_xL^2_y}\|\Delta_{k'}^h\partial_x\nabla\psi\|_{L^2_xL_y^\infty}\\
\lesssim&\ \|\nabla^2 b_x\|_{L^2}\|\partial_x^2\nabla\psi\|_{L^2}^{\theta}\||D_x|^{-\frac{4\theta-1}{2-2\theta}}D^{\frac{1}{1-\theta}}\nabla\psi\|_{L^2}^{1-\theta} \|\partial_x^2\nabla\psi\|_{L^2}^{1-\theta} \||D_x|^{-\frac{3-4\theta}{2\theta}}D^{\frac{1}{2\theta}}\nabla\psi\|_{L^2}^{\theta} \\
\lesssim&\ \|\nabla\psi\|_{H^{-\frac{4\theta-1}{2-2\theta}, \frac{1}{1-\theta}}}^{1-\theta} \|\nabla\psi\|_{H^{-\frac{3-4\theta}{2\theta}, \frac{1}{2\theta}}}^{\theta}D_1(t),
\end{split}
\end{equation}
Note the parameter $\theta\in(0,1)$ may vary in the estimates for different integrals.

For the second sum in (\ref{est-basic-bony56}), we have 
\begin{equation}\notag
\begin{split}
&\sum_{j,k}\lambda_k^2\left| \int_{\mathbb R^2}\Delta_j\Delta_k^h\nabla^2 b_x \mathcal B_{j,k}^1(\psi_y, \nabla^2 \psi) \, dxdy\right|\\
&+\sum_{j,k}\lambda_k^2\left| \int_{\mathbb R^2}\Delta_j\Delta_k^h\nabla^2 b_x \mathcal B_{j,k}^2(\psi_y, \nabla^2 \psi) \, dxdy\right|\\
&+\sum_{j,k}\lambda_k^2\left| \int_{\mathbb R^2}\Delta_j\Delta_k^h\nabla^2 b_x \mathcal B_{j,k}^3(\psi_y, \nabla^2\psi) \, dxdy\right|\\
&+\sum_{j,k}\lambda_k^2\left| \int_{\mathbb R^2}\Delta_j\Delta_k^h\nabla^2 b_x \mathcal B_{j,k}^7(\psi_y, \nabla^2\psi) \, dxdy\right|\\
&+\sum_{j,k}\lambda_k^2\left| \int_{\mathbb R^2}\Delta_j\Delta_k^h\nabla^2 b_x \mathcal B_{j,k}^9(\psi_y, \nabla^2\psi) \, dxdy\right|\\
\lesssim&  \ \|\nabla\psi\|_{H^{-4, 1}}^{\frac12}\|\nabla \psi\|_{H^4}^{\frac12} D_1(t).
\end{split}
\end{equation}
Other terms in the sum are handled in the following, for $0<\theta<1$ (which may be taken different values for different terms),
\begin{equation}\notag
\begin{split}
&\sum_{j,k}\lambda_k^2\left| \int_{\mathbb R^2}\Delta_j\Delta_k^h\nabla^2 b_x \mathcal B_{j,k}^4(\psi_y, \nabla^2 \psi) \, dxdy\right|\\
\leq &\sum_{j,k} \sum_{|j-j'|\leq 2, |k-k'|\leq 2} \|\Delta_j\Delta_k^h\nabla^2 b_x\|_{L^2}\|\Delta_{k'}^h\partial_x^2\psi_y\|_{L^2_xL_y^\infty}\|\Delta_{j'}\nabla^2\psi\|_{L_x^\infty L_y^2}\\
\lesssim&\ \|\nabla^2 b_x\|_{L^2} \|\partial_x^2\nabla\psi\|_{L^2}^\theta \| |D_x|^{-\frac{4\theta-1}{2(1-\theta)}} |D|^{\frac{1}{2(1-\theta)}}\nabla \psi\|_{L^2}^{1-\theta} \|\partial_x^2\nabla\psi\|_{L^2}^{1-\theta} \| |D_x|^{-(\frac{3}{2\theta}-2)}D^{\frac{1}{\theta}} \nabla \psi\|_{L^2}^{\theta}\\
\lesssim& \ \|\nabla\psi\|_{H^{-\frac{4\theta-1}{2(1-\theta)}, \frac{1}{2(1-\theta)}}}^{1-\theta}\|\nabla\psi\|_{H^{-\frac{3}{2\theta}+2, \frac{1}{\theta}}}^{\theta} D_1(t),
\end{split}
\end{equation} 
\begin{equation}\notag
\begin{split}
&\sum_{j,k}\lambda_k^2\left| \int_{\mathbb R^2}\Delta_j\Delta_k^h\nabla^2 b_x \mathcal B_{j,k}^5(\psi_y, \nabla^2 \psi) \, dxdy\right|\\
\lesssim &\sum_{j,k} \sum_{|j-j'|\leq 2, |k-k'|\leq 2} \|\Delta_j\Delta_k^h\nabla^2 b_x\|_{L^2}\|\Delta_{j'}\Delta_{k'}^h\partial_x^2\psi_y\|_{L^2}\|\nabla^2\psi\|_{L^\infty}\\
\lesssim&\ \|\nabla^2 b_x\|_{L^2}\|\partial_x^2\nabla\psi\|_{L^2}^{\theta}\||D_x|^{-\frac{2\theta-1}{1-\theta}}\nabla\psi\|_{L^2}^{1-\theta} \|\partial_x^2\nabla\psi\|_{L^2}^{1-\theta} \||D_x|^{-\frac{2-2\theta}{\theta}}D^{\frac{3}{2\theta}}\nabla\psi\|_{L^2}^{\theta} \\
\lesssim&\ \|\nabla\psi\|_{H^{-\frac{2\theta-1}{1-\theta},0}}^{1-\theta} \|\nabla\psi\|_{H^{-\frac{2-2\theta}{\theta}, \frac{3}{2\theta}}}^{\theta}D_1(t),
\end{split}
\end{equation}
\begin{equation}\notag
\begin{split}
&\sum_{j,k}\lambda_k^2\left| \int_{\mathbb R^2}\Delta_j\Delta_k^h\nabla^2 b_x \mathcal B_{j,k}^6(\psi_y, \nabla^2 \psi) \, dxdy\right|\\
\leq &\sum_{j,k}\lambda_k^2 \sum_{j'\geq j- 2, |k-k'|\leq 2} \|\Delta_j\Delta_k^h\nabla^2 b_x\|_{L^2}\|\widetilde{\Delta}_{j'}\Delta_{k'}^h\psi_y\|_{L^2}\|\Delta_{j'}S_{k'-1}^h\nabla^2\psi\|_{L^\infty}\\
\lesssim&\ \|\nabla\psi\|_{H^{-\frac{2\theta-1}{1-\theta},0}}^{1-\theta} \|\nabla\psi\|_{H^{-\frac{2-2\theta}{\theta}, \frac{3}{2\theta}}}^{\theta}D_1(t),
\end{split}
\end{equation}
\begin{equation}\notag
\begin{split}
&\sum_{j,k}\lambda_k^2\left| \int_{\mathbb R^2}\Delta_j\Delta_k^h\nabla^2 b_x \mathcal B_{j,k}^8(\psi_y, \nabla^2 \psi) \, dxdy\right|\\
\leq &\sum_{j,k} \sum_{|j-j'|\leq 2, k'\geq k- 2}  \lambda_{k-k'}^2\|\Delta_j\Delta_k^h\nabla^2 b_x\|_{L^2}\|\Delta_{j'}\widetilde{\Delta}_{k'}^h\partial_x\psi_y\|_{L^\infty_xL^2_y}\|\Delta_{k'}^h\partial_x\nabla^2\psi\|_{L^2_xL_y^\infty}\\
\lesssim&\ \|\nabla^2 b_x\|_{L^2}\|\partial_x^2\nabla\psi\|_{L^2}^{\theta}\||D_x|^{-\frac{4\theta-1}{2-2\theta}}\nabla\psi\|_{L^2}^{1-\theta} \|\partial_x^2\nabla\psi\|_{L^2}^{1-\theta} \||D_x|^{-\frac{3-4\theta}{2\theta}}D^{\frac{3}{2\theta}}\nabla\psi\|_{L^2}^{\theta} \\
\lesssim&\ \|\nabla\psi\|_{H^{-\frac{4\theta-1}{2-2\theta}, 0}}^{1-\theta} \|\nabla\psi\|_{H^{-\frac{3-4\theta}{2\theta}, \frac{3}{2\theta}}}^{\theta}D_1(t).
\end{split}
\end{equation}
The lemma follows from taking appropriate values of $\theta$ in the estimates above.
\cbdu

\begin{Lemma}\label{le-basic-aux78}
The estimate
\begin{equation}\notag
\begin{split}
&\sum_{j,k}\lambda_k^{2}|\tilde I_{7}+\tilde I_{8}|\\
\leq&\ C\left(\|\nabla^3\psi_x\|_{L^2}+\|\nabla^3\psi_y\|_{L^2}\right.\\
&\left.+\|\nabla\psi\|_{H^{-s,-s}}+\|\nabla\psi\|_{H^{-s,6}}+\|\nabla\psi\|_{H^{6}}\right) D_1(t)
\end{split}
\end{equation}
holds for a constant $C>0$.
\end{Lemma}
\pf
Again using integration by parts, observing cancellations and applying commutators, we obtain
\begin{equation}\notag
\begin{split}
&\tilde I_7+\tilde I_8\\
=&\int_{\mathbb R^2}\Delta_j\Delta_k^h\nabla^2\psi_y \Delta_j\Delta_k^h(\nabla^2 b \psi_x) \, dxdy-\int_{\mathbb R^2}\Delta_j\Delta_k^h\nabla^2\psi_x \Delta_j\Delta_k^h(\nabla^2 b \psi_y) \, dxdy\\
&-2\int_{\mathbb R^2}\Delta_j\Delta_k^h\nabla^2\psi_y \Delta_j\Delta_k^h(\nabla b_x \nabla\psi) \, dxdy
+2\int_{\mathbb R^2}\Delta_j\Delta_k^h\nabla^2\psi_x \Delta_j\Delta_k^h(\nabla b_y \nabla \psi) \, dxdy\\
&-\int_{\mathbb R^2}\Delta_j\Delta_k^h\nabla^2\psi_y \partial_x [\Delta_j\Delta_k^h,b]\nabla \nabla\psi \, dxdy
+\int_{\mathbb R^2}\Delta_j\Delta_k^h\nabla^2\psi_x \partial_y[\Delta_j\Delta_k^h,b]\nabla \nabla\psi \, dxdy.
\end{split}
\end{equation}
The first and second integrals in $\tilde I_7+\tilde I_8$ can be analyzed analogously, while the third and forth integrals can be estimated similarly as well. Thus we only show details for the first and forth integrals. 
The first integral in $\tilde I_7+\tilde I_8$ can be further decomposed as:
\begin{equation}\notag
\lambda_k^2\int_{\mathbb R^2}\Delta_j\Delta_k^h\nabla^2\psi_y \Delta_j\Delta_k^h(\nabla^2b \psi_x) \, dxdy
=\sum_{m=1}^9 \lambda_k^2 \int_{\mathbb R^2}\Delta_j\Delta_k^h \nabla^2\psi_y \mathcal B_{j,k}^m(\nabla^2 b,  \psi_x) \, dxdy.
\end{equation}
Among the nine terms in the sum, we have
\begin{equation}\notag
\begin{split}
&\sum_{j,k}\lambda_k^2\left| \int_{\mathbb R^2}\Delta_j\Delta_k^h \nabla^2\psi_y \mathcal B_{j,k}^1(\nabla^2 b, \psi_x) \, dxdy\right|\\
&\sum_{j,k}\lambda_k^2\left| \int_{\mathbb R^2}\Delta_j\Delta_k^h\nabla^2\psi_y \mathcal B_{j,k}^2(\nabla^2b, \psi_x) \, dxdy\right|\\
&\sum_{j,k}\lambda_k^2\left| \int_{\mathbb R^2}\Delta_j\Delta_k^h\nabla^2\psi_y \mathcal B_{j,k}^3(\nabla^2b, \psi_x) \, dxdy\right|\\
\lesssim&\  \|\nabla^3\psi_y\|_{L^2} D_1(t).
\end{split}
\end{equation}
The other ones have the estimate, for some $0<\theta<1$
\begin{equation}\notag
\begin{split}
&\sum_{j,k}\lambda_k^2\left| \int_{\mathbb R^2}\Delta_j\Delta_k^h\nabla^2\psi_y \mathcal B_{j,k}^4(\nabla^2b, \psi_x) \, dxdy\right|\\
\lesssim&\sum_{j,k}\sum_{|j-j'|\leq 2, |k-k'|\leq 2} \|\psi_x\|_{L^{\infty}} \|\Delta_j\Delta_k^h\partial_x\nabla^2\psi_y\|_{L^2}\|S_{j'-1}\Delta_{k'}^h\partial_x\nabla^2b\|_{L^2}\\
\lesssim&\ \|\partial_x^2\nabla\psi\|_{L^2}^{\theta}\||D_x|^{-\frac{2\theta-1}{1-\theta}}D^{-\frac{1}{2-2\theta}}\nabla\psi\|_{L^2}^{1-\theta}
 \|\partial_x^2\nabla\psi\|_{L^2}^{1-\theta}\||D_x|^{-\frac{2-2\theta}{\theta}}D^{\frac{2}{\theta}}\nabla\psi\|_{L^2}^{\theta}\|\partial_x\nabla^2b\|_{L^2}\\
\lesssim& \ \|\nabla\psi\|_{H^{-\frac{2\theta-1}{1-\theta}, -\frac{1}{2-2\theta}}}^{1-\theta} \|\nabla\psi\|_{H^{-\frac{2-2\theta}{\theta}, \frac{2}{\theta}}}^{\theta}D_1(t),
\end{split}
\end{equation}
\begin{equation}\notag
\begin{split}
&\sum_{j,k}\lambda_k^2\left| \int_{\mathbb R^2}\Delta_j\Delta_k^h\nabla^2\psi_y \mathcal B_{j,k}^5(\nabla^2b, \psi_x) \, dxdy\right|\\
\lesssim&\sum_{j,k}\sum_{|j-j'|\leq 2, |k-k'|\leq 2} \|\psi_x\|_{L^{\infty}} \|\Delta_j\Delta_k^h \partial_x\nabla^2\psi_y\|_{L^2}\|\Delta_{j'}\Delta_{k'}^h\partial_x\nabla^2b \|_{L^2}\\
\lesssim& \ \|\nabla\psi\|_{H^{-\frac{2\theta-1}{1-\theta}, -\frac{1}{2-2\theta}}}^{1-\theta} \|\nabla\psi\|_{H^{-\frac{2-2\theta}{\theta}, \frac{2}{\theta}}}^{\theta}D_1(t),
\end{split}
\end{equation}
\begin{equation}\notag
\begin{split}
&\sum_{j,k}\lambda_k^2\left| \int_{\mathbb R^2}\Delta_j\Delta_k^h\nabla^2\psi_y \mathcal B_{j,k}^6(\nabla^2b, \psi_x) \, dxdy\right|\\
\leq &\sum_{j,k} \lambda_k^2 \sum_{j'\geq j- 2, |k-k'|\leq 2} \|\Delta_j\Delta_k^h\nabla^2\psi_y\|_{L^2}\|\widetilde{\Delta}_{j'}\Delta_{k'}^h\nabla^2b\|_{L^2}\|\Delta_{j'}S_{k'-1}^h\psi_x\|_{L^\infty}\\
\lesssim&\sum_{j,k}\sum_{j'\geq j- 2, |k-k'|\leq 2} \|\psi_x\|_{L^{\infty}} \|\Delta_j\Delta_k^h\partial_x\nabla^2\psi_y\|_{L^2}\|\widetilde{\Delta}_{j'}\Delta_{k'}^h\partial_x\nabla^2b\|_{L^2}\\
\lesssim& \ \|\nabla\psi\|_{H^{-\frac{2\theta-1}{1-\theta}, -\frac{1}{2-2\theta}}}^{1-\theta} \|\nabla\psi\|_{H^{-\frac{2-2\theta}{\theta}, \frac{2}{\theta}}}^{\theta}D_1(t),
\end{split}
\end{equation}
\begin{equation}\notag
\begin{split}
&\sum_{j,k}\lambda_k^2\left| \int_{\mathbb R^2}\Delta_j\Delta_k^h\nabla^2\psi_y \mathcal B_{j,k}^7(\nabla^2b, \psi_x) \, dxdy\right|\\
\lesssim&\sum_{j,k}\sum_{|j-j'|\leq 2, k'\geq k- 2} \lambda_{k-k'}\|\psi_x\|_{L^{\infty}} \|\Delta_j\Delta_k^h\partial_x\nabla^2\psi_y\|_{L^2}\|S_{j'-1}\widetilde{\Delta}_{k'}^h\partial_x\nabla^2b\|_{L^2}\\
\lesssim& \ \|\nabla\psi\|_{H^{-\frac{2\theta-1}{1-\theta}, -\frac{1}{2-2\theta}}}^{1-\theta} \|\nabla\psi\|_{H^{-\frac{2-2\theta}{\theta}, \frac{2}{\theta}}}^{\theta}D_1(t),
\end{split}
\end{equation}

\begin{equation}\notag
\begin{split}
&\sum_{j,k}\lambda_k^2\left| \int_{\mathbb R^2}\Delta_j\Delta_k^h\nabla^2\psi_y \mathcal B_{j,k}^8(\nabla^2b, \psi_x) \, dxdy\right|\\
\lesssim&\sum_{j,k}\sum_{|j-j'|\leq 2, k'\geq k- 2} \lambda_{k-k'}\|\psi_x\|_{L^{\infty}} \|\Delta_j\Delta_k^h\partial_x\nabla^2\psi_y\|_{L^2}\|\Delta_{j'}\widetilde{\Delta}_{k'}^h\partial_x\nabla^2b\|_{L^2}\\
\lesssim& \ \|\nabla\psi\|_{H^{-\frac{2\theta-1}{1-\theta}, -\frac{1}{2-2\theta}}}^{1-\theta} \|\nabla\psi\|_{H^{-\frac{2-2\theta}{\theta}, \frac{2}{\theta}}}^{\theta}D_1(t),
\end{split}
\end{equation}

\begin{equation}\notag
\begin{split}
&\sum_{j,k}\lambda_k^2\left| \int_{\mathbb R^2}\Delta_j\Delta_k^h\nabla^2\psi_y \mathcal B_{j,k}^9(\nabla^2b, \psi_x) \, dxdy\right|\\
\lesssim&\sum_{j,k}\sum_{j'\geq j- 2, k'\geq k- 2} \lambda_{k-k'}\|\psi_x\|_{L^{\infty}}\|\Delta_j\Delta_k^h\partial_x\nabla^2\psi_y\|_{L^2}\|\widetilde{\Delta}_{j'}\widetilde{\Delta}_{k'}^h\partial_x\nabla^2b\|_{L^2}\\
\lesssim& \ \|\nabla\psi\|_{H^{-\frac{2\theta-1}{1-\theta}, -\frac{1}{2-2\theta}}}^{1-\theta} \|\nabla\psi\|_{H^{-\frac{2-2\theta}{\theta}, \frac{2}{\theta}}}^{\theta}D_1(t).
\end{split}
\end{equation}

The forth integral in $\tilde I_7+\tilde I_8$ can be further decomposed as:
\begin{equation}\notag
\lambda_k^2\int_{\mathbb R^2}\Delta_j\Delta_k^h\nabla^2\psi_x \Delta_j\Delta_k^h(\nabla b_y \nabla\psi) \, dxdy
=\sum_{m=1}^9 \lambda_k^2 \int_{\mathbb R^2}\Delta_j\Delta_k^h \nabla^2\psi_x \mathcal B_{j,k}^m(\nabla b_y,  \nabla\psi) \, dxdy.
\end{equation}
Again the first three terms in the sum are simple,
\begin{equation}\notag
\begin{split}
&\sum_{j,k}\lambda_k^2\left| \int_{\mathbb R^2}\Delta_j\Delta_k^h \nabla^2\psi_x \mathcal B_{j,k}^1(\nabla b_y, \nabla\psi) \, dxdy\right|\\
&+\sum_{j,k}\lambda_k^2\left| \int_{\mathbb R^2}\Delta_j\Delta_k^h \nabla^2\psi_x \mathcal B_{j,k}^1(\nabla b_y, \nabla\psi) \, dxdy\right|\\
&+\sum_{j,k}\lambda_k^2\left| \int_{\mathbb R^2}\Delta_j\Delta_k^h \nabla^2\psi_x \mathcal B_{j,k}^1(\nabla b_y, \nabla\psi) \, dxdy\right|\\
\lesssim&\  \|\nabla^3\psi_x\|_{L^2} D_1(t).
\end{split}
\end{equation}
The other terms in the sum are estimated as
\begin{equation}\notag
\begin{split}
&\sum_{j,k}\lambda_k^2\left| \int_{\mathbb R^2}\Delta_j\Delta_k^h\nabla^2\psi_x \mathcal B_{j,k}^4(\nabla b_y, \nabla\psi) \, dxdy\right|\\
\lesssim&\sum_{j,k}\sum_{|j-j'|\leq 2, |k-k'|\leq 2} \|\nabla\psi\|_{L^{\infty}} \|\Delta_j\Delta_k^h\partial_x\nabla^2\psi_x\|_{L^2}\|S_{j'-1}\Delta_{k'}^h\partial_x\nabla b_y\|_{L^2}\\
\lesssim&\ \|\partial_x^2\nabla\psi\|_{L^2}^{\theta}\||D_x|^{-\frac{2\theta}{1-\theta}}D^{\frac{1}{2-2\theta}}\nabla\psi\|_{L^2}^{1-\theta}
 \|\partial_x^2\nabla\psi\|_{L^2}^{1-\theta}\|D_x^{2}D^{\frac{1}{\theta}}\nabla\psi\|_{L^2}^{\theta}\|\partial_x\nabla b_y\|_{L^2}\\
\lesssim& \ \|\nabla\psi\|_{H^{-\frac{2\theta}{1-\theta}, \frac{1}{2-2\theta}}}^{1-\theta} \|\nabla\psi\|_{H^{2, \frac{1}{\theta}}}^{\theta}D_1(t),
\end{split}
\end{equation}
\begin{equation}\notag
\begin{split}
&\sum_{j,k}\lambda_k^2\left| \int_{\mathbb R^2}\Delta_j\Delta_k^h\nabla^2\psi_x \mathcal B_{j,k}^5(\nabla b_y, \nabla\psi) \, dxdy\right|\\
\lesssim&\sum_{j,k}\sum_{|j-j'|\leq 2, |k-k'|\leq 2} \|\nabla\psi\|_{L^{\infty}} \|\Delta_j\Delta_k^h \partial_x\nabla^2\psi_x\|_{L^2}\|\Delta_{j'}\Delta_{k'}^h\partial_x\nabla b_y \|_{L^2}\\
\lesssim& \ \|\nabla\psi\|_{H^{-\frac{2\theta}{1-\theta}, \frac{1}{2-2\theta}}}^{1-\theta} \|\nabla\psi\|_{H^{2, \frac{1}{\theta}}}^{\theta}D_1(t),
\end{split}
\end{equation}
\begin{equation}\notag
\begin{split}
&\sum_{j,k}\lambda_k^2\left| \int_{\mathbb R^2}\Delta_j\Delta_k^h\nabla^2\psi_x \mathcal B_{j,k}^6(\nabla b_y, \nabla\psi) \, dxdy\right|\\
\leq &\sum_{j,k} \lambda_k^2 \sum_{j'\geq j- 2, |k-k'|\leq 2} \|\Delta_j\Delta_k^h\nabla^2\psi_x\|_{L^2}\|\widetilde{\Delta}_{j'}\Delta_{k'}^h\nabla b_y\|_{L^2}\|\Delta_{j'}S_{k'-1}^h\nabla\psi\|_{L^\infty}\\
\lesssim&\sum_{j,k}\sum_{j'\geq j- 2, |k-k'|\leq 2} \|\nabla\psi\|_{L^{\infty}} \|\Delta_j\Delta_k^h\partial_x\nabla^2\psi_x\|_{L^2}\|\widetilde{\Delta}_{j'}\Delta_{k'}^h\partial_x\nabla b_y\|_{L^2}\\
\lesssim& \ \|\nabla\psi\|_{H^{-\frac{2\theta}{1-\theta}, \frac{1}{2-2\theta}}}^{1-\theta} \|\nabla\psi\|_{H^{2, \frac{1}{\theta}}}^{\theta}D_1(t),
\end{split}
\end{equation}
\begin{equation}\notag
\begin{split}
&\sum_{j,k}\lambda_k^2\left| \int_{\mathbb R^2}\Delta_j\Delta_k^h\nabla^2\psi_x \mathcal B_{j,k}^7(\nabla b_y, \nabla\psi) \, dxdy\right|\\
\lesssim&\sum_{j,k}\sum_{|j-j'|\leq 2, k'\geq k- 2} \lambda_{k-k'}\|\nabla\psi\|_{L^{\infty}} \|\Delta_j\Delta_k^h\partial_x\nabla^2\psi_x\|_{L^2}\|S_{j'-1}\widetilde{\Delta}_{k'}^h\partial_x\nabla b_y\|_{L^2}\\
\lesssim& \ \|\nabla\psi\|_{H^{-\frac{2\theta}{1-\theta}, \frac{1}{2-2\theta}}}^{1-\theta} \|\nabla\psi\|_{H^{2, \frac{1}{\theta}}}^{\theta}D_1(t),
\end{split}
\end{equation}
\begin{equation}\notag
\begin{split}
&\sum_{j,k}\lambda_k^2\left| \int_{\mathbb R^2}\Delta_j\Delta_k^h\nabla^2\psi_x \mathcal B_{j,k}^8(\nabla b_y, \nabla\psi) \, dxdy\right|\\
\lesssim&\sum_{j,k}\sum_{|j-j'|\leq 2, k'\geq k- 2} \lambda_{k-k'}\|\nabla\psi\|_{L^{\infty}} \|\Delta_j\Delta_k^h\partial_x\nabla^2\psi_x\|_{L^2}\|\Delta_{j'}\widetilde{\Delta}_{k'}^h\partial_x\nabla b_y\|_{L^2}\\
\lesssim& \ \|\nabla\psi\|_{H^{-\frac{2\theta}{1-\theta}, \frac{1}{2-2\theta}}}^{1-\theta} \|\nabla\psi\|_{H^{2, \frac{1}{\theta}}}^{\theta}D_1(t),
\end{split}
\end{equation}
\begin{equation}\notag
\begin{split}
&\sum_{j,k}\lambda_k^2\left| \int_{\mathbb R^2}\Delta_j\Delta_k^h\nabla^2\psi_x \mathcal B_{j,k}^9(\nabla b_y, \nabla\psi) \, dxdy\right|\\
\lesssim&\sum_{j,k}\sum_{j'\geq j- 2, k'\geq k- 2} \lambda_{k-k'}\|\nabla\psi\|_{L^{\infty}}\|\Delta_j\Delta_k^h\partial_x\nabla^2\psi_x\|_{L^2}\|\widetilde{\Delta}_{j'}\widetilde{\Delta}_{k'}^h\partial_x\nabla b_y\|_{L^2}\\
\lesssim& \ \|\nabla\psi\|_{H^{-\frac{2\theta}{1-\theta}, \frac{1}{2-2\theta}}}^{1-\theta} \|\nabla\psi\|_{H^{2, \frac{1}{\theta}}}^{\theta}D_1(t).
\end{split}
\end{equation}

Applying commutator estimate we deduce for the fifth and sixth integrals in $\tilde I_7+\tilde I_8$, for  $0<\theta<1$ (again may be different in different terms)
\begin{equation}\notag
\begin{split}
&\sum_{j,k}\lambda_k^2\int_{\mathbb R^2}\left|\Delta_j\Delta_k^h\nabla^2\psi_y \partial_x [\Delta_j\Delta_k^h,b]\nabla \nabla\psi \right| \, dxdy\\
\leq&\sum_{j,k} \lambda_k^2 \|\Delta_j\Delta_k^h\nabla^2\psi_y\|_{L^2}\|\Delta_j\Delta_k^h\partial_x\nabla b\|_{L^2}\|\nabla \psi\|_{L^\infty}\\
\lesssim&\sum_{j,k} \|\Delta_j\Delta_k^h\partial_x\nabla^2\psi_y\|_{L^2}\|\Delta_j\Delta_k^h\partial_x^2\nabla b\|_{L^2}\|\nabla \psi\|_{L^\infty}\\
\lesssim& \ \|\nabla\psi\|_{H^{-\frac{2\theta}{1-\theta}, \frac{1}{2-2\theta}}}^{1-\theta} \|\nabla\psi\|_{H^{-\frac{2-2\theta}{\theta}, \frac{2}{\theta}}}^{\theta}D_1(t),
\end{split}
\end{equation}
\begin{equation}\notag
\begin{split}
&\sum_{j,k}\lambda_k^2\int_{\mathbb R^2}\left|\Delta_j\Delta_k^h\nabla^2\psi_x \partial_y [\Delta_j\Delta_k^h,b]\nabla \nabla\psi \right| \, dxdy\\
\leq&\sum_{j,k} \lambda_k^2 \|\Delta_j\Delta_k^h\nabla^2\psi_x\|_{L^2}\|\Delta_j\Delta_k^h\partial_y\nabla b\|_{L^2}\|\nabla \psi\|_{L^\infty}\\
\lesssim&\sum_{j,k} \|\Delta_j\Delta_k^h\partial_x\nabla^2\psi_x\|_{L^2}\|\Delta_j\Delta_k^h\partial_x\nabla b_y\|_{L^2}\|\nabla \psi\|_{L^\infty}\\
\lesssim& \ \|\nabla\psi\|_{H^{-\frac{2\theta}{1-\theta}, \frac{1}{2-2\theta}}}^{1-\theta} \|\nabla\psi\|_{H^{2, \frac{1}{\theta}}}^{\theta}D_1(t).
\end{split}
\end{equation}

Summarizing the estimates above and taking appropriate values of $\theta\in(0,1)$ concludes the proof of the lemma.

\cbdu

\begin{Lemma}\label{le-basic-aux1112}
The estimates
\begin{equation}\notag
\begin{split}
&\sum_{j,k}\lambda_k^{2}|\tilde I_{9}+\tilde I_{10}|\\
\leq&\ C\left(\|\nabla \psi\|_{H^{-s,s}}+\|\nabla \psi\|_{H^{-s,6}}+\|\nabla^2\psi_x\|_{L^2}\right.\\
&\left.+\|\nabla^2\psi_y\|_{L^2}+\|\nabla\psi_y\|_{L^2}\right) D_1(t)\\
&\sum_{j,k}\lambda_k^{2}|\tilde I_{11}+\tilde I_{12}|\\
\leq&\ C\left(\|\nabla b_y\|_{L^2}+\|\nabla b_x\|_{L^2}+\|\nabla\psi_x\|_{L^2}+\|\nabla\psi_y\|_{L^2}\right) D_1(t)
\end{split}
\end{equation}
hold for a constant $C>0$.
\end{Lemma}
\pf
With integration by parts and cancellations, we have
\begin{equation}\notag
\begin{split}
&\tilde I_9+\tilde I_{10}\\
=&\ \varepsilon_1 \int_{\mathbb R^2}\Delta_j\Delta_k^h\psi_x \Delta_j\Delta_k^h(\psi_y \Delta\psi_x) \, dxdy
-\varepsilon_1 \int_{\mathbb R^2}\Delta_j\Delta_k^h\psi_x \Delta_j\Delta_k^h(\psi_x \Delta\psi_y) \, dxdy\\
=&-\varepsilon_1 \int_{\mathbb R^2}\Delta_j\Delta_k^h\nabla\psi_x \Delta_j\Delta_k^h(\psi_y \nabla\psi_x) \, dxdy
+\varepsilon_1 \int_{\mathbb R^2}\Delta_j\Delta_k^h\nabla\psi_x \Delta_j\Delta_k^h (\psi_x\nabla\psi_y) \, dxdy,
\end{split}
\end{equation}
\begin{equation}\notag
\begin{split}
&\tilde I_{11}+\tilde I_{12}\\
=&\ \varepsilon_1 \int_{\mathbb R^2}\Delta_j\Delta_k^hb \Delta_j\Delta_k^h\partial_x(b_y\psi_x) \, dxdy
-\varepsilon_1 \int_{\mathbb R^2}\Delta_j\Delta_k^hb \Delta_j\Delta_k^h\partial_x(b_x \psi_y) \, dxdy\\
=&- \varepsilon_1 \int_{\mathbb R^2}\Delta_j\Delta_k^hb_x \Delta_j\Delta_k^h(b_y\psi_x) \, dxdy
+\varepsilon_1 \int_{\mathbb R^2}\Delta_j\Delta_k^hb_x \Delta_j\Delta_k^h(b_x \psi_y) \, dxdy.
\end{split}
\end{equation}
Note these integrals involve lower order terms compared to previous terms, like the ones in $\tilde I_{7}+\tilde I_{8}$. The estimates are easier and hence details are omitted.
\cbdu

\bigskip

\section*{Acknowledgement}
The author is grateful for the hospitality of the Institute for Advanced Study and Princeton University where the work was completed.

\end{document}